\documentclass[12pt,english]{article}
\usepackage{amsfonts,babel,amssymb,amsmath,amsthm}
\usepackage{graphicx}
\usepackage{framed}
\usepackage{float}

\numberwithin{equation}{section}

\newcommand{\ds}{\displaystyle}
\newcommand{\smallfrac}[2]{{\textstyle\frac{#1}{#2}}}
\newcommand{\jump}[1]{[\![#1]\!]}
\newcommand{\ave}[1]{\{\!\!\{#1\}\!\!\}}

\title{Convolution Quadrature for Wave Simulations}

\author{Matthew Hassell \& Francisco--Javier Sayas \\
Department of Mathematical Sciences, University of Delaware\\
{\tt \{mhassell,fjsayas\}@math.udel.edu}
}

\date{\today}

\begin{document}

\maketitle


The following document contains the notes prepared for a course to be delivered by the second author at the {\em XVI Jacques-Louis Lions Spanish-French School on Numerical Simulation in Physics and Engineering}, in Pamplona (Spain), September 2014. We will not spend much time with the introduction. Let it be said that this is a course on how to approximate causal convolutions and convolution equations, that is, expressions 
\[
\int_0^t f(t-\tau) g(\tau)\mathrm d\tau=h(t),
\]
where either $g$ or $h$ is unknown. This seems like a very small problem to be working on when it is presented in this flippant form. The truth is hidden in what the convolution integral means. We will be dealing with operator valued distributions $f$ in convolution with a function valued distribution $g$. A large set of problems related to the scattering of linear waves (acoustic, elastic, electromagnetic) can be written in this form after being moved to the boundary of the scatterer. We will focus on acoustic waves (we will actually restrict all our attention to a single model problem) and on a particular class of discretization methods, the so-called Convolution Quadrature method, introduced by Christian Lubich in the mid-eighties. 

These notes will emphasize the introduction of concepts and algorithms in a rigorous language, while avoiding proofs. As a matter of fact, we will not state a single theorem explicitly. (This is not due to us not liking theorems, but with the goal of keeping a more narrative tone.) However, clear results will be stated as part of the text. 

The mathematics of the field of time domain boundary integral equations (which include many important and highly non-trivial examples of convolution equations) involve many interesting and deep analytic concepts, as the reader will be able to ascertain from these notes. The deeper mathematical structure of this field is explored in the lecture notes \cite{Sayas:2013}, in a step-by-step traditional mathematical fashion, with no much time for computation. From that point of view, these notes represent the algorithmic counterpart to \cite{Sayas:2013}. Convolution Quadrature is not the only method to approximate convolution equations that appear in wave propagation phenomena. Galerkin and collocation methods compete with CQ in interest, applicability, and good properties. We will not discuss or compare methods, especially because much is still to be explored both in theory and practice. We will not comment on existing literature on CQ for scattering problems on elastic and electromagnetic waves either.

Before we start, let us take some time for acknowledgements. Our research is partially funded by the National Science Foundation (grant DMS 1216356). We now become I. I (FJS) want to thank the organizers of the EHF2014 for the invitation to participate in the school. It is actually my second time in this series (the first one was in Laredo, so many years ago). Since then, the school has made its name even longer by honoring the extraordinary Jacques-Louis Lions, who happens to be my academic great-grandfather. Much of what I know on time domain integral equations and CQ has been a consequence of readings and discussions with Lehel Banjai and Christian Lubich. Both of them are an inspiration for practitioners of serious numerical analysis. Matthew Hassell and I have been working on these notes for several months, trying to give a wide perspective of the mathematical and computational aspects of the problem. We hope the readers will enjoy them as much as we did writing them.

\tableofcontents

\section{Causal convolutions and Laplace transforms}

In this section we introduce what will be the scope of this course: convolution of causal functions and distributions. Causal convolution operators will be often recognized through their Laplace transforms (which will be called their transfer functions). As a first step towards a precise determination of the kind of functions we will be dealing with, let us define the term {\bf causal}. A function $f:\mathbb R \to X$ (where $X$ is any vector space) is said to be an $X$-valued causal function when $f(t)=0$ for all $t<0$. The reader might already wonder what the point is to have functions defined in $\mathbb R$ when all we want from them is their restriction to $[0,\infty)$. However, this will be a fundamental distinction in our way of presenting wave propagation operators, where the vanishing past of the function will be often used, and where a non-zero value at time $t=0$ will be considered a jump discontinuity. The attentive reader will have already seen that the independent variable $t$ will often be called time.

We are going to spend this entire section giving a very general justification to causal convolutions and their Laplace transforms:
\begin{framed}
\begin{eqnarray*}
(f*g)(t) &=&\int_0^t f(t-\tau)g(\tau)\mathrm d\tau,\\
\mathrm F(s) &=& \int_0^\infty e^{-st} f(t)\mathrm dt.
\end{eqnarray*}
\end{framed}

\subsection{Causal functions and convolutions}\label{sec:1.1}

The {\bf causal convolution} of two causal functions $f:\mathbb R\to \mathbb R$ and $g:\mathbb R\to \mathbb R$ is defined as
\begin{equation}\label{eq:1.1}
(f*g)(t):=\int_0^t f(t-\tau)g(\tau)\mathrm d\tau.
\end{equation}
This definition makes sense, for instance, if both functions are continuous. It also makes sense if one of the functions is continuous and the other one is integrable. (It can be extended to many other cases, but we will wait for this.) Note that this definition coincides with the more traditional form of the convolution of functions
\[
(f*g)(t)=\int_{-\infty}^\infty f(t-\tau)g(\tau)\mathrm d\tau,
\]
when $f$ and $g$ are causal. The first extension we will need to consider is when $f:\mathbb R\to \mathbb R^{n\times m}$ and $g:\mathbb R\to \mathbb R^m$. In this case, \eqref{eq:1.1} defines a causal function $f*g:\mathbb R\to \mathbb R^n$. In this more general definition (where the convolution integrals are easily defined component by component), it is clear that we cannot even discuss commutativity of the convolution operator \eqref{eq:1.1}. The second big generalization involves two Hilbert spaces $X$ and $Y$ and the space
\[
\mathcal B(X,Y):=\{ A:X\to Y\,:\, \mbox{$A$ linear and bounded}\}.
\]
We can then start with a causal continuous function $f:\mathbb R\to \mathcal B(X,Y)$ and a causal function $g:\mathbb R\to X$ and obtain through convolution \eqref{eq:1.1} a causal function $f*g:\mathbb R\to Y$. Because all functions involved have been assumed to be continuous, the integration in \eqref{eq:1.1} can be easily understood to be defined in the sense of a Riemann integral for each value of $t$. Finally, a {\bf causal convolution equation} is an equation of the form
\begin{equation}\label{eq:1.2}
(f*g)(t)=h(t) \qquad \forall t,
\end{equation}
where $h:\mathbb R\to Y$ is causal, $f:\mathbb R\to \mathcal B(X,Y)$ is causal, and we look for a causal $X$-valued function $g$.

\paragraph{Some preliminary examples.} The simplest possible example of convolution is the causal antiderivative
\[
\int_0^t g(\tau)\mathrm d \tau,
\]
corresponding to the convolution with the Heaviside function:
\[
H(t):=\left\{ \begin{array}{ll} 1, & t\ge 0,\\ 0, & t< 0.\end{array}\right.
\]
A slightly more general operator is given by the expression
\begin{equation}\label{eq:1.3}
\int_0^t e^{\lambda (t-\tau)} g(\tau)\mathrm d\tau.
\end{equation}
Note that if we define
\begin{equation}\label{eq:1.4}
y(t)=\int_0^t e^{\lambda (t-\tau)} g(\tau)\mathrm d\tau,
\end{equation}
then $y$ is the only causal solution to the equation
\[
\dot y-\lambda\,y=g,
\]
or, in the more traditional language of ordinary differential equations, $y$ satisfies
\begin{equation}\label{eq:1.5}
\dot y-\lambda\,y=g, \qquad \mbox{in $[0,\infty)$}, \qquad y(0)=0,
\end{equation}
and has been extended by zero to the negative real axis. The formula \eqref{eq:1.4} is the variation of parameters formula (or Duhamel principle) for the initial value problem \eqref{eq:1.5}. Equation \eqref{eq:1.5} shows our first use of the dot as the symbol for time-differentiation. Similarly
\[
y(t)=\lambda^{-1}\int_0^t \sin(\lambda(t-\tau))g(\tau)\mathrm d\tau
\]
is the operator that yields the unique causal solution to
\[
\ddot y+\lambda^2 y=g.
\]
The latter example can be easily extended to cover some sort of discrete wave equations. We start with a discrete version of the second derivative in one (implicitly given) space dimension:
\[
\mathbb R^{N\times N}\ni \Delta_N:=\frac1{(N+1)^2} 
\left[ \begin{array}{ccccc}
-2 & 1 \\ 
1 & -2 & 1 \\
&\ddots & \ddots & \ddots \\
& & 1 & -2 & 1 \\
& & & 1 & -2
\end{array}\right].
\]
Since $\Delta_N$ is symmetric and it is not terribly complicated to see that it is negative definite, we can find a unitary matrix $D_N$ such that
\[
\Delta_N =D_N \left[\begin{array}{cccc} -\lambda_1^2 \\ & -\lambda_2^2 \\
& & \ddots \\ & & & -\lambda_N^2\end{array}\right] D_N^\top, \qquad \lambda_j > 0 \quad\forall j.
\]
We can then define
\[
(-\Delta_N)^{1/2} =D_N \left[\begin{array}{cccc} \lambda_1 \\ & \lambda_2 \\
& & \ddots \\ & & & \lambda_N\end{array}\right] D_N^\top,
\]
and
\[
\sin (t(-\Delta_N)^{1/2})=D_N \left[\begin{array}{cccc} 
	\sin(\lambda_1 t)\\ 
	& \sin(\lambda_2 t) \\
	& & \ddots \\ & & & \sin(\lambda_N t)
\end{array}\right] D_N^\top,
\]
and finally the vector-valued operator
\[
\mathbf y(t):=(-\Delta_N)^{-1/2} \int_0^t \sin((t-\tau)(-\Delta_N)^{1/2}) \mathbf g(\tau)\mathrm d \tau,
\]
which yields the solution of
\[
\ddot{\mathbf y}-\Delta_N \mathbf y=\mathbf g \quad\mbox{in $[0,\infty)$}, \qquad \mathbf y(0)=\dot{\mathbf y}(0)=\mathbf 0,
\]
extended by zero to negative $t$. This is a discrete version (after using low order finite differences) of the solution operator to the one-dimensional wave equation
\begin{alignat*}{6}
\partial_t^2 y -\partial_x^2 y = g && \qquad & \mbox{in $[0,\infty)\times [0,1]$},\\  
y(\cdot,0)=\partial_t y(\cdot,0)\equiv 0 &&& \mbox{in $[0,1]$},\\
y(0,\cdot)=y(1,\cdot)=0 & && \mbox{in $[0,\infty)$}.
\end{alignat*}

\paragraph{A convolution equation.}
One example of causal convolution equation is the Abel integral equation
\begin{equation}\label{eq:1.6}
\int_0^t \frac{g(\tau)}{\sqrt{t-\tau}}\mathrm d\tau=h(t) \qquad t\ge 0.
\end{equation}
The weakly singular operator in the left-hand-side of \eqref{eq:1.6} is actually related to the anti-differentiation operator. If we define
\[
y(t)=\int_0^t \left(\int_0^r \frac{g(\tau)}{\sqrt{r-\tau}\sqrt{t-r}}\mathrm d\tau
\right)\mathrm d r,
\]
then it is easy to see (it requires some patience and the use of the Euler Beta function) that
\[
y(t)=\pi \int_0^t g(\tau)\mathrm d\tau,
\]
that is, the Abel integral operator can be consider as a square root of the antiderivative and, therefore,
\begin{equation}\label{eq:1.7}
\frac1{\sqrt\pi}\int_0^t \frac{\dot g(\tau)}{\sqrt{t-\tau}}\mathrm d\tau
\end{equation}
is a square root of the differentiation operator. This is one of the Caputo fractional derivatives. Note that at this time we cannot yet understand the operator in \eqref{eq:1.7} as a convolution operator.

\paragraph{A much more complicated example.} Just to give a better flavor of operators to come, let us show one related to propagation of linear waves in the plane. Let $X=Y=L^2(\Gamma)$, where $\Gamma$ is a simple closed curve in the plane. Consider now the operator
\[
(\mathcal V(t) \xi)(\mathbf x):=\frac1{2\pi}\int_\Gamma 
\frac{H(t-|\mathbf x-\mathbf y|)}{\sqrt{t^2-|\mathbf x-\mathbf y|^2}}
\xi(\mathbf y)\mathrm d\Gamma(\mathbf y),
\]
where $H$ is the Heaviside function.  
For given $t$, this is a well defined operator $X\to Y$. We can then define the convolution of $\mathcal V$ with a function $g:\mathbb R\to L^2(\Gamma)$ (which can better be understood as a function $g(\mathbf y,t)$ with $\mathbf y\in \Gamma$ such that $g(\cdot,t)\equiv 0$ for $t< 0$), leading to the expression
\[
\frac1{2\pi} \int_\Gamma \left(\int_0^{t-|\mathbf x-\mathbf y|}
\frac{g(\mathbf y,\tau)}{\sqrt{(t-\tau)^2-|\mathbf x-\mathbf y|^2}}\mathrm d\tau\right)\mathrm d\Gamma(\mathbf y).
\]

\subsection{Causal distributions}

It is customary to keep on extending the definitions of the previous subsection to `functions' $f$ that include Dirac deltas or their derivatives. These extensions can be carried out with better or worse justified limiting processes. Instead of doing that, we will go the whole nine yards and deal with some {\em very elementary concepts of function- and operator-valued distributions.}

\begin{framed}
\noindent
{\bf On notation.} We will always keep in mind some Hilbert spaces $X,Y,Z,...$ and the spaces of bounded linear operators between pairs of them $\mathcal B(X,Y),..$. When we want to refer at the same time to the Hilbert spaces or to the spaces of operators, we will just refer to a general Banach space $\mathbb X$.
\end{framed}

\paragraph{The test space.} We consider the space of smooth compactly supported functions
\[
\mathcal D(\mathbb R):=\{ \psi\in \mathcal C^\infty(\mathbb R)\,:\, \psi\equiv 0 \mbox{ outside $[-M,M]$ for some $M$}\}.
\]
We will not need a precise definition of support, but here it is just in case
\[
\mathrm{supp}\, \psi=\mbox{closure of }\{ t\,:\,\psi(t)\neq 0\}.
\]
A sequence $\{\psi_n\}\subset\mathcal D(\mathbb R)$ is said to converge to $\psi\in \mathcal D(\mathbb R)$ when the support of all the elements of the sequence and of the limiting function is contained in a bounded interval $[-M,M]$, and
\[
\mbox{for all $m\ge 0$}, \quad
\psi_n^{(m)} \to \psi^{(m)} \, \mbox{ uniformly in $\mathbb R$}.
\]

\paragraph{Causal vector-valued distributions.} A distribution with values in a Banach space $\mathbb X$ is a functional $f:\mathcal D(\mathbb R) \to \mathbb X$ that is sequentially continuous, that is, such that it transforms convergent sequences to convergent sequences. An $\mathbb X$-valued distribution $f$ such that
\[
\langle f,\psi\rangle=0\qquad \mbox{whenever supp $\psi\subset (-\infty,0)$}
\]
is called a causal distribution. Note how we have used the angled bracket for the action of the distribution $f$ on the test function $\psi$. The simplest example of causal distributions are causal functions. If $f:\mathbb R\to \mathbb X$ is a continuous function, then we can define
\[
\langle f,\psi\rangle:=\int_0^\infty f(t)\psi(t)\mathrm d t,
\]
with the integral defined in the sense of a Riemann integral, or, even simpler, as a limit
\[
\int_0^\infty f(t)\psi(t)\mathrm d t=\lim_{N\to\infty}\frac1{N}
\sum_{n=0}^\infty f(n/N) \psi(n/N).
\]
Then $f$ defines a causal distribution. If $x\in \mathbb X$, we can define
\[
\langle x\otimes \delta_0,\psi\rangle:=\psi(0)\,x,
\]
which is a causal $\mathbb X$-valued distribution. More generally, if $t_0\ge 0$, then
\[
\langle x\otimes \delta_{t_0},\psi\rangle:=\psi(t_0)x
\]
is a causal distribution. If we take $t_0<0$, we still have a distribution, but it is not causal anymore. 

\paragraph{Steady state operators.}
Once we have a causal distribution $f$ with values in the space $X$, any bounded linear operator $A:X\to Y$ allows us to define the distribution
\[
\langle Af,\psi\rangle:=A\langle f,\psi\rangle,
\]
with values in $Y$. When we make a linear operator $A:X\to Y$ act on $X$-valued causal distributions, we will say that we have used a {\bf steady state} operator.
In particular, if $X\subset Y$ with bounded inclusion, then every $X$-valued distribution can be understood as a $Y$-valued distribution.

\paragraph{Differentiation.} It is very simple to see why if $\psi_n$ converges to $\psi$ in $\mathcal D(\mathbb R)$, then $\dot\psi_n$ converges to $\dot\psi$ in $\mathcal D(\mathbb R)$. Therefore, if $f$ is a causal $X$-valued distribution, then
\[
\langle \dot f,\psi\rangle:=-\langle f,\dot\psi\rangle
\]
is a causal $X$-valued distribution too. For instance, consider $x\in X$ and let $H$ be the Heaviside function. The derivative of
\[
\langle x\otimes H,\psi\rangle=\left(\int_0^\infty \psi(t)\mathrm dt\right)x
\]
is
\[
-\langle x\otimes H,\dot\psi\rangle=\psi(0)\,x=\langle x\otimes \delta_0,\psi\rangle.
\]
The derivative of $x\otimes \delta_0$ is
\[
\langle x\otimes \dot\delta_0,\psi\rangle=-\dot\psi(0)x.
\]

\subsection{Laplace transforms}\label{sec:1.3}

\paragraph{Why?} Let us first go back to Section \ref{sec:1.1}. Let us admit that all the following formal manipulations actually make sense:
\begin{alignat*}{6}
\int_0^\infty e^{-s t} (f*g)(t)\mathrm d t 
	& = \int_0^\infty e^{-st} \left(\int_0^t f(t-\tau)g(\tau)\mathrm d\tau\right)\mathrm d t\\
	& = \int_0^\infty \left(\int_\tau^\infty e^{-s t} f(t-\tau) g(\tau)\mathrm d t\right)\mathrm d \tau\\
	& = \int_0^\infty e^{-s\tau}\left(\int_\tau^\infty e^{-s(t-\tau)} f(t-\tau)\mathrm d t\right) g(\tau)\mathrm d \tau\\
	& =\left(\int_0^\infty e^{-st} f(t)\mathrm d t\right)\left(\int_0^\infty e^{-st} g(t)\mathrm d t\right).
\end{alignat*}
This computation shows how the {\bf Laplace transform}
\[
h\mapsto \int_0^\infty e^{-s t} h(t)\mathrm d t
\]
maps convolutions to `products'. Let us try to first give a precise meaning to the Laplace transform.

\paragraph{The Schwartz class.} The expression
\[
\mathcal S(\mathbb R):=\{ \psi\in \mathcal C^\infty(\mathbb R)\,:\, p_m\,\smallfrac{\mathrm d}{\mathrm dt^m}\psi\in L^\infty(\mathbb R)
\quad\forall m\ge 0\},\qquad p_m(t):=1+t^{2m}, 
\]
gives an abbreviated definition of the set of all smooth functions whose derivatives of all orders can be bounded by non-vanishing rational functions of all possible decays at infinity. This set obviously contains $\mathcal D(\mathbb R)$. Convergence in $\mathcal S(\mathbb R)$ is defined as follows: the sequence $\{\psi_n\}\subset \mathcal S(\mathbb R)$ converges to $\psi\in \mathcal S(\mathbb R)$, when
\[
\mbox{for all $m\ge 0$}\qquad
p_m\,\psi_n^{(m)}\to p_m\,\psi^{(m)}\,
\mbox{uniformly in $\mathbb R$}.
\]
Using some simple cut-off arguments, it is possible to show that every element of $\mathcal S(\mathbb R)$ is the limit of a sequence of elements of $\mathcal D(\mathbb R)$, that is, $\mathcal D(\mathbb R)$ is dense in the Schwartz class.

\paragraph{A careful construction of the Laplace transform.} We start with a smooth version of the Heaviside function $H$. Consider a function $h:\mathbb R\to \mathbb R$ with the following properties:
\begin{equation}\label{eq:1.8}
h\in \mathcal C^\infty(\mathbb R), \qquad 0\le h\le 1, \qquad h\equiv 1 \mbox{ in $[-\smallfrac12,\infty)$}, \qquad h\equiv 0 \mbox{ in $(-\infty,-1]$}.
\end{equation}
This function can be easily constructed using the antiderivative of a positive smooth compactly supported function, conveniently displaced and scaled.
Let now
\[
s\in \mathbb C_+:=\{ s\in \mathbb C\,:\, \mathrm{Re}\,s>0\}.
\]
Then, the function
\begin{equation}\label{eq:1.9}
\psi_s(t):= h(t) e^{-s\,t}
\end{equation}
is an element of the Schwartz class. Therefore, there exists a sequence $\{\psi_{n,s}\}\subset \mathcal D(\mathbb R)$ converging to $\psi_s$ in $\mathcal S(\mathbb R)$.
We will say that the causal $\mathbb X$-valued distribution $f$ has a Laplace transform when
\[
\lim_{n\to \infty}\langle f,\psi_{n,s}\rangle \mbox{ exists (in $\mathbb X$) for all $s\in \mathbb C_+$.}
\]
We then define the function $\mathrm F:\mathbb C_+\to \mathbb X$ given by the limit
\begin{equation}\label{eq:1.10}
\mathrm F(s):=\lim_{n\to \infty}\langle f,\psi_{n,s}\rangle.
\end{equation}
It is common to use all the following expressions:
\[
\mathrm F(s)=\mathcal L\{f\}(s)=\langle f,\psi_s\rangle=\langle f,e^{-s\cdot}\rangle.
\]
In the last one we act on the understanding that because $f$ is causal, the function $h$ is actually invisible.
Actually,  if we fix $s$ but define any other $h$ with the same properties, the difference between the corresponding $\psi_s$ is an element of $\mathcal D(\mathbb R)$ with support in $[-1,-\smallfrac12]$ and thus completely invisible for every causal distribution. In other words, the definition \eqref{eq:1.10} is independent of the particular choice of $h$, as long as this smooth Heaviside function satisfies \eqref{eq:1.8} or a similar expression with $[-1,-\smallfrac12]$ substituted by any closed interval contained in the negative real axis.

\paragraph{Some examples.} The Laplace transform of
\[
x\otimes \delta_{t_0}
\]
is
\[
\langle x\otimes\delta_{t_0},e^{-s \cdot}\rangle=\langle \delta_{t_0}, e^{-s\cdot}\rangle\, x=
e^{-s t_0} x.
\]
The Laplace transform of
\[
x\otimes H
\]
is
\[
\langle x\otimes H,e^{-s\cdot}\rangle = \left(\int_0^\infty e^{-st}\mathrm dt\right) x=\frac1{s} x.
\]
If the Laplace transform of an $X$-valued distribution $f$ exists and $A:X\to Y$ is a steady-state operator (bounded linear), then the Laplace transform of $Af$ is
\[
\mathcal L\{ A f\} (s)=\langle A f,e^{-s\cdot}\rangle = A \langle f,e^{-s\cdot}\rangle =A\mathrm F(s).
\]

\paragraph{Laplace transform and differentiation.} A simple computation shows that $\dot\psi_s=-s\psi_s+\varphi_s$, where $\varphi_s\in \mathcal D(\mathbb R)$ is supported in $(-\infty,0)$. Therefore, if $f$ has a Laplace transform, then $\dot f$ has a Laplace transform and
\begin{alignat*}{6}
\mathcal L\{ \dot f\} 
	&=\langle \dot f,\psi_s\rangle =-\langle f,-s\psi_s+\varphi_s\rangle\\
	&=\langle f,s\psi_s\rangle = s\mathrm F(s),
\end{alignat*}
since $\langle f,\varphi_s\rangle=0$, due to the causality of $f$. This is the differentiation theorem. 

\paragraph{A remark.} Readers used to classical Laplace transforms might be missing the value of $f$ at $t=0$ subtracted from the right-hand-side. It is missing in our formulation, because the derivative $\dot f$ includes also a Dirac delta at $t=0$ caused by a non-zero $f(0)$. More precisely: assume that $f$ is a fastly decaying smooth function $f:[0,\infty)\to \mathbb X$, extended by zero to the negative real axis. Let $f':[0,\infty)\to \mathbb X$ be its classical derivative, which we also assume to be fastly decaying. Then
\[
\dot f=f(0)\otimes \delta_0+f',
\]
that is, the distributional derivative contains a functional part $f'$ and a Dirac delta term related to the initial value of $f$. Then
\[
\mathcal L\{ f'\} (s)=\mathcal L\{ \dot f\}(s)-\mathcal L\{ f(0)\otimes \delta_0\}= s\mathrm F(s)-f(0),
\]
which is the formula we typically learn in an introduction to ODE class. The lesson to learn here is to be exactly aware of what we understand by the derivative. In our case, the derivative will always be defined for a causal function, and will include anything happening at $t=0$.

\paragraph{Warning.} The reader might think that we are developing quite a general theory. This is however not the case. The Laplace transform can sometimes be defined in a half plane of the complex plane $\mathbb C$, larger or smaller than $\mathbb C_+$. It is sometimes possible to extend it to the entire complex plane or to a larger region of the complex plane. The kind of operators that we want to focus on (arising from hyperbolic evolutionary equations) does not require to go any further than what we have done here. Parabolic equations produce transforms that are defined in the complementary set of a sector in the negative complex half plane. The theory of Laplace transforms is often defined for non-causal distributions, leading to double-sided Laplace transforms, which end up being analytic extensions of Fourier transforms. As already mentioned, we will try to keep the theory quite close to our interests.

\paragraph{Analiticity.} One of the nicest surprises of advanced analysis is the fact that many results in the theory of analytic functions of a complex variable can be extended almost word by word to analytic functions of a complex variable with values in a Banach space $\mathbb X$. In particular, it is quite easy to prove that if $f$ is causal and Laplace transformable, with values in $\mathbb X$, then the function 
\[
\mathbb C_+\ni s\longmapsto \mathrm F(s)\in \mathbb X
\]
is differentiable in $s$ and therefore it is an analytic $\mathbb X$-valued function. 

\paragraph{Real or complex spaces?} There is some fine detail missing in our treatment of the definition of the Laplace transform. We have implicitly assumed that $\mathcal D(\mathbb R)$ and $\mathcal S(\mathbb R)$ are spaces of real-valued functions. In this case, $\psi_s= h \exp(-s\cdot)$ is complex-valued, and we should be careful in defining
\[
\mathrm F(s)=
\langle f,\psi_s\rangle:=\langle f,\mathrm{Re}\,\psi_s\rangle+\imath \langle f,\mathrm{Im} \psi_s\rangle.
\]
For the kind of applications we have in mind $\mathbb X$ is a real Banach space and $\mathrm F(s)$ takes values in the complexification $\mathbb X+\imath \mathbb X$. (There is no problem when $\mathbb X$ is a complex space though.) A reader who is not comfortable with the idea of complexification should think of
\[
L^2(\Omega) :=\{ f:\Omega \to \mathbb R\,:\, \int_\Omega |f|^2 <\infty\}
\]
becoming
\[
\{ f:\Omega \to \mathbb C\,:\, \int_\Omega |f|^2 <\infty\}.
\]
Here is an important consequence that we will exploit in Section \ref{sec:3}: if $\mathbb X$ is a real space, then its complexification is a complex vector space that admits a conjugation operator; in this case
\begin{equation}\label{eq:1.300}
\mathrm F(\overline s)=\overline{\mathrm F(s)}.
\end{equation}

\subsection{Transfer functions and convolution operators}\label{sec:1.4}

\paragraph{A very weak definition of convolution.} Let $f$ be a $\mathcal B(X,Y)$-valued causal Laplace transformable distribution and let $g$ be an $X$-valued causal Laplace transformable distribution. The convolution $f*g$ is defined to be the $Y$-valued causal distribution whose Laplace transform satisfies
\begin{equation}\label{eq:1.11}
\mathcal L\{ f* g\}(s)=\mathrm F(s)\mathrm G(s).
\end{equation}
The right-hand-side of \eqref{eq:1.11} is the action of the operator $\mathrm F(s)\in \mathcal B(X,Y)$ on $\mathrm G(s)\in X$. The justification of how this formula makes sense follows from a theorem for the inversion of the Laplace transform given at the end of this section. For the moment being, we will accept this as a definition. A similar argument can be invoked to substitute the $X$-valued distribution $g$ by a $\mathcal B(Z,X)$-valued distribution. In this case, the convolution defines a new distribution which is again operator-valued, this time in $\mathcal B(Z,Y)$.

\paragraph{Examples.}
The first new case of a convolution operator is differentiation:
\[
 \dot\delta_0*f=\dot f \qquad\longleftrightarrow\qquad s \mathrm F(s)=\mathcal L\{ \dot f\}.
\]
This is also valid when $f$ takes values on $X$, by convoluting with
\[
\smallfrac{\mathrm d}{\mathrm dt}(I_X\otimes \delta_0)=I_X\otimes \dot\delta_0,
\]
where $I_X$ is the identity operator in $X$. This is one of the few cases where we can talk about commutativity: if $f$ takes values in $B(X,Y)$, we can write
\[
(I_Y\otimes \dot\delta_0)* f=\dot f=f*(I_X\otimes \dot\delta_0),
\]
which is the time domain form of
\[
(sI_Y)\mathrm F(s)=s\mathrm F(s)=\mathrm F(s)(s I_X).
\]
Delays are also causal convolution operators
\[
	(I\otimes \delta_{t_0})*f=f(\cdot-t_0)\qquad\longleftrightarrow  \qquad e^{-st_0}\mathrm F(s).
\]
In this case either $f$ is $X$-valued and $I$ is the identity in $X$, or $f$ is $\mathcal B(X,Y)$-valued and $I$ is the identity in $Y$. We can then combine derivatives and delays in the form of an equation: given a causal function $g$, we look for a causal function $y$ satisfying
\[
\dot y=g-g(\cdot-t_1), \qquad \mbox{where $t_1>0$}.
\]
This is then equivalent to
\[
\mathrm Y(s)=\frac{1-e^{-st_1}}{s} \mathrm G(s).
\]
The inversion theory will allow us to see how the solution operator $g\mapsto y$ is also a convolution operator.

\paragraph{Terminology.} If $f$ is a $\mathcal B(X,Y)$-valued causal Laplace transformable distribution, its Laplace transform $\mathrm F:\mathbb C_+\to\mathcal B(X,Y)$ is often referred to as its {\bf transfer function}. In more analytic contexts, the transfer function is also called the {\bf symbol} of the operator.

\paragraph{Inversion theory.} There is an unfortunate fact in the theory of Laplace transforms, stemming from the fact that few classical smoothness properties of a distribution (the distribution equals a function which is continuous or smoother) can be directly seen in their Laplace transforms. However, by restricting the set of distributions and their Laplace transforms, we can find two sets that correspond one to one. At the distributional sense, the set is:
\begin{framed}\noindent
All causal continuous functions $f:\mathbb R\to \mathbb X$ with polynomial growth as $t\to \infty$ and their distributional derivatives.
\end{framed}
Polynomial growth of $f$ means that $\| f(t)\|$ can be bounded by a polynomial in the variable $t$. The set of all Laplace transforms of these distribution can be characterized as follows:
\begin{framed}\noindent
All analytic functions $\mathrm F:\mathbb C_+\to \mathbb X$ such that
\begin{equation}\label{eq:1.100}
\| \mathrm F(s)\| \le C(\mathrm{Re}\,s)|s|^\mu \qquad \forall s\in \mathbb C_+,
\end{equation}
where $\mu \in \mathbb R$ and $C:(0,\infty)\to (0,\infty)$ is a non-increasing function such that
\begin{equation}\label{eq:1.101}
\sigma^m C(\sigma)\le C_0, \qquad \forall\sigma\in (0,1].
\end{equation}
\end{framed}

\paragraph{Convolution justified.} If $f$ is a $\mathcal B(X,Y)$-valued distribution whose Laplace transform satisfies properties \eqref{eq:1.100}-\eqref{eq:1.101}, and $g$ is an $X$-valued distribution whose Laplace transform satisfies the same type of bound (with different $\mu$ and $C$), from where it is clear that $\mathrm F(s)\mathrm G(s)$ is also the Laplace transform of a causal distribution, which is what we call $f*g$. 

\paragraph{A final complicated example.} For this example we need the Sobolev space $H^1_0(\Omega)$, which can be defined as the closure of the set of $\mathcal C^\infty(\Omega)$ functions that vanish in a neighborhood of $\partial\Omega$, with respect to the Sobolev norm
\[
\|u\|_{1,\Omega}^2:=\| \nabla u\|_\Omega^2+\| u\|_\Omega^2=\int_\Omega |\nabla u|^2+\int_\Omega |u|^2.
\]
We now start with $g\in L^2(\Omega)$ and look for
\begin{equation}\label{eq:1.12}
u\in H^1_0(\Omega)\qquad \mbox{s.t.}\qquad (\nabla u,\nabla v)_\Omega+s^2(u,v)_\Omega=(g,v)_\Omega \quad\forall v\in H^1_0(\Omega).
\end{equation}
With some help of the Lax-Milgram lemma, it is not complicated to see that the solution operator $L^2(\Omega)\ni g\mapsto u\in H^1_0(\Omega)$ is well defined and satisfies
\[
\| u\|_{1,\Omega}\le C(\mathrm{Re}\,s) \| g\|_\Omega, \qquad C(\sigma):=\frac1{\sigma\min\{1,\sigma\}}.
\]
It is also easy to see how the solution operator is an analytic function of the parameter $s\in \mathbb C_+$. Therefore, there exists a $\mathcal B(L^2(\Omega),H^1_0(\Omega))$-valued causal distribution $f$ whose Laplace transform is the operator $\mathrm F(s)$ that solves problem \eqref{eq:1.12}. We keep the same letter $g$ for a causal function $g:\mathbb R\to L^2(\Omega)$ with polynomial growth. Then $\mathrm F(s)\mathrm G(s)$ is analytic in $\mathbb C_+$, with values in $H^1_0(\Omega)$. It is also the Laplace transform of a causal $H^1_0(\Omega)$-valued distribution $u=f*g$. This is the very weak form of the problem looking for a causal $H^1_0(\Omega)$-valued distribution $u$ such that
\[
(\nabla u,\nabla v)_\Omega+(\ddot u,v)_\Omega=(g,v)_\Omega\qquad\forall v\in H^1_0(\Omega).
\]
If $u$ were a smooth function of time (it might not be depending on how smooth $g$ is), then we could characterize $u:[0,\infty)\to H^1_0(\Omega)$ in a stronger form
\[
(\nabla u(t),\nabla v)_\Omega+(\ddot u(t),v)_\Omega=(g(t),v)_\Omega \qquad \forall v\in H^1_0(\Omega), \qquad t\ge 0,
\]
with vanishing initial conditions
\[
u(0)=0, \qquad \dot u(0)=0.
\]
Readers acquainted with basic Sobolev space theory (the one used for the most elementary elliptic equations) will recognize a weak form of the wave equation
\[
\partial_{tt} u=\Delta u+g
\]
with homogeneous Dirichlet boundary conditions and homogeneous initial condition. The moral of the story is that the solution operator for this problem is a causal distributional operator.

\subsection*{Credits}

Laplace transforms of (scalar) distributions are an important part of the original theory designed by Laurent Schwartz. A readable, while quite general, introduction is given in \cite[Chapter 6]{Schwartz:1961}. Note that the general case allows for the Laplace transform to be defined in any right semiplane of $\mathbb C$ (we only pay attention here to the one with positive real part). The problem of defining convolutions of distributions is also a classic in modern analysis \cite[Chapter 3]{Schwartz:1961}, especially because many different pairs of distributions can be convoluted, but not every pair. However, the causal case is much simpler, since any pair of causal distributions can be put in convolution. There are no easy references to learn vector-valued distributions. A very comprehensive treatment is given in \cite{Treves:1965}, and a much more concise introduction can be found in \cite{DaLi:1992}. For an introduction tailored to our needs, the reader is referred to \cite{Sayas:2013}. The inversion theorem for the Laplace transform is part of the general theory, but the class of Laplace transforms that we present here, as well as their time-domain representatives, is inherited from \cite{Sayas:2013}, as a further refinement of a class introduced in \cite{LaSa:2009a}. This class of symbols is the one that appears systematically in the treatment of exterior problems for the wave equation \cite{LaSa:2009a}. 


\section{Multistep convolution quadrature}

The goal of this section is the presentation and justification of some discrete quadrature approximations of convolutions
\[
y(t)=\int_0^t f(\tau) g(t-\tau)\mathrm d\tau
\]
and convolution equations
\[
\int_0^t f(\tau) g(t-\tau)\mathrm d\tau=h(t).
\]
\begin{framed}\noindent
The discretization will be developed on a uniform grid of time-step $\kappa>0$
\[
t_n:=n\kappa \quad n\ge 0.
\]
Data will always been dealt with in the time domain. This will lead to discrete convolutions
\[
y(t_n)\approx y_n:=\sum_{m=0}^n \omega_m^{\mathrm F}(\kappa) g(t_{n-m}) \qquad n\ge 0,
\]
and discrete convolution equations
\[
\sum_{m=0}^n \omega_m^{\mathrm F}(\kappa) g_{n-m} = h(t_n) \qquad n\ge 0.
\]
\end{framed}

\subsection{A backward differentiation approach}

\paragraph{A very simple model problem.} Given a causal function $g:\mathbb R\to \mathbb R$ and $c>0$, we look for causal $y$ such that
\begin{equation}\label{eq:2.1}
\dot y-c y=g.
\end{equation}
The causality of $y$ means exactly that $y(t)=0$ for $t<0$. If we want to impose an initial condition $y(0)=y_0$, the way to go is to incorporate it to the right-hand-side
\[
\ddot y-c y=g+y_0\delta_0.
\]
We will not deal with this case, since we will insist on the right-hand-side of \eqref{eq:2.1}. The reader might think it is quite a stupid notion to deal only with homogeneous initial conditions. (It is also true that \eqref{eq:2.1} is a linear equation for with the exact solution can be given a closed formula.) The emphasis of this exposition is in picking up a very simple example to try to understand what we will do for more complicated problems where time-stepping strategies are far from obvious.

\paragraph{Backward Euler differentiation.} The starting point for applying the Backward Euler (BE) method to \eqref{eq:2.1} is the simple discrete differentiation formula:
\begin{equation}\label{eq:2.2}
\dot y(t_n)\approx \smallfrac1\kappa (y(t_n)-y(t_{n-1})).
\end{equation}
We can then define the BE approximation to \eqref{eq:2.1} by
\begin{equation}\label{eq:2.3}
\smallfrac1{\kappa} (y_n-y_{n-1})-c y_n=g(t_n).
\end{equation}
We want the sequence $\{ y_n\}$ to be causal, which means that $y_n=0$ for $n<0$ (but not for $n=0$). Since the discrete differentiation formula \eqref{eq:2.2} only uses one point in the past, we will only need to set $y_{-1}=0$. In this case we will have that $g(0)=g(t_0)=0$ implies $y_0=0$. If $g(0)\neq 0$, then the derivative of the solution of \eqref{eq:2.1} needs to jump at $t=0$ and therefore \eqref{eq:2.2} is not a very good approximation of a quantity that just does not exist at the point $t_0$. Once again, this is not our concern, because the kind of problems we will be dealing with start smoothly from zero. We can write \eqref{eq:2.3} in a more explicit form
\[
y_n =\frac1{1-\kappa c}y_{n-1}+\frac\kappa{1-\kappa c}g(t_n),
\]
and work backwards by induction until we reach $n=0$ to get
\begin{equation}\label{eq:2.4}
y_n=\kappa \sum_{m=0}^n \frac1{(1-\kappa c)^{m+1}} g(t_{n-m}).
\end{equation}
The proper BE method would stop the sum at $m=1$ and never use the value $g(0)$, i.e., \eqref{eq:2.4} is the proper BE method only when $g(0)=0$.

\paragraph{Remark.} The expression \eqref{eq:2.4} is somewhat worrying for $1-\kappa c$ might be zero. In principle one could always think of taking $\kappa$ small enough to ensure that \eqref{eq:2.4} makes sense. (Actually it only fails if $\kappa c=1$.) But you might keep on wondering, were we not using an implicit method? How come this might not work? The reason is actually that when we think of implicit methods and stability issues, we automatically assume $c<0$ (we look for stable approximations of stable problems!), so there is nothing to worry about. We need this form for some future computations.

\paragraph{A more convolutional point of view.} Let me remind you that the causal solution of \eqref{eq:2.1} is given by the convolution expression
\[
y(t)=\int_0^t e^{c\tau} g(t-\tau)\mathrm d\tau.
\]
We focus on the point $t_n$ and write
\begin{equation}\label{eq:2.5}
y(t_n)=\sum_{m=0}^n \int_{t_{m-1}}^{t_m} e^{c\tau} g(t_n-\tau)\mathrm d\tau.
\end{equation}
(The reader should not be too concerned with us including the interval $(t_{-1},t_0)$, where everything should be zero. This is part of the method and we have already warned that when $g(0)\neq 0$ the method should be modified.) The BE method makes the slighly strange approximation
\[
\int_{t_{m-1}}^{t_m} e^{c\tau} g(t_n-\tau)\mathrm d\tau\approx \frac{\kappa}{(1-\kappa c)^{m+1}} g(t_n-t_m),
\]
which we will not try to justify.

\paragraph{Challenges we will not accept.} We can also think of our very primitive wave equation
\[
\ddot y+c^2 y=g
\]
when we approximate
\[
\ddot y(t_n)=\frac1{\kappa^2}(y(t_n)-2y(t_{n-1})+y(t_{n-2})),
\]
which is the first order approximation of the second derivative corresponding to \eqref{eq:2.2}. The associated discrete method is
\[
\frac1{\kappa^2}(y_n-2y_{n-1}+y_{n-2}) +c^2 y_n=g(t_n),
\]
starting with causal conditions $y_{-2}=y_{-1}=0$. While it is possible to find a formula
\begin{equation}\label{eq:2.6}
y_n=\kappa^2\sum_{m=0}^n \alpha_m (c\kappa) g(t_{n-m}),
\end{equation}
we will not waste our time trying to figure out these coefficients. What is important is the fact that \eqref{eq:2.6} approximates
\[
y(t_n)=c^{-1} \int_0^{t_n} \sin(c\tau) g(t_n-\tau)\mathrm d\tau.
\]
If we go back to \eqref{eq:2.1} but now approximate the derivative using a double backward differentiation formula
\[
\dot y(t_n)\approx\smallfrac1\kappa (\smallfrac32 y(t_n)-2y(t_{n-1})+\smallfrac12 y(t_{n-2})),
\]
we end up with the BDF2 method
\[
\smallfrac1\kappa(\smallfrac32 y_n-2y_{n-1}+\smallfrac12 y_{n-2})+cy_n=g(t_n).
\]
This recurrence can also be solved to obtain a formula ressembling \eqref{eq:2.4}. These formulas become really cumbersome to obtain (especially when it is unclear why we need them), so instead of hitting our heads against the wall repeatedly, we are going to move on to understand the language of finite difference equations.

\subsection{The language of $\zeta$ transforms}\label{sec:2.2}

\paragraph{In two words.} The $\zeta$ transform is to initial value problems for difference equations (recurrences) what the Laplace transform is to initial value problems for ODEs. However, it is mainly a formal transform, involving series whose convergence will not concern us at all.

\paragraph{The $\zeta$ transform.} Our input is a sequence $\{ y_n\}$, which we can consider to start at $n=0$ or to be causal ($y_n=0$ for $n\le -1$). We then associate the formal series
\[
\mathrm Y(\zeta):=\sum_{n=0}^\infty y_n \zeta^n.
\]
This series is to be understood in a purely algebraic way, as a sort of polynomial with infinitely many coefficients. A simple operation to be performed with causal sequences is the displacement to the right
\[
(y_0,y_1,y_2,\ldots) \quad\longmapsto \quad (0,y_0,y_1,\ldots)
\]
which we can understand with the full causal sequence as $\{ y_n\}\mapsto \{ y_{n-1}\}$. This operation is very easy to describe with $\zeta$ series:
\[
\mathrm Y(\zeta)\quad\longmapsto\quad \zeta\mathrm Y(\zeta)=\sum_{n=1}^\infty y_{n-1}\zeta^n.
\]
The second important operation that $\zeta$ transforms describe in a simple way is discrete convolution. The convolution
\[
\sum_{m=0}^n a_m b_{n-m}
\]
corresponds to the product of the $\zeta$ transforms
\[
\mathrm A(\zeta)\mathrm B(\zeta)=\sum_{n=0}^\infty \left(\sum_{m=0}^n a_m b_{n-m} \right)\zeta^n.
\]
In these expressions we might be thinking of scalar quantities (sequences taking values in $\mathbb R$ or $\mathbb C$) or more complicated situations, where, for instance, $\{ y_n\}$ is a sequence in a vector space $X$ and the sequences $\{a_n\}$ and $\{ b_n\}$ take values on spaces of operators where the multiplications $a_mb_{n-m}$ are meaningful.

\paragraph{BE again.} The discrete recurrence 
\[
\smallfrac1\kappa (y_n-y_{n-1}) - cy_n=g_n, \qquad g_n:=g(t_n),
\]
is transformed into
\[
\left(\frac{1-\zeta}\kappa - c\right) \mathrm Y(\zeta)=\mathrm G(\zeta)
\]
and thus yields
\begin{equation}\label{eq:2.7}
\mathrm Y(\zeta)=\frac1{\frac{1-\zeta}\kappa-c}\mathrm G(\zeta).
\end{equation}
Working through the algebra, we can easily write
\begin{equation}\label{eq:2.8}
\frac1{\frac{1-\zeta}\kappa-c}=\frac{\kappa}{(1-\kappa c)-\zeta}=\frac{\kappa}{1-\kappa c} \left(1-\frac\zeta{1-\kappa c}\right)^{-1}=
\sum_{n=0}^\infty \frac{\kappa}{(1-\kappa c)^{n+1}}\zeta^n,
\end{equation}
which means that \eqref{eq:2.7} and \eqref{eq:2.8} are just encoding \eqref{eq:2.4}.

\paragraph{BDF2.} The $\zeta$ transform of the BDF2 equations
\[
\smallfrac1\kappa(\smallfrac32 y_n-2y_{n-1}+\smallfrac12 y_{n-2})+cy_n=g_n, \qquad g_n:=g(t_n)
\]
is
\[
\left(\frac1\kappa\left(\frac32-2\zeta+\frac12\zeta^2\right)-c\right)\mathrm Y(\zeta)=\mathrm G(\zeta),
\]
or in explicit form
\begin{equation}\label{eq:2.9}
\mathrm Y(\zeta)=\frac1{\frac1\kappa\left(\frac32-2\zeta+\frac12\zeta^2\right)-c} \mathrm G(\zeta).
\end{equation}

\paragraph{Implicit differentiation by the trapezoidal rule.} The equation
\[
\dot y-cy=g
\]
can also be approximated using the trapezoidal rule
\begin{equation}\label{eq:2.10}
\frac1\kappa(y_n-y_{n-1})-\frac{c}2 (y_n+y_{n-1})=\frac12 (g_n+g_{n-1}).
\end{equation}
In this case we are not starting from a backward discretization of the derivative, but have found the method working directly on the differential equation. With $\zeta$ transforms, \eqref{eq:2.10} is written as
\[
\left(\frac{1-\zeta}\kappa-c\frac{1+\zeta}2\right)\mathrm Y(\zeta)=\frac{1+\zeta}2\mathrm G(\zeta),
\]
which can be given explicitly as
\begin{equation}\label{eq:2.11}
\mathrm Y(\zeta)=\frac1{\frac1\kappa 2 \frac{1-\zeta}{1+\zeta}-c} \mathrm G(\zeta)
\end{equation}

\paragraph{Summary.} It is instructive to pay attention to the expressions of the recurrences \eqref{eq:2.7}, \eqref{eq:2.9}, and \eqref{eq:2.11}. All of them share the form
\[
\mathrm Y(\zeta)=\frac1{\frac1\kappa \delta(\zeta)-c}\mathrm G(\zeta),
\]
with
\[
\delta(\zeta):=\left\{ \begin{array}{ll} 
	1-\zeta, & \mbox{(BE)},\\
	\frac32-2\zeta+\frac12\zeta^2, & \mbox{(BDF2)},\\
	2\frac{1-\zeta}{1+\zeta}, & \mbox{(TR)}.
\end{array}\right.
\]
We already understood $\delta(\zeta)$ as a backward differentiation formula for BE and BDF2. We could also think of higher order BDF methods
\[
\delta(\zeta):=\sum_{\ell=1}^p \frac1\ell(1-\zeta)^\ell.
\]
While these methods are useful in the parabolic world, we will not be able to use them in wave propagation problems, so we will quietly put them aside. The trapezoidal approximation deserves some special attention. The function
\[
\delta(\zeta)=2\frac{1-\zeta}{1+\zeta}=2(1-\zeta)\sum_{n=0}^\infty (-\zeta)^n=2+4\sum_{n=1}^\infty (-1)^n \zeta^n 
\]
does not represent a simple difference formula for approximation of the derivative of a function. Instead, it can be seen as a long memory approximation of the derivative of a causal function
\[
\dot y(t_n)\approx \frac1{\kappa}\big(2 y(t_n)-4y(t_{n-1})+4y(t_{n-2})-4y(t_{n-3})+\ldots\big)
\]
(the sum is finite because $y$ is causal), or as an implicit method to approximate the derivative. If we know $y$, then $h=\dot y$ is approximated by solving recurrently 
\[
\frac12(h_n+h_{n-1}) =\frac1\kappa (y_n-y_{n-1}), \qquad y_n=y(t_n).
\]
This formula might look surprising to the reader. Just think that differentiating the function is a sort of opposite process to solving a differential equation, so now the data are being differentiated and appear in the right-hand-side.

\subsection{Moving from $s$ to $\zeta$}\label{sec:2.3}

\paragraph{Back to the Laplace transform.} The causal solution to
\[
\dot y-c y=g
\]
is given through its Laplace transform by
\[
\mathrm Y(s) = \frac1{s-c}\mathrm G(s).
\]
The discrete versions we have obtained in the previous subsection fit in the general frame
\[
\mathrm Y(\zeta)=\frac1{\frac1\kappa\delta(\zeta)-c}\mathrm G(\zeta),
\]
or also in implicit form
\begin{equation}\label{eq:2.12}
\frac1\kappa\delta(\zeta)\mathrm Y(\zeta)-c\mathrm Y(\zeta)=\mathrm G(\zeta),
\end{equation}
where we see the approximation of the differentiation operator by a discrete derivative operator, given in $\zeta$ transformed style by $\frac1\kappa\delta(\zeta)$. Our transfer function was $(s-c)^{-1}$. The discrete transfer function is
\[
\left(\frac1\kappa\delta(\zeta)-c\right)^{-1}.
\]
We do not need the expansion of this function into coefficients, because the method itself is built from a recurrence we can just apply. However, just for the sake of the argument, let us formally expand:
\[
\left(\frac1\kappa\delta(\zeta)-c\right)^{-1}=\sum_{n=0}^\infty \omega_n^c(\kappa)\zeta^n.
\]
Then, the recurrence hidden in \eqref{eq:2.12} can be given an explicit form
\[
y_n=\sum_{m=0}^n \omega_n^c(\kappa) g(t_{n-m}).
\]

\paragraph{A first generalization.} Let us think again of our primitive wave equation
\[
\ddot y+c^2 y =g.
\]
In the Laplace domain this is
\[
\mathrm Y(s)=\frac1{s^2+c^2}\mathrm G(s),
\]
which corresponds to
\[
y(t)=c^{-1}\int_0^t \sin(c\tau)g(t-\tau)\mathrm d\tau.
\]
If we use one of our time-stepping methods, we are essentially applying the recurrence
\[
\left(\frac1\kappa\delta(\zeta)\right)^2 \mathrm Y(\zeta)+c^2\mathrm Y(\zeta)=\mathrm G(\zeta),
\]
or in explicit form
\[
\mathrm Y(\zeta)=\frac1{\left(\frac1\kappa\delta(\zeta)\right)^2+c^2}\mathrm G(\zeta).
\]
If we are able to expand
\[
\frac1{\left(\frac1\kappa\delta(\zeta)\right)^2+c^2}=\sum_{n=0}^\infty \omega_n^{c^2}(\kappa)\zeta^n,
\]
then the recurrence becomes a discrete convolution
\[
y_n=\sum_{m=0}^n \omega_n^{c^2}(\kappa)g(t_{n-m}).
\]
Even if we are just rewriting very simple solvers for linear differential equations, we can see how the final discrete convolution formula uses the Laplace transform of the solution operator (in this case $(s^2+c^2)^{-1}$), while data are sampled in the time domain.

\paragraph{Convolutions become discrete.} We next move to a general convolution
\begin{equation}\label{eq:2.13}
y=f*g.
\end{equation}
Here $f$ and $g$ are known. We need $g$ to be a causal function $\mathbb R\to X$. On the other hand $f$ can be any causal Laplace transformable distribution with values in $\mathcal B(X,Y)$. We actually do not need to know $f$, but are happy enough with its Laplace transform $\mathrm F(s)$. The Laplace transform of \eqref{eq:2.13} is 
\[
\mathrm Y(s)=\mathrm F(s)\mathrm G(s).
\]
The discrete equations are then
\begin{equation}\label{eq:2.14}
\mathrm Y(\zeta)=\mathrm F(\smallfrac1\kappa\delta(\zeta))\mathrm G(\zeta), \qquad \mathrm G(\zeta):=\sum_{n=0}^\infty g(t_n)\zeta^n.
\end{equation}
If we are able to expand
\begin{equation}\label{eq:2.15}
\mathrm F(\smallfrac1\kappa\delta(\zeta))=\sum_{n=0}^\infty \omega_n^{\mathrm F}(\kappa)\zeta^n,
\end{equation}
then \eqref{eq:2.14} becomes a discrete convolution
\begin{equation}\label{eq:2.16}
y_n=\sum_{m=0}^n \omega_m^{\mathrm F}(\kappa) g(t_{n-m}),\qquad n\ge 0
\end{equation}
which uses time-domain readings of the data function $g$ combined with a discrete transfer function that arises from the original transfer function and the approximation of the differentiation operator. The discrete convolution \eqref{eq:2.16}, with coefficients computed following \eqref{eq:2.15}, is the {\bf Convolution Quadrature} method to approximate the convolution \eqref{eq:2.13} at discrete time steps. Recall that our hypotheses for $\mathrm F$ included $\mathrm F$ to be analytic in $\mathbb C_+$. If $\delta(0)\in \mathbb C_+$, then the function $\zeta\mapsto\mathrm F(\frac1\kappa\delta(\zeta))$ is analytic in a neighborhood of $\zeta=0$ and the expansion \eqref{eq:2.15} is just a Taylor expansion.

\paragraph{Convolution equations.} A convolution equation can be seen (in the Laplace domain) as
\begin{equation}\label{eq:2.17}
\mathrm F(s)\mathrm G(s)=\mathrm H(s) \qquad \mbox{or}\qquad \mathrm G(s)=\mathrm F(s)^{-1}\mathrm H(s).
\end{equation}
The version in the right-hand-side is an explicit convolution, with the caveat that $\mathrm F(s)^{-1}$ might not be known. The Convolution Quadrature method is then given by the recurrence that is equivalent to the $\zeta$ transform equation:
\[
\mathrm F(\smallfrac1\kappa\delta(\zeta)) \mathrm G(\zeta)=\mathrm H(\zeta), \qquad \mathrm H(\zeta):=\sum_{n=0}^\infty h(t_n)\zeta^n.
\]
In discrete times, this is equivalent to
\begin{equation}\label{eq:2.18}
\sum_{m=0}^n \omega_m^{\mathrm F}(\kappa) g_{n-m} = h(t_n), \qquad n\ge 0,
\end{equation}
or, in a slighly more explicit form, to
\[
\omega_0^{\mathrm F}(\kappa)g_n=h(t_n)-\sum_{m=1}^n \omega_m^{\mathrm F}(\kappa) g_{n-m} =
h(t_n)-\sum_{m=0}^{n-1} \omega_{n-m}^{\mathrm F}(\kappa) g_m.
\]
The coefficients $\{\omega_n^{\mathrm F}(\kappa)\}\subset \mathcal B(X,Y)$ are computed using \eqref{eq:2.15}. The method is well defined when $\omega_0^{\mathrm F}(\kappa)$ is invertible. 

\paragraph{Just in case.} The way we have written the convolution equation \eqref{eq:2.17} raises an interesting issue. Assume that we happen to know both $\mathrm F(s)$ and $\mathrm F(s)^{-1}$, and that we can expand
\begin{equation}\label{eq:2.19}
\mathrm F(\smallfrac1\kappa\delta(\zeta))=\sum_{n=0}^\infty \omega_n^{\mathrm F}(\kappa)\zeta^n,\qquad
\mathrm F(\smallfrac1\kappa\delta(\zeta))^{-1}=\sum_{n=0}^\infty \omega_n^{\mathrm F^{-1}}(\kappa)\zeta^n.
\end{equation}
We wonder whether \eqref{eq:2.18} (the discretization of the convolution equation) is the same method as
\[
g_n=\sum_{m=0}^n \omega_m^{\mathrm F^{-1}}(\kappa) h(t_{n-m}).
\]
The answer is yes.
The key to understanding why lies in the fact that both expansions \eqref{eq:2.19} are Taylor expansions. In fact, there is an underlying assumption in the background: we assume that $\mathrm F$ is a transfer function in the conditions of Section \ref{sec:1.4}, and so is $\mathrm F^{-1}$. This means that for all $s\in \mathbb C_+$, we need $\mathrm F(s)$ to be invertible and we want the function $\mathbb C_+\ni s\mapsto F(s)^{-1}\in \mathcal B(Y,X)$ to be analytic (we get this for free) and to admit a bound like \eqref{eq:1.100}. 
In particular, these hypotheses imply (see Section \ref{sec:1.4}) that the convolution equation $f*g=h$ is solved with another convolution equation $g=f^{-1}*g$, where $\mathcal L\{ f^{-1}\}(s)=\mathrm F(s)^{-1}$.
In this case
\[
\mathrm F(\smallfrac1\kappa\delta(\zeta))\mathrm F(\smallfrac1\kappa\delta(\zeta))^{-1}=I 
\]
for $\zeta$ small enough, which implies that
\[
\omega_0^{\mathrm F}(\kappa)^{-1}=\omega_0^{\mathrm F^{-1}}(\kappa), \qquad \sum_{m=0}^n \omega_m^{\mathrm F}(\kappa)\omega_{n-m}^{\mathrm F^{-1}}(\kappa)=0 \quad n\ge 1.
\]
This proves that
\[
\mathrm F(\smallfrac1\kappa\delta(\zeta)) \mathrm G(\zeta)=\mathrm H(\zeta)
\qquad
\mbox{and}
\qquad
 \mathrm G(\zeta)=\mathrm F(\smallfrac1\kappa\delta(\zeta))^{-1}\mathrm H(\zeta)
\]
deliver the same sequence $\{ g_n\}$.

\paragraph{Why implicit.} It is tempting to think whether an explicit approximation of the derivative would work. While in principle there might not be any problem for it (you just need to compute what $\delta(\zeta)$ is, plug it in \eqref{eq:2.15}, and expand), it is not very clear what we would gain from it. As already emphasized, we are not dealing with linear differential equations. We are merely using them to motivate the methods. The actual transfer functions are analytic functions of the variable $s$ and using one $\delta(\zeta)$ or another might not make much of a difference at the point of implementing the method. It has to be said, nevertheless, that for applications in wave propagation, using the discrete differentiation operators associated to explicit approximations of the derivative does not work. We can also just wash our hands and claim that we are dealing with causal problems (causal convolutions), and we want to keep everything causal. Explicit (forward) differentiation breaks causality by looking into the future.

\subsection{A bold aproach in the Laplace domain}

\paragraph{A new and fast introduction of BE-CQ.} Let us just focus on causal convolutions
\[
y=f*g \qquad \mbox{i.e.}\qquad \mathrm Y(s)=\mathrm F(s)\mathrm G(s).
\]
A particular case is $\mathrm F(s)=s$, corresponding to differentiation $y=\dot g$. We fix the time step $\kappa$ but do not pay special attention to the discrete times. The Euler discrete backward derivative approximation to $y=\dot g$ is
\[
y\approx \smallfrac1\kappa(g-g(\cdot-\tau))=:y_\kappa.
\]
The Laplace domain form of this equation is
\[
\mathrm Y_\kappa(s)=\smallfrac1\kappa(1-e^{-s\kappa}) \mathrm G(s).
\]
This can be written as
\begin{equation}\label{eq:2.20}
\mathrm Y_\kappa(s)=s_\kappa \mathrm G(s),
\end{equation}
where
\[
s_\kappa=\smallfrac1\kappa(1-e^{-s\kappa})=\smallfrac1\kappa\delta(e^{-s\kappa}), \qquad \delta(\zeta)=1-\zeta.
\]

\paragraph{BDF2 and TR.} We can try the same ideas with our other two particular methods introduced in Sections \ref{sec:2.2} and \ref{sec:2.3}. The BDF2 approximation of the derivative applied to $y =\dot g$ is
\[
y \approx y_\kappa:=\smallfrac1\kappa(\smallfrac32 g-2g(\cdot-\kappa)+\smallfrac12 g(\cdot-2\kappa)),
\]
which can be written as
\[
\mathrm Y_\kappa (s)=\underbrace{\smallfrac1\kappa(\smallfrac32-2e^{-s\kappa}+\smallfrac12e^{-2s\kappa})}_{s_\kappa}\mathrm G(s).
\]
In the trapezoidal rule $y_\kappa$ is approximated implicitly
\[
\smallfrac12 (y_\kappa+y_\kappa(\cdot-\kappa))=\smallfrac1\kappa(g-g(\cdot-\kappa))
\]
and we similarly obtain \eqref{eq:2.20} with
\[
s_\kappa:=\frac2\kappa \frac{1-e^{-s\kappa}}{1+e^{-s\kappa}}=\frac1\kappa \delta(e^{-s\kappa}). 
\]

\paragraph{Convolution Quadrature again.} We then extend the idea of approximating
\[
\mathrm Y(s)=s\mathrm G(s)\qquad \mbox{by}\qquad \mathrm Y_\kappa(s)=s_\kappa \mathrm G(s),
\]
for more general transfer functions, approximating then
\[
\mathrm Y(s)=\mathrm F(s)\mathrm G(s)\qquad \mbox{by}\qquad \mathrm Y_\kappa(s)=\mathrm F(s_\kappa) \mathrm G(s).
\]
We next have to figure out what this method is. Recall that we are using
\[
\mathrm F(s_\kappa), \qquad\mbox{with}\qquad s_\kappa=\smallfrac1\kappa\delta(e^{-s\kappa}).
\]
We start by Taylor expanding
\[
\mathrm F(\smallfrac1\kappa\delta(\zeta))=\sum_{n=0}^\infty \omega_n^{\mathrm F}(\kappa)\zeta^n,
\]
and then substituting
\begin{equation}
\mathrm F(s_\kappa)=
\sum_{n=0}^\infty \omega_n^{\mathrm F}(\kappa) e^{-n\kappa s}
=\sum_{n=0}^\infty \omega_n^{\mathrm F}(\kappa) e^{-t_n s}.
\end{equation}
We finally recognize that $\mathrm F(s_\kappa)$ is the Laplace transform of the distribution
\[
\sum_{n=0}^\infty \omega_n^{\mathrm F}(\kappa)\otimes \delta_{t_n}
\]
and therefore $Y_\kappa(s)=\mathrm F(s_\kappa) \mathrm G(s)$ is the Laplace transform of 
\begin{equation}\label{eq:2.22}
y_\kappa =\sum_{m=0}^\infty (\omega_m^{\mathrm F}(\kappa)\otimes \delta_{t_m}) * g =
\sum_{m=0}^\infty \omega_m^{\mathrm F}(\kappa)g(\cdot-t_m). 
\end{equation}
In this form $g$ is allowed to be discontinuous in time and $y_\kappa$ is obtained as a function of continuous time. The sum in \eqref{eq:2.22} is finite for any given value of the time variable. For instance
\[
y_\kappa(t)=\sum_{m=0}^n \omega_m^{\mathrm F}(\kappa)g (t-t_m), \qquad t_{n}\le t < t_{n+1}.
\]
In particular, with $t=t_n$ we obtain the discrete Convolution Quadrature method of Section \ref{sec:2.3}
\[
y_n:=y_\kappa(t_n)=\sum_{m=0}^n \omega_m^{\mathrm F}(\kappa)g(t_{n-m}).
\]
This means that the discrete sequence $\{ y_n\}$ obtained in Section \ref{sec:2.3} is actually composed of time samples of a continuous-in-time function $y_\kappa$. Of course, if we want to evaluate $y_\kappa$ at a point that is not on the time grid, we have to use $g$ at different points, instead of only at the points $t_m$. The case of convolution equations follows similar lines: to approximate
\[
\mathrm F(s)\mathrm G(s)=\mathrm H(s)
\] 
with
\[
\mathrm F(s_\kappa)\mathrm G_\kappa(s)=\mathrm H(s)
\]
is equivalent to writing the equation
\begin{equation}\label{eq:2.23}
\sum_{m=0}^\infty \omega_m^{\mathrm F}(\kappa) g_\kappa(\cdot-t_m)=h.
\end{equation}
Ideally we could solve \eqref{eq:2.23} progressively in intervals:
\[
\omega_0^{\mathrm F}(\kappa)g(t)=h(t)-\sum_{m=1}^n \omega_m^{\mathrm F}(\kappa)g(t-t_m), \qquad t\in [t_{n-1},t_n),\quad n\ge 0.
\]
This would give $g_\kappa$ one interval after another. What we will focus on is the values $g_n:=g_\kappa(t_n)$, which satisfy
\[
\omega_0^{\mathrm F}(\kappa)g_n=h(t_n)-\sum_{m=1}^n \omega_m^{\mathrm F}(\kappa)g(t_{n-m}), \qquad n\ge 0,
\]
sending us back to the CQ method of Section \ref{sec:2.3}.

\paragraph{Lubich's notation.} We have favored some sort of convolutional notation in these notes. In his expositions of the method, Christian Lubich has adopted a clever operational notation that is often useful. The idea is simple. Differentiation 
\[
y=\partial g=\dot g
\]
corresponds to the operator `multiply by $s$'. We then admit the notation
\[
y=\mathrm F(\partial) g=f*g, \qquad \mathrm F=\mathcal L\{ f\},
\]
that makes the transfer function apparent even in the time domain. In CQ, there is a discrete differentiation operator
\[
\partial_\kappa g \approx \dot g, \qquad \mathcal L\{ \partial_\kappa g\}(s)=s_\kappa\mathrm G(s), \qquad s_\kappa=\smallfrac1\kappa\delta(e^{-s\kappa}),
\]
and the CQ method is denoted
\[
\mathrm F(\partial_\kappa)g=\sum_{m=0}^\infty\omega_m^{\mathrm F}(\kappa) g_\kappa(\cdot-t_m)=f_\kappa * g,
\qquad f_\kappa:=\sum_{m=0}^\infty \omega_m^{\mathrm F}(\kappa) \otimes \delta_{t_m},
\]
which makes $\mathcal L\{ f_\kappa\}(s)=\mathrm F(s_\kappa)$.

\subsection{Another point of view}\label{sec:2.5}

\paragraph{Aims and tools.}
In this section we are going to give Lubich's original introduction to the BE-CQ method. We will restrict to some particular transfer functions $\mathrm F:\mathbb C_+\to \mathcal B(X,Y)$ satisfying
\[
\| \mathrm F(s)\| \le C(\mathrm{Re}\,s) |s|^\mu \qquad \forall s\in \mathbb C_+, \quad \mu<-1. 
\]
This condition is enough to show that the Laplace inversion formula
\begin{equation}
f(t)=\frac1{2\pi\imath}\int_{\sigma-\imath\infty}^{\sigma+\imath\infty} e^{st}\mathrm F(s)\mathrm ds=
\frac1{2\pi}\int_{-\infty}^\infty e^{(\sigma+\imath\omega)\,t} \mathrm F(\sigma+\imath\omega)\mathrm d\omega
\end{equation}
defines a causal continuous function whose Laplace transform is $\mathrm F$, independently of the value of $\sigma$ taken for the inversion path.
The hypotheses also ensure that if $a>\sigma$, then
\begin{equation}\label{eq:2.25}
\frac1{n!}\mathrm F^{(n)}(a)=-\frac1{2\pi\imath} \int_{\sigma-\imath\infty}^{\sigma+\imath\infty} \frac1{(s-a)^{n+1}}\mathrm F(s)\mathrm d s.
\end{equation}
Equation \eqref{eq:2.25} is a somewhat uncommon expression deriving from the Cauchy integral formulas. It can be proved using the Cauchy formulas: integrating along a rectangular countour with the left edge in the line $\sigma+\imath\mathbb R$, the right edge in $R+\imath \mathbb R$ and the upper and bottom edges in $\mathbb R\pm\imath R$, and then sending $R\to \infty$.

\paragraph{A mixed formulation of causal convolution.} All the following steps can be justified: just stay with the flow
\begin{eqnarray*}
(f*g)(t) 
	&=&\int_0^t f(\tau) g(t-\tau)\mathrm d\tau\\
	&=& \int_0^t \left(\frac1{2\pi\imath}\int_{\sigma-\imath\infty}^{\sigma+\imath\infty} e^{s\tau}\mathrm F(s)\mathrm ds\right)g(t-\tau)
		\mathrm d\tau\\
	&=& \frac1{2\pi\imath}\int_{\sigma-\imath\infty}^{\sigma+\imath\infty}\mathrm F(s)
		\underbrace{\left(\int_0^t e^{s\tau}g(t-\tau)\mathrm d\tau\right)}_{y(t)=y(t;s)}\mathrm ds.
\end{eqnarray*}
We have now reached a point where $g$ appears in a scalar convolution equation, corresponding to the solution of
\[
\dot y=sy+g,
\]
with the implicit assumption of causality that can also be phrased as $y(0)=0$.

\paragraph{Introducing the ODE solver.}
We restrict our attention to $t=t_n$ and approximate
\[
\int_0^{t_n} e^{s\tau}g(t_n-\tau)\mathrm d\tau\approx \sum_{m=0}^n \frac{\kappa}{(1-\kappa s)^{m+1}} g(t_{n-m}),
\]
that is, we approximate $y(t_n)$ using the BE method. This leads to
\[
(f*g)(t_n) \approx \sum_{m=0}^n 
\underbrace{\left(\frac1{2\pi\imath}\int_{\sigma-\imath\infty}^{\sigma+\imath\infty}
			\frac{\kappa}{(1-\kappa s)^{m+1}}\mathrm F(s)\mathrm ds\right)}_{\omega_m^{\mathrm F}(\kappa)} g(t_{n-m})
\]
The final step consists of the figuring out what the coefficients $\omega_n^{\mathrm F}(\kappa)$ are. From \eqref{eq:2.25} and some elementary arguments, we obtain
\begin{eqnarray*}
\omega_n^{\mathrm F}(\kappa) 
	&=& \frac{(-1)^n}{\kappa^n}\left(
			-\frac1{2\pi\imath}\int_{\sigma-\imath\infty}^{\sigma+\imath\infty} \frac1{(s-1/\kappa)^{n+1}}\mathrm F(s)\mathrm ds\right)\\
	&=& \frac{(-1)^n}{\kappa^n}\frac1{n!} \mathrm F^{(n)}(1/\kappa)\\
	&=&\frac1{n!} \frac{\mathrm d^n}{\mathrm d\zeta^n}\left(\mathrm F\Big(\frac{1-\zeta}\kappa\Big)\right)\Bigg|_{\zeta=0}.
\end{eqnarray*}
Therefore
\[
\mathrm F\left(\frac{1-\zeta}\kappa\right)=\sum_{n=0}^\infty \omega_n^{\mathrm F}(\kappa) \zeta^n,
\]
and we have found the BE-CQ method in a slightly different way.

\subsection{A discrete operational calculus}

\paragraph{Rephrasing.} What we have done so far can be understood in many different ways, but it can be synthesized with several simple statements. Given a causal $\mathcal B(X,Y)$-valued distribution with Laplace transform $\mathrm F:\mathbb C_+\to\mathcal B(X,Y)$, we 
\begin{equation}\label{eq:2.40}
\mbox{approximate } f \approx \sum_{n=0}^\infty \omega_n^{\mathrm F}(\kappa)\otimes \delta_{t_n}, 
\mbox{ where } \mathrm F(\smallfrac1{\kappa}\delta(\zeta))=\sum_{n=0}^\infty \omega_n^{\mathrm F}(\kappa)\zeta^n.
\end{equation}
The function $\delta(\zeta)$ can be understood as a generator for an approximation for the differentiation operator, so that
\begin{equation}\label{eq:2.41}
\smallfrac1{\kappa}\delta(\zeta)=\smallfrac1\kappa\sum_{n=0}^\infty \delta_n \zeta^n
\end{equation}
is the symbolic form of the operator
\[
\partial_\kappa g:=\smallfrac1\kappa \sum_{n=0}^\infty \delta_n g(\cdot-n\kappa).
\]
Equations \eqref{eq:2.40} and \eqref{eq:2.41} can also be written in the Laplace domain, via the definition of
\[
s_\kappa:=\smallfrac1\kappa\delta(e^{-s\kappa}) =\smallfrac1\kappa \sum_{n=0}^\infty \delta_n e^{-s\,t_n}\approx s,
\]
to give approximations
\begin{equation}\label{eq:2.42}
\mathrm F(s_\kappa)=\sum_{n=0}^\infty \omega_n^{\mathrm F}(\kappa) e^{-s\, t_n} \approx \mathrm F(s).
\end{equation}
As already mentioned, there are different ways of notating the convolution operator defined by \eqref{eq:2.40}. One way is to name
\[
f_\kappa:=\sum_{n=0}^\infty \omega_n^{\mathrm F}(\kappa)\otimes \delta_{t_n}
\]
and then note that $f_\kappa * g\approx f* g$. The other way is to name the operator directly $\mathrm F(\partial_\kappa) g$. {\em There is a caveat to this condensed introduction.} Even if the CQ expressions $\mathrm F(\partial_\kappa)g=f_\kappa* g$ exist in continuous time, they are only computed in discrete times, so knowledge of the output of these operators is purely given at discrete time steps.

\paragraph{Associativity.} Let now $\mathrm F_1:\mathbb C_+\to \mathcal B(X,Y)$ and $\mathrm F_2:\mathbb C_+\to\mathcal B(Z,X)$. It is clear from the Laplace transform of the convolution quadrature operators that
\[
\mathrm F_1(s_\kappa) \big( \mathrm F_2(s_\kappa)\mathrm G(s)\big)=\big(\mathrm F_1(s_\kappa) \mathrm F_2(s_\kappa)\big) \mathrm G(s),
\]
which is another way of writing
\begin{equation}\label{eq:2.43}
\mathrm F_1(\partial_\kappa)\big(\mathrm F_2(\partial_\kappa) g\big) = 
\big(\mathrm F_1(\partial_\kappa)\mathrm F_2(\partial_\kappa)\big) g,
\end{equation}
or equivalently
\[
f_{1,\kappa} * (f_{2,\kappa}*g)= (f_{1,\kappa}*f_{2,\kappa}) * g.
\]
The latter expression emphasizes the convolutional form of all operators, from where associativity is clear. The interest of \eqref{eq:2.43} might not be apparent at first sight. In practice, we are going to apply the left-hand-side of \eqref{eq:2.43}, that is, first one discrete convolution process, and then the second one. However, this is equivalent to applying a simple CQ process, which we will never compute (it would mean composing sequences of operators), but we can use for analysis. 

\paragraph{Forward convolutions and equations.} The associativity idea can be further exploited to see why analyzing CQ for convolution equations is equivalent to analyzing CQ for the inverse operator, even if this is only known at the theoretical level. If $\mathrm F:\mathbb C_+\to \mathcal B(X,Y)$ and $\mathrm F^{-1}:\mathbb C_+\to \mathcal B(Y,X)$ are operators such that
\[
\mathrm F(s)\mathrm F^{-1}(s)=\mathrm I_Y, \qquad \mathrm F^{-1}(s)\mathrm F(s)=I_X, \qquad \forall s\in \mathbb C_+,
\]
(or, in other words, $\mathrm F^{-1}(s)=(\mathrm F(s))^{-1}$), then
\[
\mathrm F(s_\kappa)\mathrm F^{-1}(s_\kappa)=\mathrm I_Y, \qquad \mathrm F^{-1}(s_\kappa)\mathrm F(s_\kappa)=I_X,
\]
which means that the discrete operators associated to $\mathrm F(s_\kappa)$ and $\mathrm F^{-1}(s_\kappa)$ are inverse of each other. Therefore, the analysis of the equation $f_\kappa *g=h$ (here $g$ is unknown) is equivalent to the analysis of the inverse operator $g=f_\kappa^{-1}*h$.

\subsection{Discrete transfer functions}

\paragraph{Motivation.} Just for the sake of it, let us explore the problem of inverting discrete convolution equations. Let $\{ f_n\} \subset \mathcal B(X,Y)$ and $\{ h_n\}\subset Y$ be given sequences. We look for a sequence $\{ g_n\} \subset X$ satisfying
\begin{equation}\label{eq:2.100}
\sum_{m=0}^n f_m g_{n-m}=h_n.
\end{equation}
It is clear that if $f_0$ is invertible, then the sequence $\{g_n\}$ can be computed with the semi-implicit recurrence
\begin{equation}\label{eq:2.101}
g_n=f_0^{-1}(h_n-\sum_{m=1}^n f_m g_{n-m}).
\end{equation}
With $\zeta$ transforms the convolution equation \eqref{eq:2.100} can be written as
\begin{equation}\label{eq:2.102}
\mathrm F(\zeta)\mathrm G(\zeta)=\mathrm H(\zeta).
\end{equation}

\paragraph{The discrete transfer function.}
At least formally we can try to invert $\mathrm F(\zeta)$ and send it to the right-hand-side:
\begin{equation}\label{eq:2.103}
\mathrm G(\zeta)=\mathrm R(\zeta)\mathrm H(\zeta),
\end{equation}
hoping to find a sequence $\{ r_n\}$ that rewrites \eqref{eq:2.101} in the explicit form
\[
g_n=\sum_{m=0}^n r_m h_{n-m}.
\]
The simplest way to find $\{ r_n\}$ is to use a discrete unit impulse in the right-hand-side of \eqref{eq:2.100}. This unit impulse at discrete time zero is the sequence $\{ \delta_{0,n}\}$ (written in terms of Kronecker deltas) and its $\zeta$ transform is $\mathrm D_0(\zeta)\equiv 1$. Clearly this is a way of formally obtaning $\mathrm R(\zeta)$ by plugging $\mathrm H(\zeta)\equiv 1$ as data in \eqref{eq:2.103} or as right-hand-side of \eqref{eq:2.102}. Before we do it, let us just try to be more precise: in this case the impulse is $\{\delta_{0,n} I\}$ where $I$ is the identity operator in the space $Y$ and $\{r_n\}$ will be a sequence in $\mathcal B(Y,X)$. The sequence $\{ r_n\}$ (or $\mathrm R(\zeta)=\mathrm F(\zeta)^{-1}$) can be found with a simple recurrence
\[
r_0=f_0^{-1}, \qquad r_n=-f_0^{-1} \big(\sum_{m=1}^n f_m r_{n-m}\big), \quad n\ge 1.
\]
This is the discrete transfer function for the solution operator associated to the convolution equation \eqref{eq:2.100}.

\paragraph{Throwing analyticity into the mix.} Let us go further. In principle
\begin{equation}\label{eq:2.104}
\mathrm F(\zeta)=\sum_{n=0}^\infty f_n \zeta^n
\end{equation}
is just a formal series. If 
\[
\| f_n\| \le R^n,
\]
then the power series \eqref{eq:2.104} converges in the space $\mathcal B(X,Y)$ for $|\zeta|< R^{-1}$. It thus defines an analytic function from the disk centered at the origin with radius $1/R$ to $\mathcal B(X,Y)$.
Let us write \eqref{eq:2.104} in a slightly different way:
\[
\mathrm F(\zeta)=f_0\big( \sum_{n=0}^\infty f_0^{-1} f_n \zeta^n\big),
\]
so that we can focus on operators $\{ f_0^{-1}f_n\}\subset \mathcal B(X,X)$. We then worry about invertibility of the operator
\begin{equation}
I_X+\mathrm Q(\zeta), \qquad \mathrm Q(\zeta):=\sum_{n=1}^\infty f_0^{-1} f_n \zeta^n.
\end{equation}
A back of the envelope calculation yields
\[
\Big\| \sum_{n=1}^\infty f_0^{-1} f_n \zeta^n\Big\| \le \| f_0\|^{-1} \frac{1}{1-|\zeta| R},\qquad |\zeta|< 1/R.
\]
If we are lucky enough to have $\| f_0^{-1}\|< 1$, then, values of $\zeta$ satisfying
\[
|\zeta|< \frac{1-\| f_0^{-1}\|}{R},
\]
provide $\|\mathrm Q(\zeta)\|<1$ and therefore $I_X+\mathrm Q(\zeta)$ is invertible. In this case, the analytic function $\mathrm F(\zeta)$ is invertible in a disk around $\zeta=0$ and 
\[
\mathrm F(\zeta)^{-1}=\sum_{n=0}^\infty r_n \zeta^n
\]
is also a convergent power series corresponding to the discrete transfer function. Note that $\mathrm F(\zeta)^{-1}$ means exactly what it is written: it is the inverse operator for $\mathrm F(\zeta)$ for a given value of $\zeta$.

\subsection*{Credits}

The Convolution Quadrature method and its associated discretized operational calculus were created by Christian Lubich in the late eighties \cite{Lubich:1988}, as a surprising computational extension of some classical ideas by Liouville. Lubich's introduction is the one shown in Section \ref{sec:2.5}. 
A key paper in what respects to CQ applied to the wave equation is \cite{Lubich:1994}, which is the first reference to the use of CQ for a wave propagation problem. Any introduction to the numerical analysis of discretization methods for Ordinary Differential Equations contains basic material on multistep methods, and more specifically on the BDF formulas. Even if we have not seen it yet, A-stability will play an important role in how multistep methods become CQ solvers for hyperbolic problems. To learn about A-stability, see \cite{HaWa:2010}. The language of $\zeta$-series is common to the analysis of ODEs, even if sometimes is just implicitly used. A highly readable presentation is given in Henrici's classic introduction to applied complex analysis \cite{Henrici:1988}. The point of view of substituting a continuous symbol $\mathrm F(s)$ by a discrete symbol $\mathrm F(s_\kappa)$ is more or less implicit to \cite{Lubich:1994}. It was made more apparent in Antonio Laliena's thesis \cite{Laliena:2011} and it is the center of the theory for multistep CQ in \cite{Sayas:2013}. 


\section{Implementation}\label{sec:3}

In this section we are going to show how to compute discrete convolutions and convolution equations with CQ. 

\subsection{The Discrete Fourier Transform}

\paragraph{DFT and IDFT.} Given a vector $\mathbf x:=(x_0,\ldots,x_M)\in \mathbb C^{M+1}$, we consider the vector $\widehat{\mathbf x}\in \mathbb C^{M+1}$ given by
\begin{framed}
\begin{equation}\label{eq:3.1}
\widehat x_\ell:=\sum_{n=0}^M x_n \zeta_{M+1}^{-\ell n}, \quad \ell=0,\ldots,M+1,\qquad \mbox{where } \qquad \zeta_{M+1}:=e^{\frac{2\pi\imath}{M+1}}.
\end{equation}
\end{framed}
The expression \eqref{eq:3.1} is the well known definition of the Discrete Fourier Transform (DFT). The Fast Fourier Transform (FFT) is a very optimized algorithm to compute exactly \eqref{eq:3.1} when $M+1=2^p$ for an integer $p$. The Inverse Discrete Fourier Transform (IDFT) is given by
\begin{equation}\label{eq:3.2}
x_n:=\frac1{M+1} \sum_{\ell=0}^M \widehat x_\ell \zeta^{\ell n}_{M+1}, \qquad n=0,\ldots,M.
\end{equation}
The DFT and the IDFT are inverse operators. Moreover, the IDFT can be computed exactly as the DFT: first conjugate the vector $\widehat{\mathbf x}$, then apply the DFT, finally conjugate the result again and divide by $M+1$. This implies that any fast algorithm for the  DFT can also be applied to the IDFT.

\paragraph{The DFT and periodic discrete convolutions.} Let $\mathbf x,\mathbf y\in \mathbb C^{M+1}$. We define their periodic discrete convolution $\mathbf x*_{\mathrm{per}}\mathbf y\in \mathbb C^{M+1}$ with the formula
\begin{equation}\label{eq:3.3}
(\mathbf x *_{\mathrm{per}} \mathbf y)_n:=\sum_{m=0}^n x_m y_{n-m}+\sum_{m=n+1}^M x_m y_{M+1+n-m}, \qquad n=0,\ldots,M.
\end{equation}
The formula in the definition \eqref{eq:3.3} is not very inspiring. It is much easier to think that $\mathbf x,\mathbf y\in \ell^0(\mathbb Z)$ are sequences indexed by $n\in \mathbb Z$ that are $(M+1)$-periodic. Then
\begin{equation}
(\mathbf x *_{\mathrm{per}} \mathbf y)_n=\sum_{m=0}^M x_m y_{n-m} \qquad n\in \mathbb Z
\end{equation}
is also an $(M+1)$-periodic sequence coinciding with \eqref{eq:3.3} in its main entries, indexed from $n=0$ to $n=M$. The DFT diagonalizes periodic convolutions in the following sense:
\begin{equation}\label{eq:3.5}
\widehat{(\mathbf x *_{\mathrm{per}} \mathbf y)}_\ell=\widehat x_\ell \widehat y_\ell.
\end{equation}
Therefore, a fast way of computing periodic convolutions is to take the DFT of the two vectors, multiply them component by component, and then take the IDFT of the result.

\paragraph{Discrete convolutions made periodic.} Let now $\mathbf x,\mathbf y$ be causal sequences. We are interested in computing the first $N$ components of their discrete convolution
\begin{equation}\label{eq:3.6}
(\mathbf x*\mathbf y)_n=\sum_{m=0}^n x_m y_{n-m}, \qquad n=0,\ldots,N.
\end{equation}
This formula obviously involves only the vectors $(x_0,\ldots,x_N), (y_0,\ldots,y_N)\in \mathbb C^{N+1}$. Let then
\begin{equation}\label{eq:3.7}
\mathbf x^{\mathrm{ext}}:=(x_0,\ldots,x_N,\underbrace{0,\ldots,0}_{N+1}),\quad
\mathbf y^{\mathrm{ext}}:=(y_0,\ldots,y_N,\underbrace{0,\ldots,0}_{N+1})\in \mathbb C^{2N+2}
\end{equation}
be the result of cutting the vectors to the components that are needed for \eqref{eq:3.6} and then extending them with zeros, thus doubling the number of components of the vectors. A simple argument shows that
\begin{equation}\label{eq:3.8}
(\mathbf x*\mathbf y)_n=(\mathbf x^{\mathrm{ext}}*_{\mathrm{per}}\mathbf y^{\mathrm{ext}})_n, \qquad n=0,\ldots,N.
\end{equation}
Formulas \eqref{eq:3.7} and \eqref{eq:3.8} give an algorithm to compute the beginning of a discrete convolution of sequences \eqref{eq:3.6}.

\begin{framed}
\paragraph{Algorithm 3.I (discrete convolutions by the DFT).} We aim to compute
\[
(\mathbf x*\mathbf y)_n=\sum_{m=0}^n x_m y_{n-m}, \qquad n=0,\ldots,N.
\]
\begin{itemize}
\item[(a)] Keep $N+1$ components of the sequences $\mathbf x$ and $\mathbf y$ and extend them with $N+1$ zeros at the end \eqref{eq:3.7}.
\item[(b)] Take the DFT of the extended vectors.
\item[(c)] Multiply the resulting vectors component by component
\item[(d)] Take the IDFT of the result of (c)
\item[(e)] Keep only the first $N+1$ components of the result of (d).
\end{itemize}
\end{framed}

\paragraph{Symmetry arguments.} It will be common in our computations that roughly half of the coefficients of a vector will be conjugated from the other half. Before we meet this, let us introduce some notation. Assume that $\mathbb X$ is a Banach space where we can conjugate. Let $\mathbf x=(x_0,\ldots,x_N)\in \mathbb X^{N+1}$ be a vector of elements of $\mathbb X$. If
\[
x_{N+1-\ell}=\overline x_\ell \qquad  \ell=1,\ldots,N,
\]
we will say that the vector $\mathbf x$ is Hermitian (note that this is not a standard definition). This is equivalent to the following construction: extend first $\mathbf x$ to a sequence $x_n$ with $n\in \mathbb Z$ which is $(N+1)$-periodic; then $\mathbf x$ is Hermitian whenever $x_\ell=\overline x_\ell$ for all $\ell\ge 1$. In later examples we will also have that $\overline x_0=x_0$, which would yield a better definition of Hermiticity. A clear example of a Hermitian vector is the DFT of a vector such that $\overline x_n=x_n$ for all $n$. In forthcoming algorithms to compute components $x_\ell$ for $\ell=0,\ldots,N$, we will say that the algorithm can be symmetrized when we know in advance that $\mathbf x$ is Hermitian. Then we will:
\begin{itemize}
\item Compute $x_\ell$ for $\ell=0,\ldots,\lfloor\frac{N+1}2\rfloor$.
\item Copy the missing components $x_{N+1-\ell}=\overline x_\ell$ for $\ell=1,\ldots\lfloor\frac{N}2\rfloor.$
\end{itemize}

\subsection{Computation of CQ weights}\label{sec:3.2}

\paragraph{Our next goal.} We now prepare the way to compute some of the weights of the CQ process:
\begin{equation}\label{eq:3.9}
\omega_n^{\mathrm F}(\kappa):=\frac1{n!} \frac{\mathrm d^n}{\mathrm d\zeta^n} 
\left(\mathrm F\left(\frac1\kappa \delta(\zeta)\right)\right)\Big|_{\zeta=0}, \qquad n=0,\ldots,N.
\end{equation}
One of the coefficients is straightfoward:
\begin{equation}\label{eq:3.10}
\omega_0^{\mathrm F}(\kappa)=\mathrm F(\smallfrac1\kappa\delta(0)).
\end{equation}
For the other ones, we are going to use Cauchy's Formula for the computation of derivatives of an analytic function, using a contour
\[
C_R:=\{ \zeta\in \mathbb C\,:\, |\zeta|=R\}
\]
for some $R<1$ that we will later determine. Since the resulting integral can be written as the integral of a periodic function, the best way to approximate it will be with the trapezoidal rule. We will use as many points in the trapezoidal rule as coefficients we want to compute in \eqref{eq:3.9}. The process is synthesized in the following formulas:
\begin{alignat*}{6}
\omega_n^{\mathrm F}(\kappa) 
	&=\frac1{2\pi\imath}\oint_{C_R} \zeta^{-n-1} \mathrm F(\smallfrac1\kappa\delta(\zeta))\mathrm d\zeta
		&\qquad &\mbox{(Cauchy formula)}\\
	&=R^{-n} \int_0^1 e^{-2\pi\imath n\theta}
		\mathrm F(\smallfrac1\kappa\delta(R e^{2\pi\imath\theta}))\mathrm d\zeta
		&&\mbox{(parametrization $\zeta=Re^{2\pi\imath\theta}$)}\\
	&\approx \frac{R^{-n}}{N+1}\sum_{\ell=0}^N \zeta_{N+1}^{-n \ell} 
		\mathrm F(\smallfrac1\kappa \delta(R \zeta_{N+1}^\ell))
		&&\mbox{(trapezoidal rule: $\zeta_{N+1}^\ell=e^{\frac{2\pi\imath \ell}{N+1}}$)}\\
	& =\frac{R^{-n}}{N+1}\sum_{\ell=0}^N \zeta_{N+1}^{n \ell} 
		\mathrm F(\smallfrac1\kappa \delta(R \zeta_{N+1}^{-\ell})).
	 	&&\mbox{(reindexing)}
\end{alignat*}
This formula gives us an algorithm to compute \eqref{eq:3.9}. 

\begin{framed}\noindent
The formula
\[
\omega_n^{\mathrm F}(\kappa) \approx \frac{R^{-n}}{N+1}\sum_{\ell=0}^N \zeta_{N+1}^{n \ell} 
		\widehat{\mathrm F}_\ell,
		\quad \mbox{where}\quad
		\widehat{\mathrm F}_\ell:=\mathrm F(\smallfrac1\kappa \delta(R\zeta_{N+1}^{-\ell})),
\]
will be chosen as an approximation of the CQ coefficients. The chosen value of the radius of the integration path is
\begin{equation}\label{eq:R}
R=\epsilon^{\frac1{2(N+1)}},
\end{equation}
where $\epsilon$ is proportional to the machine epsilon. If $\mathrm F(\overline s)=\overline{\mathrm F(s)}$ (see the end of Section \ref{sec:1.3}) and $\delta(\overline\zeta)=\overline{\delta(\zeta)}$, then the vector of evaluations $\widehat{\mathrm F}_\ell$ is Hermitian and its computation can be reduced by symmetry.
\end{framed}

\begin{framed}
\paragraph{Algorithm 3.II (computation of CQ coefficients).} In some of the algorithms to follow we will write {\bf (Par+sym)} when the computation in that step can be reduced by symmetry and parallelized.
\begin{itemize}
\item[(a)] {\bf (Par+sym)} Evaluate
\[
\widehat{\mathrm F}_\ell:=\mathrm F(\smallfrac1\kappa \delta(R\zeta_{N+1}^{-\ell})), \qquad \ell=0,\ldots,N.
\]
\item[(b)] Apply the IDFT
\[
\mathrm F_n:=\frac1{N+1}\sum_{\ell=0}^N \widehat{\mathrm F}_{\ell}\zeta_{N+1}^{n\ell},
\qquad n=0,\ldots,N.
\]
When $\mathrm F$ is matrix valued, this has to be done component by component.
\item[(c)] Scale the result
\[
\omega_n^{\mathrm F}(\kappa)=R^{-n}\mathrm F_n, \qquad n=0,\ldots,N.
\]
At this step it is also possible to substitute the approximation of $\omega_0^{\mathrm F}(\kappa)$ by its exact value $\omega_0^{\mathrm F}(\kappa)=\mathrm F(\frac1\kappa\delta(0))$.
\end{itemize}
\end{framed}

\paragraph{Forward convolutions and convolution equations.} To compute now
\begin{equation}\label{eq:3.11}
\sum_{m=0}^n \omega_m^{\mathrm F}(\kappa) g_{n-m} \qquad n=0,\ldots,N,
\end{equation}
we can either apply a naive implementation  formula \eqref{eq:3.11} or use the zero padding strategy of Algorithm 3.I. Note that one of the first things we will do in Algorithm 3.I will be to take the DFT of the sequence of coefficients $\omega_n^{\mathrm F}(\kappa)$. In Section \ref{sec:3.3}, we will explore a more direct way of doing this. A similar problem arises in the solution of convolution equations
\begin{equation}\label{eq:3.12}
\omega_0^{\mathrm F}(\kappa)g_n=h_n-\sum_{m=1}^n \omega_m^{\mathrm F}(\kappa) g_{n-m}.
\end{equation}
The right-hand-side of \eqref{eq:3.12} contains a discrete convolution that can be evaluated using Algorithm 3.I by artificially including a term $\omega_0^{\mathrm F}(\kappa)=0$. Solving the CQ equations by the forward substitution scheme \eqref{eq:3.12} is a particular case of what people in the know call {\bf marching-on-in-time} (MoT) schemes.

\subsection{All-steps-at-once CQ computation}\label{sec:3.3}

\paragraph{Development of the algorithm.} Our next goal is the computation of
\[
\sum_{m=0}^n \omega_{n-m}^{\mathrm F}(\kappa) g_{m}, \qquad n=0,\ldots,N,
\]
assuming that we have approximated the CQ coefficients by
\[
\omega_n^{\mathrm F}(\kappa)\approx
\frac{R^{-n}}{N+1}\sum_{\ell=0}^N \widehat{\mathrm F}_\ell\zeta_{N+1}^{\ell n}, \qquad
\widehat{\mathrm F}_\ell:=\mathrm F(\smallfrac1\kappa\delta(R \zeta^{-\ell}_{N+1})).
\]
We are actually going to compute something slightly different. Because of Cauchy's Theorem
\[
\omega_n^{\mathrm F}(\kappa)=
\oint_{C_R}\zeta^{-n+1} \mathrm F(\zeta)\mathrm d\zeta = 0, \qquad \forall n \le -1,
\]
we are going to define (using the same idea as in Section \ref{sec:3.2})
\[
\widetilde\omega_n^{\mathrm F}(\kappa):=\frac{R^{-n}}{N+1}\sum_{\ell=0}^N \widehat{\mathrm F}_\ell\zeta_{N+1}^{\ell n} \qquad n\in \mathbb Z,
\]
knowing that the approximated coefficients for $n\le -1$ converge to zero. Then
\begin{alignat*}{6}
u_n 
	&:= \sum_{m=0}^n \omega_{n-m}^{\mathrm F}(\kappa)g_m
		&\qquad& \mbox{(exact convolution)}\\
	& = \sum_{m=0}^N \omega_{n-m}^{\mathrm F}(\kappa) g_m 
		& & \mbox{($\omega_n^{\mathrm F}(\kappa)=0$ for $n\le -1$)}\\
	& \approx \sum_{m=0}^N \left(\frac{R^{m-n}}{N+1}
		\sum_{\ell=0}^N \widehat{\mathrm F}_\ell \zeta_{N+1}^{\ell(n-m)}\right) g_m
		& & \mbox{($\omega_n^{\mathrm F}(\kappa)\approx \widetilde\omega_n^{\mathrm F}(\kappa)$)}\\
	& = R^{-n} \left( \frac1{N+1} \sum_{\ell=0}^N \widehat{\mathrm F}_\ell
		\left(\sum_{m=0}^N R^m g_m \zeta_{N+1}^{-m\ell}\right) \zeta_{N+1}^{\ell n}\right).
\end{alignat*}
These formulas give the next algorithm.

\begin{framed}
\paragraph{Algorithm 3.III (all-steps-at-once forward convolution).} We sample the input at given discrete times
\[
g_n:=g(t_n), \qquad n=0,\ldots,N
\]
and aim to compute an approximation of
\[
u_n=\sum_{m=0}^n \omega_{n-m}^{\mathrm F}(\kappa)g_m, \qquad n=0,\ldots,N.
\]
\begin{itemize}
\item[(a)] Scale the data
\[
h_m:=R^m g_m, \qquad m=0,\ldots,N.
\]
\item[(b)] Compute the DFT of the previous sequence
\[
\widehat h_\ell:=\sum_{m=0}^N h_m \zeta_{N+1}^{-m \ell}, \qquad \ell=0,\ldots,N.
\]
\item[(c)] {\bf (Par+sym)} Apply the transfer functions in the transformed domain
\[
\widehat v_\ell:=\widehat{\mathrm F}_\ell \widehat h_\ell, \qquad 
\widehat{\mathrm F}_\ell:=\mathrm F(\smallfrac1\kappa\delta(R \zeta_{N+1}^{-\ell})).
\]
\item[(d)] Compute the IDFT of the result of (c)
\[
v_n:=\frac1{N+1}\sum_{\ell=0}^N \widehat v_\ell \zeta_{N+1}^{\ell n}, \qquad n=0,\ldots,N.
\]
\item[(e)] Scale the result back
\[
u_n=R^{-n} v_n.
\]
\end{itemize}
\end{framed}

\paragraph{Modification for convolution equations.} When the goal is to solve a convolution equation
\[
f*g=h,
\]
we can use exactly the same ideas applied to the operator $\mathrm F^{-1}$. To compute
\[
g_n=\sum_{m=0}^n \omega_m^{\mathrm F^{-1}}(\kappa) h(t_{n-m}), \qquad n=0,\ldots,N,
\]
we approximate the coefficients
\[
\omega_n^{\mathrm F^{-1}}(\kappa) \approx \frac{R^{-n}}{N+1} \sum_{\ell=0}^N \widehat{\mathrm F}_\ell^{-1} \zeta_{N+1}^{\ell n}, \qquad 
\widehat{\mathrm F}_\ell^{-1}=\mathrm F(\smallfrac1\kappa\delta(R\zeta_{N+1}^{-\ell}))^{-1}.
\]
The final expression is
\[
g_n\approx R^{-n} \left( \frac1{N+1} \sum_{\ell=0}^N \widehat{\mathrm F}_\ell^{-1}
		\left(\sum_{m=0}^N R^m h_m \zeta_{N+1}^{-m\ell}\right) \zeta_{N+1}^{\ell n}\right).
\]
from where it is clear that the inverses $\widehat{\mathrm F}_\ell^{-1}$ do not need to be computed. Instead, some linear systems will be solved.

\begin{framed}
\paragraph{Algorithm 3.IV (all-steps-at-once convolution equation).} We sample the input at given discrete times
\[
h_n:=h(t_n), \qquad n=0,\ldots,N
\]
and aim to compute an approximation of the solution of
\[
\sum_{m=0}^n \omega_{n-m}^{\mathrm F}(\kappa)g_m=h_n, \qquad n=0,\ldots,N.
\]
\begin{itemize}
\item[(a)] Scale the data
\[
v_m:=R^m h_m, \qquad m=0,\ldots,N.
\]
\item[(b)] Compute the DFT of the previous sequence $\widehat v_\ell$, $\ell=0,\ldots,N$.
\item[(c)] {\bf (Par+sym)} Solve equations in the transformed domain
\[
\widehat{\mathrm F}_\ell\widehat w_\ell:= \widehat v_\ell, \qquad 
\widehat {\mathrm F}_\ell:=\mathrm F(\smallfrac1\kappa\delta(R \zeta_{N+1}^{-\ell})).
\]
\item[(d)] Compute the IDFT of the result of (c) $w_\ell$, $\ell=0,\ldots,N$.
\item[(e)] Scale the result back
\[
g_n=R^{-n} w_n, \qquad n=0,\ldots,N.
\]
\end{itemize}
\end{framed}

\subsection{Computing pieces of a discrete convolution}\label{sec:3.4}

\paragraph{Pieces of a discrete convolution.} Our next goal is the computation of quantities
\begin{equation}\label{eq:3.13}
g_n:=\sum_{m=0}^Q \widetilde\omega_{n-m}^{\mathrm F}(\kappa) u_m, \qquad n=Q+1,\ldots,M,
\end{equation}
where
\begin{equation}\label{eq:3.14}
\widetilde\omega_n^{\mathrm F}(\kappa)=\frac{R^{-n}}{N+1}\sum_{\ell=0}^N \widehat{\mathrm F}_\ell \zeta_{N+1}^{n \ell},
\qquad \widehat{\mathrm F}_\ell:=\mathrm F(\smallfrac1\kappa \delta(R\zeta_{N+1}^{-\ell})),
\end{equation}
for a given $N\ge M$, which can be chosen as $N=2 M$ for example. Note that \eqref{eq:3.13} uses $u_0,\ldots,u_Q$ and $\widetilde\omega_n^{\mathrm F}(\kappa)$ for $n=1,\ldots,M$.
We then proceed as follows:
\begin{alignat*}{6}
\widetilde g_k 
	&:= g_{k+Q+1} 
		\qquad\qquad\qquad\qquad\qquad (k=0,\ldots,M-Q-1)\\
	&=\sum_{m=0}^Q \left(\frac{R^{m-k-Q-1}}{N+1}
			\sum_{\ell=0}^N\widehat{\mathrm F}_\ell\zeta_{N+1}^{\ell(k+Q+1-m)}\right) u_m \\
	&=R^{-k-Q-1} \left(\frac1{N+1}\sum_{\ell=0}^N \zeta_{N+1}^{\ell(Q+1)}
			\widehat{\mathrm F}_\ell \left( \sum_{m=0}^Q \zeta_{N+1}^{-\ell m} R^m u_m\right) 
			\zeta_{N+1}^{\ell k}\right)\\
	&=R^{-k-Q-1} \left(\frac1{N+1}\sum_{\ell=0}^N \left(\zeta_{N+1}^{\ell(Q+1)}
			\widehat{\mathrm F}_\ell \left( \sum_{m=0}^N \zeta_{N+1}^{-\ell m} w_m\right)\right) 
			\zeta_{N+1}^{\ell k}\right),
\end{alignat*}
where
\begin{equation}\label{eq:3.15}
w_m:=\left\{ \begin{array}{ll} R^m u_m, & 0\le m\le Q,\\ 0, & Q+1\le m\le N.\end{array}\right.
\end{equation}

\begin{framed}
\paragraph{Algorithm 3.V (computation of a piece of a convolution).} Our goal is to compute \eqref{eq:3.13}. The value $N$ in \eqref{eq:3.14} is given as a parameter. The data $u_m$, $m=0,\ldots,Q$ has already been sampled.
\begin{itemize}
\item[(a)] Scale the data and add zeros following \eqref{eq:3.15}.
\item[(b)] Compute the DFT of the vector in (a) $\widehat w_\ell$, $\ell=0,\ldots,N$.
\item[(c)] {\bf (Par+sym)} Apply the operators
\[
\widehat h_\ell:=\zeta_{N+1}^{\ell(Q+1)} \widehat{\mathrm F}_\ell\widehat w_\ell, \qquad \ell=0,\ldots,N.
\]
\item[(d)] Compute the IDFT of the sequence in (c) $h_k$, $k=0,\ldots,N$.
\item[(e)] Scale and chop the resulting sequence
\[
\widetilde g_k:=R^{-k-Q-1} h_k, \qquad k=0,\ldots,M-Q-1.
\]
\item[(f)] Change indices
\[
g_n:=\widetilde g_{n-Q-1}, \qquad n=Q+1,\ldots,M.
\]
\end{itemize}
\end{framed}

If we skip step (f), what we have computed is 
\[
\widetilde g_k:=\sum_{m=0}^Q \widetilde\omega_{k+Q+1-m}^{\mathrm F}(\kappa) u_m, \qquad k=0,\ldots,M-Q-1, 
\]
which can be understood as the formal product with a piece of a Toeplitz matrix
\[
\left[\begin{array}{c}
	\widetilde g_0 \\ \widetilde g_1 \\ \vdots \\ \widetilde g_{M-Q-1}
\end{array}\right]=
\left[ \begin{array}{cccc}
	\widetilde\omega_{Q+1}^{\mathrm F}(\kappa) & \widetilde\omega_{Q}^{\mathrm F}(\kappa) & \ldots &
	\widetilde\omega_{1}^{\mathrm F}(\kappa)\\
	\widetilde\omega_{Q+2}^{\mathrm F}(\kappa) & \widetilde\omega_{Q+1}^{\mathrm F}(\kappa) & \ldots &
	\widetilde\omega_{2}^{\mathrm F}(\kappa)\\
	\widetilde\omega_{Q+3}^{\mathrm F}(\kappa) & \widetilde\omega_{Q+2}^{\mathrm F}(\kappa) & \ldots &
	\widetilde\omega_{3}^{\mathrm F}(\kappa)\\
	\vdots & \vdots & \ddots & \vdots\\
	\widetilde\omega_{M}^{\mathrm F}(\kappa) & \widetilde\omega_{M-1}^{\mathrm F}(\kappa) & \ldots &
	\widetilde\omega_{M-Q}^{\mathrm F}(\kappa)
\end{array}\right]
\left[\begin{array}{c} 
	u_0 \\ u_1 \\ \vdots \\ u_Q
\end{array}\right].
\]

\begin{figure}[htb]
\begin{center}
\includegraphics[width=7cm]{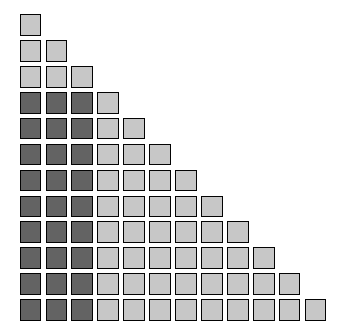}
\includegraphics[width=7cm]{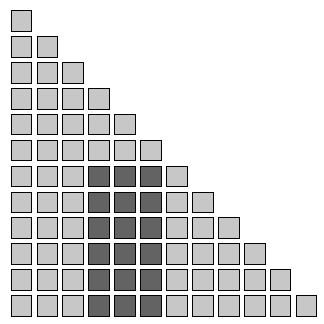}
\end{center}
\caption{The leftmost figure shows a typical convolution piece as computed by Algorithm 3.V. Three elements of the sequence of data are input and nine are output. The figure on the right is a simple variation. By a shifting of indices in the input data a different block of the Toeplitz matrix can be accessed. Note that each of the elements is a matrix (or an operator if we have not discretized the original operator).}
\end{figure}

\subsection{Recursive strategies}

\paragraph{Small systems.} We are going to next include a new parameter $\underline{n}$ corresponding to sizes of lower triangular systems (convolution equations) that we want to solve directly:
\begin{equation}\label{eq:3.16}
\sum_{m=0}^n \omega_{n-m} g_m=u_n, \qquad n=0,\ldots,\underline n.
\end{equation}
From now on, we are going to ignore where the coefficients $\omega_n$ are coming from. They are operators, $\omega_0$ being invertible, or simply square matrices. We can write in matrix form with help of the matrix
\[
\Omega_{\underline n}:=\left[\begin{array}{cccc}
	\omega_0 \\
	\omega_1 & \omega_0 \\
	\vdots & \ddots & \ddots \\
	\omega_{\underline n-1} & \ldots & \omega_1 & \omega_0
\end{array}\right].
\]

\paragraph{Look-ahead method.} With this strategy, we solve groups of $\underline n$ equations and then subtract the contribution of the newly computed unknowns from the right-hand-side of the system. The method can then be understood as block recursive forward substitution.

\begin{framed}
\paragraph{Algorithm 3.VI (look-ahead solution of CQ equations).} The goal is the solution of
\[
\sum_{m=0}^n \omega_{n-m} g_m=u_n, \qquad n=0,\ldots,N.
\]
We give a parameter $\underline n$ for the size of the small systems and break the list of indices in the form
\[
\underbrace{0,\ldots,\underline n-1}_{b(1)},
\underbrace{\underline n,\ldots,2\underline n-1}_{b(2)},\ldots,
\underbrace{(k-1)\underline n,\ldots,k\underline n-1}_{b(k)},
\underbrace{k\underline n,\ldots,N}_{\mbox{remainder}},
\]
so that
\[
b(i)=\{ (i-1)\underline n,\ldots,i\underline n-1\} \qquad i=1,\ldots, k, \quad \mbox{where }k =\lfloor N/\underline n\rfloor,
\]
is a generic block. The algorithm loops in $i=1,\ldots,k$ in the following form:
\begin{itemize}
\item[(a)] Make a copy of a block of data $\boldsymbol v=\boldsymbol u_{b(i)}$ Solve
\[
\Omega_{\underline n}\boldsymbol h=\boldsymbol v,
\]
and copy the result $\boldsymbol q_{b(i)}=\boldsymbol h$. 
\item[(b)] Compute the following piece of convolution
\[
r_n=\sum_{m=0}^{\underline n-1} \omega_{n-m} h_m, \qquad n=\underline n,\ldots,N-(i-1)\underline n.
\]
\item[(c)] Correct the right-hand-side
\[
u_n=u_n-r_{n-(i-1) \underline n} \qquad n=i \underline n,\ldots, N.
\]
\end{itemize}
Finally, after the loop is finished, we have to solve for the tail (because of $N+1$ not being a multiple of $\underline n$):
\[
\sum_{m=0}^n \omega_{n-m} g_{k\underline n+m}=u_{k\underline n+m}, \qquad m=0,\ldots,N-k\underline n.
\]
\end{framed}

\begin{figure}[htb]
\begin{center}
\includegraphics[width=7cm]{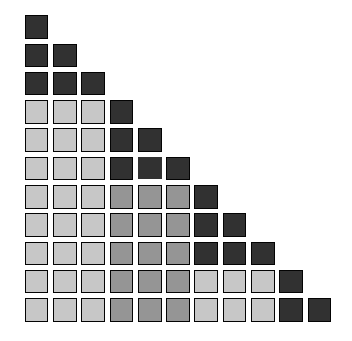}
\includegraphics[width=7cm]{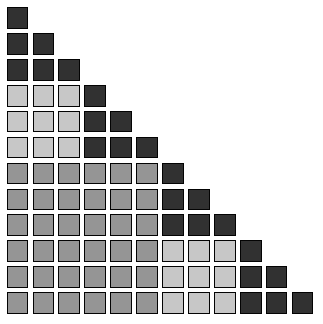}
\end{center}
\caption{This cartoon shows the basic ideas of a look-ahead and a recursive strategy side-by-side. The black blocks of equations are inverted directly by forward substitution. In the look-ahead strategy, once a group of unknowns has been computed, its influence on the right-hand-side is immediately taken into account. The recursive idea breaks down the triangular system into smaller triangular systems joined with square blocks.}
\label{fig:2}
\end{figure}

\paragraph{Recursive methods.} Another option to deal with the CQ equations is the development of a recursive algorithm, based on the break down of a triangular system into two pieces of the same size and the square block that interconnects them. We will not develop this algorithm any further. A pictorial representation, side by side with the look-ahead algorithm, is given in Figure \ref{fig:2}.

\subsection*{Credits}

An excellent introduction to the algorithms for CQ is given in \cite{BaSc:2012}. The idea of using the trapezoidal rule on a Cauchy integral representation in order to compute the CQ coefficients appears already in the original papers by Lubich \cite{Lubich:1988, Lubich:1994}. The simultaneous computation of all time steps is introduced in a paper by Lehel Banjai and Stefan Sauter \cite{BaSa:2008} for the specific problem of the wave equation (see Section \ref{sec:4}), based on algorithms developed in \cite{HaLuSc:1985}. The choice of $R$ in the integration contour \eqref{eq:R} is justified in \cite{Lubich:1988b}.


\section{Integral equations for waves}\label{sec:4}

In this section we are going to explain how the ideas on Laplace transforms, convolution operators and equations, and CQ discretizations can be used for a problem of scattering of acoustic waves by an obstacle. This section is considerably heavier with respects to Sobolev spaces. The reader should have at least some working knowledge of the Sobolev space $H^1(\mathcal O)$ for an open set $\mathcal O$ and of the spaces $H^{\pm 1/2}(\partial\mathcal O)$ for an open set $\mathcal O$ which lies on one side of its Lipschitz boundary. 

\begin{framed}
\noindent
The space
\[
H^1_\Delta(\mathcal O):=\{ u\in H^1(\mathcal O)\,:\, \Delta u \in L^2(\mathcal O)\},
\]
endowed with its natural norm, will play a key role in this section. Note that from this space we can define three bounded operators:
\[
\gamma : H^1_\Delta(\mathcal O)\to H^{1/2}(\partial\mathcal O),\qquad
\partial_\nu : H^1_\Delta(\mathcal O) \to H^{-1/2}(\mathcal O),\qquad
\Delta : H^1_\Delta(\mathcal O)\to L^2(\mathcal O),
\]
the first two of which correspond to the trace and the normal derivative.
\end{framed}

\subsection{A scattering problem}\label{sec:4.1}

\paragraph{Incident waves.} Let $\Omega_-\subset \mathbb R^3$ be a connected domain, lying on one side of its connected Lipschitz boundary $\Gamma$. Let $\Omega_+:=\mathbb R^3\setminus\overline{\Omega_-}$. At time $t=0$ an incident wave approaches the obstacle $\Gamma$. The incident wave can be described as a solution of the wave equation in free space. For instance, a plane wave
\begin{equation}\label{eq:4.1}
u^{\mathrm{inc}}(\mathbf z,t):=\psi(c (t-t_{\mathrm{lag}})-\mathbf z\cdot\mathbf d), \qquad |\mathbf d|=1, \quad c>0,
\end{equation}
solves
\begin{equation}\label{eq:4.2}
c^{-2}\ddot u^{\mathrm{inc}}=\Delta u^{\mathrm{inc}} \qquad \mbox{ in $\mathbb R^3\times (-\infty,\infty)$}.
\end{equation}
Equation \eqref{eq:4.2} is satisfied in a classical way when the signal is smooth enough, for instance, when $\psi\in \mathcal C^2(\mathbb R)$. When $\psi$ is just continuous, and even when it is only locally integrable, \eqref{eq:4.2} can be understood in a weak sense of distributions in four variables (an approach we have decided not to take), or even more surprisingly in a sense of distributions of the time variable, with values in some functions spaces. We will discuss this later. We need to make some assumptions on the plane wave \eqref{eq:4.1} in order for the scattering process to be physically meaningful. For instance, we can assume that $\psi$ is causal, $\Omega_-\subset B(\mathbf 0;R):=\{ \mathbf z\in \mathbb R^3\,:\, |\mathbf z|<R\}$ and $t_{\mathrm{lag}}> R/c$. This ensures that at time $t=0$, the wave, which is moving in the direction $\mathbf d$, has not reached the area of the space where the obstacle $\Omega_-$ is placed. A spherical incident wave
\begin{equation}\label{eq:4.3}
u^{\mathrm{inc}}(\mathbf z,t):=
\frac{\psi(ct-|\mathbf z-\mathbf z_{\mathrm{sc}}|)}{4\pi |\mathbf z-\mathbf z_{\mathrm{sc}}|},
\end{equation}
can also be defined for any causal $\psi$ and source point $\mathbf z_{\mathrm{sc}}\not\in \overline\Omega_-$. Instead of \eqref{eq:4.2}, this $u^{\mathrm{inc}}$ is a causal solution of an equation
\begin{equation}\label{eq:4.4}
c^{-2}\ddot u^{\mathrm{inc}}=\Delta u^{\mathrm{inc}}+f,
\end{equation}
where $f$ is formally a distribution that is supported in the source point. Let us warn the reader that, even if it is possible to give a completely rigorous  presentation of the meaning of the equation satisfied by the incident wave, we will only need its values on $\Gamma$ and we can make much simpler while rigorous arguments for those values.

\paragraph{Transient scattering by a sound-soft obstacle.} Let us then assume that we have an incident wave satisfying \eqref{eq:4.4} for all times and for a given $f$. We do not need to know $f$. Knowing that $f$ exists is enough. The total wave is a function $u^{\mathrm{tot}}$ satisfying
\begin{subequations}
\begin{equation}
c^{-2} \ddot u^{\mathrm{tot}}=\Delta u^{\mathrm{tot}}+f \qquad \mbox{in $\Omega_+\times (-\infty,\infty)$},
\end{equation}
with a boundary condition 
\begin{equation}
u^{\mathrm{tot}}=0 \qquad \mbox{on $\Gamma\times (-\infty,\infty)$},
\end{equation}
and the assumption
\begin{equation}
u^{\mathrm{tot}}-u^{\mathrm{inc}} \quad\mbox{is causal.}
\end{equation}
\end{subequations}
Instead of working with the total wave field, we will think in terms of the {\bf scattered wave field}
\[
u=u^{\mathrm{tot}}-u^{\mathrm{inc}},
\]
which satisfies a homogeneous wave equation
\begin{subequations}\label{eq:4.6}
\begin{equation}
c^{-2} \ddot u =\Delta u  \qquad \mbox{in $\Omega_+\times (-\infty,\infty)$},
\end{equation}
a non-homogeneous boundary condition 
\begin{equation}
u+u^{\mathrm{inc}}=0 \qquad \mbox{on $\Gamma\times (-\infty,\infty)$},
\end{equation}
and the assumption
\begin{equation}
u \quad\mbox{is causal.}
\end{equation}
\end{subequations}
It is interesting to note that our assumptions for the incident wave (causality of the signal in \eqref{eq:4.3}, and more complicated geometric assumptions on \eqref{eq:4.1}) were destined to ensure that the restriction of $u^{\mathrm{inc}}$ to $\Gamma$ is a causal function. Therefore, instead of \eqref{eq:4.6}, we can think of looking for a causal distribution with values in
\[
H^1_\Delta(\Omega_+):=\{ v\in H^1(\Omega_+)\,:\, \Delta v\in L^2(\Omega_+)\}.
\]
satisfying
\begin{equation}\label{eq:4.7}
c^{-2} \ddot u =\Delta u \quad \mbox{in $H^1_\Delta(\Omega_+)$}, \qquad \gamma u=\beta\quad \mbox{in $H^{1/2}(\Gamma)$}.
\end{equation}
Let us first clarify the meaning of equations \eqref{eq:4.7}. The steady state operator $\Delta:H^1_\Delta(\Omega_+)\to L^2(\Omega_+)$ can be applied to any $H^1_\Delta(\Omega_+)$-valued distribution $u$, producing an $L^2(\Omega_+)$-valued distribution $\Delta u$. Similarly, since $u$ is $H^1_\Delta(\Omega_+)$-valued, so is $\ddot u$, and therefore, we can understand $\ddot u$ as an $L^2(\Omega_+)$-valued distribution, by using the steady-state embedding operator $H^1_\Delta(\Omega_+)\to L^2(\Omega_+)$. The second part of \eqref{eq:4.7} arises from taking the trace operator $\gamma:H^1_\Delta(\Omega_+)\to H^{1/2}(\Gamma)$ and comparing the $H^{1/2}(\Gamma)$-valued distribution $\gamma u$ with a given causal $H^{1/2}(\Gamma)$-valued distribution $\beta$. At this time, it is clear what we need from the incident wave: as long as the boundary values of $u^{\mathrm{inc}}$ define a causal distribution with values in $H^{1/2}(\Gamma)$, we can give meaning to equations \eqref{eq:4.7}.

\subsection{The acoustic single layer potential}\label{sec:4.2}

\paragraph{A what-it-does definition.} We are first going to formally introduce the $d$-dimensional acoustic single layer potential and two related boundary integral operators not through mathematical expressions, but through their properties. To do this we need two jump operators
\[
\jump{\gamma u}:=\gamma^- u-\gamma^+ u, \qquad \jump{\partial_\nu u}:=\partial_\nu^+ u-\partial_\nu^- u.
\]
Both of them are well defined in $H^1_\Delta(\mathbb R^d\setminus\Gamma)$. Let now $\eta$ be a causal $H^{-1/2}(\Gamma)$-valued distribution, which will be referred to as a {\bf density}. We then consider the following wave propagation problem in free space
\begin{subequations}\label{eq:4.8}
\begin{alignat}{6}
\label{eq:4.8a}
\ddot u=\Delta u &\qquad &&\mbox{in $L^2(\mathbb R^d\setminus\Gamma)$},\\
\jump{\gamma u} = 0 &&&\mbox{in $H^{1/2}(\Gamma)$},\\
\jump{\partial_\nu u}=\eta &&& \mbox{in $H^{-1/2}(\Gamma)$}.
\end{alignat}
\end{subequations}
The mathematically savvy reader will undoubtedly find silly that we have written $L^2(\mathbb R^d\setminus\Gamma)$ instead of $L^2(\mathbb R^d)$ in \eqref{eq:4.8a}. We do it just to emphasize that we are applying the Laplacian as an operator $\Delta:H^1_\Delta(\mathbb R^d\setminus\Gamma)\to L^2(\mathbb R^d\setminus\Gamma)$. Note also that in \eqref{eq:4.8a} we consider the wave equation on both sides of the boundary of the scatterer, at speed $c=1$.  For the moment being, we are going to admit that problem \eqref{eq:4.8} is solvable, and we then (formally) denote
\begin{equation}\label{eq:4.9}
u =:  \mathcal S*\eta,\qquad
\gamma u =:\mathcal V*\eta=\gamma^\pm (\mathcal S*\eta),\qquad
\smallfrac12(\partial_\nu^-u+\partial_\nu^+u)=:\mathcal J*\eta.
\end{equation}
The convolution symbol in the definitions \eqref{eq:4.9} will be justified below. By definition
\begin{equation}\label{eq:4.10}
\partial_\nu^+ u =-\jump{\partial_\nu u}+\smallfrac12(\partial_\nu^-u+\partial_\nu^+u)=-\smallfrac12\eta+\mathcal J*\eta.
\end{equation}

\paragraph{The Laplace transform way.} The simplest way to justify problem \eqref{eq:4.8}, and the three definitions \eqref{eq:4.9} that follow from it, is to take the Laplace transform. Let then $\mathrm E=\mathcal L\{\eta\}$ and consider the problem of finding $\mathrm U(s)\in H^1_\Delta(\mathbb R^d\setminus\Gamma)$ satisfying
\begin{subequations}\label{eq:4.11}
\begin{alignat}{6}
s^2 \mathrm U(s)=\Delta \mathrm U(s) 
	&\qquad && \mbox{in $\mathbb R^d\setminus\Gamma$, (i.e., as functions of $L^2(\mathbb R^d\setminus\Gamma))$}\\
\jump{\gamma\mathrm U(s)}=0, & && \mbox{(as functions in $H^{1/2}(\Gamma)$)}\\
\jump{\partial_\nu \mathrm U(s)}=\mathrm E(s). & & & \mbox{(as functions in $H^{-1/2}(\Gamma)$)}
\end{alignat}
\end{subequations}
Equations \eqref{eq:4.11} are uniquely solvable for values $s\in \mathbb C_+$. Actually the theory of time-harmonic acoustic layer potentials gives an explicit solution to \eqref{eq:4.11}. We first introduce the fundamental solution of the Laplace resolvent operator $\Delta -s^2$,
\[
\Phi(s;\mathbf r):=
\left\{\begin{array}{ll} 
	\frac\imath4 H^{(1)}_0(-\imath s |\mathbf r|), & \mbox{when $d=2$},\\
	\ds \frac{e^{-s|\mathbf r|}}{4\pi|\mathbf r|}, &\mbox{when $d=3$}.
\end{array}\right.
\]
With it, and for a given density $\xi\in H^{-1/2}(\Gamma)$, we define the potential
\begin{subequations}\label{eq:4.12}
\begin{equation}\label{eq:4.12a}
\mathrm S(s) \xi 
:=\int_\Gamma \Phi(s;|\cdot-\mathbf y|)\xi(\mathbf y)\mathrm d\Gamma(\mathbf y) 
\,:\, \mathbb R^d\setminus\Gamma\to \mathbb C
\end{equation}
and two integral operators
\begin{alignat}{6}
\label{eq:4.12b}
\mathrm V(s) \xi
:=& \int_\Gamma \Phi(s;|\cdot-\mathbf y|)\xi(\mathbf y)\mathrm d\Gamma(\mathbf y) 
& \,:\, \Gamma\to \mathbb C,\\
\label{eq:4.12c}
\mathrm J(s) \xi 
:= & \boldsymbol\nu \cdot \int_\Gamma \nabla \Phi(s;|\cdot-\mathbf y|) \xi(\mathbf y)\mathrm d\Gamma(\mathbf y) 
& \,:\, \Gamma\to \mathbb C.
\end{alignat}
\end{subequations}
Given the overall lack of regularity, the definitions \eqref{eq:4.12} have to be taken with a grain of salt. For instance, for a given $\mathbf z\in \mathbb R^d\setminus\Gamma$, the definition \eqref{eq:4.12a} of $(\mathrm S(s)\xi)(\mathbf z)$ can always be understood as a duality product of $\xi\in H^{-1/2}(\Gamma)$ with $\Phi(s;|\mathbf z-\cdot|)\in H^{1/2}(\Gamma)$, which admits an integral expression like \eqref{eq:4.12a} when $\xi\in L^2(\Gamma)$. The definitions \eqref{eq:4.12b} and \eqref{eq:4.12c} are completely reasonable for smooth enough $\xi$ and $\Gamma$. Otherwise they need to be extended using density arguments. What really matters now is that
\begin{equation}\label{eq:4.13}
\mathrm U(s)=\mathrm S(s)\mathrm E(s), \qquad \gamma\mathrm U(s)=\mathrm V(s)\mathrm E(s),\qquad
\smallfrac12(\partial_\nu^- \mathrm U(s)+\partial_\nu^+ \mathrm U(s))=\mathrm J(s)\mathrm E(s)
\end{equation}
is the solution to \eqref{eq:4.11}. For future reference (and mimicking \eqref{eq:4.10}) we are going to care about the operator
\begin{equation}
\mathrm E(s) \longmapsto 
-\smallfrac12 \mathrm E(s)+\mathrm J(s)\mathrm E(s).
\end{equation}

\paragraph{Some Sobolev space notation.} Given an open set $\mathcal O$, we will write
\[
\| u\|_{\mathcal O}^2:=\int_{\mathcal O}|u|^2, \qquad \| u\|_{1,\mathcal O}^2:=\| u\|_{\mathcal O}^2+\| \nabla u\|_{\mathcal O}^2.
\]
The norm in the trace space $H^{1/2}(\Gamma)$ will be denoted $\|\cdot\|_{1/2,\Gamma}$, and its dual norm is
\[
\| \eta\|_{-1/2,\Gamma}:=\sup_{0\neq \phi\in H^{1/2}(\Gamma)} \frac{|\langle \eta,\phi\rangle_\Gamma|}{\| \phi\|_{1/2,\Gamma}},
\]
where $\langle \eta,\phi\rangle_\Gamma$ is the duality bracket.

\paragraph{From Laplace to time.}
It does not take very long (the Lax-Milgram theorem with some carefully crafted bounds taking care of the parameter $s$) to see that
for every $s\in \mathbb C_+$, the solution of \eqref{eq:4.11} exists, and we can bound 
\begin{subequations}\label{eq:4.15}
\begin{alignat}{6}
\| \mathrm U(s) \|_{1,\mathbb R^d\setminus\Gamma} &\le C \frac{|s|}{\sigma \underline\sigma^2} \| \mathrm E(s)\|_{-1/2,\Gamma},\\
\| \partial_\nu^{\pm} \mathrm U(s)\|_{-1/2,\Gamma} & \le C \frac{|s|^{3/2}}{\sigma\underline\sigma^{3/2}}\| \mathrm E(s)\|_{-1/2,\Gamma},
\end{alignat}
\end{subequations}
where
\[
\sigma:=\mathrm{Re}\,s, \qquad \underline\sigma:=\min\{\sigma,1\}.
\]
From \eqref{eq:4.13} and \eqref{eq:4.15} we can derive the following bounds (depending on $s$) 
\begin{subequations}\label{eq:4.16}
\begin{alignat}{6}
\| \mathrm S(s)\|_{H^{-1/2}(\Gamma) \to H^1(\mathbb R^d)}
	&\le C \frac{|s|}{\sigma\underline\sigma^2},\\
\| \mathrm V(s) \|_{H^{-1/2}(\Gamma) \to H^{1/2}(\Gamma)}
	& \le C \frac{|s|}{\sigma\underline\sigma^2},\\
\| \mathrm J(s) \|_{H^{-1/2}(\Gamma)\to H^{-1/2}(\Gamma)}
	& \le C \frac{|s|^{3/2}}{\sigma\underline\sigma^{3/2}}.
\end{alignat}
\end{subequations}
Using the results of Section \ref{sec:1.4} (see specifically conditions \eqref{eq:1.100} and \eqref{eq:1.101}) we can show that there exist causal operator valued distributions $\mathcal S$, $\mathcal V$, and $\mathcal J$, whose Laplace transforms are $\mathrm S$, $\mathrm V$, and $\mathrm J$ respectively. This long process through the Laplace domain gives a precise meaning to the time-domain convolution operators \eqref{eq:4.9} and \eqref{eq:4.10}. 

\paragraph{Time domain expressions.} Explicit expressions for the operators $\mathcal S*\eta$, $\mathcal V*\eta$, and $\mathcal J*\eta$ are available in two and three dimensions. Note that with the CQ approach, there is no need for them, since only their Laplace transform is ever used. For instance in three dimensions
\begin{equation}
(\mathcal S*\eta)(t) :=\int_\Gamma \frac{\eta(\mathbf y;t-|\cdot-\mathbf y|)}{4\pi|\cdot-\mathbf y|} \mathrm d\Gamma(\mathbf y)
: \mathbb R^3\setminus\Gamma \to \mathbb R,
\end{equation}
is the retarded (or Huygens) potential. It is interesting to remark that while we know an integral expression for $\mathcal S*\eta$, it is not that easy to understand what $\mathcal S$ is. Just for the sake of illustration, let us give an idea how $\mathcal S$ looks. Consider the operator
\[
(\mathcal H(t) \xi)(\mathbf z) :=\int_{\Gamma \cap B(\mathbf z,t)} \frac{\xi(\mathbf y)}{4\pi|\mathbf z-\mathbf y|}
\mathrm d\Gamma(\mathbf y),
\]
where
\[
B(\mathbf z,t):=\{ \mathbf y\in \mathbb R^3\,:\, |\mathbf y-\mathbf z|< t\}.
\]
Using arguments in Fourier analysis it is possible to show that $\mathcal H$ is a continuous causal function of $t$ with values in the space of bounded linear operators from $H^{-1/2}(\Gamma)$ to $H^1(\mathbb R^3)$. Its distributional time derivative is $\mathcal S$. 

\subsection{A boundary integral equation for scattering}\label{sec:4.3}

\paragraph{The potential ansatz.} We want to find causal distributional solutions to
\begin{subequations}\label{eq:4.18}
\begin{alignat}{6}
c^{-2} \ddot u=\Delta u &\qquad && \mbox{in $L^2(\Omega_+)$}\\
\gamma u=\beta &&& \mbox{in $H^{1/2}(\Gamma).$}
\end{alignat}
\end{subequations}
In order to write an integral representation we need to modify the speed of the wave operators of Section \ref{sec:4.2}. We thus define
\[
\mathcal S_c:=\mathcal L^{-1}\{ \mathrm S(s/c)\}, \qquad 
\mathcal V_c:=\mathcal L^{-1} \{ \mathrm V(s/c)\}=\gamma \mathcal S_c,\qquad
\mathcal J_c:=\mathcal L^{-1} \{\mathrm J(s/c)\}=\ave{\partial_\nu\cdot}\mathcal S_c.
\]
We represent the solution of \eqref{eq:4.18} by means of a single layer potential
\begin{subequations}\label{eq:4.19}
\begin{equation}
u=\mathcal S_c*\eta
\end{equation}
for a causal density $\eta$ (an $H^{-1/2}(\Gamma)$-valued distribution) to be determined. We then impose the boundary condition to obtain an integral equation for the density
\begin{equation}
\mathcal V_c*\eta=\beta.
\end{equation}
The exterior normal derivative can then be computed as a postprocessing of $\eta$:
\begin{equation}
\lambda=-\smallfrac12\eta+\mathcal J_c*\eta=\partial_\nu^+ u=(\partial_\nu^+ \mathcal S_c)* \eta.
\end{equation}
\end{subequations}

\paragraph{The transfer functions.} From the point of view of data, equations \eqref{eq:4.19} involve three transfer functions
\begin{subequations}\label{eq:4.20}
\begin{alignat}{6}
\mathrm V(s/c)^{-1}, &\qquad &&\mbox{(convolution equation)}\\
\mathrm S(s/c)\mathrm V(s/c)^{-1}, & & & \mbox{(conv eqn followed by a forward conv)}\\
(-\smallfrac12\mathrm I+\mathrm J(s/c))\mathrm V(s/c)^{-1}. & & & \mbox{(same)}
\end{alignat}
\end{subequations}
We will come back to how these transfer functions behave as functions of $s\in \mathbb C_+$. For the moment, let us just understand what they do in the Laplace domain. If we focus on the problem
\[
\left(\smallfrac{s}{c}\right)^2  \mathrm U= \Delta U, \qquad \gamma \mathrm U=\mathrm B,
\]
and think of the integral representation $\mathrm U=\mathrm S(s/c)\mathrm E$, 
the associated operators are:
\begin{alignat*}{6}
\mathrm V(s/c)^{-1}, &\qquad &&\mbox{(solution of integral equation)}\\
\mathrm S(s/c)\mathrm V(s/c)^{-1}, & & & \mbox{(solution of exterior Dirichlet problem)}\\
(-\smallfrac12\mathrm I+\mathrm J(s/c))\mathrm V(s/c)^{-1}. & & & \mbox{(Dirichlet to Neumann operator)}
\end{alignat*}

\paragraph{Galerkin semidiscretization in space.} Let us now fix a finite dimensional subspace
\[
X_h \subset H^{-1/2}(\Gamma).
\]
For instance, if $\Gamma$ is a polygon, we can admit any triangulation of $\Gamma$ and the space of piecewise constant functions with respect to this triangulation. The semidiscrete version of \eqref{eq:4.19} is the search for a causal $X_h$-valued distribution $\eta^h$ satisfying
\begin{subequations}\label{eq:4.21}
\begin{equation}\label{eq:4.21a}
\langle \mu^h,\mathcal V_c*\eta^h-\beta\rangle_\Gamma =0 \qquad \forall \mu^h \in X_h,
\end{equation}
posprocessed to obtain a potential
\begin{equation}
u^h:=\mathcal S_c*\eta^h
\end{equation}
and the associated Neumann data
\begin{equation}
\lambda^h:=-\smallfrac12\eta^h+\mathcal J_c*\eta^h.
\end{equation}
\end{subequations}
It is inherent to the fact that we are discretizing using a boundary integral representation to see that $u^h$ is exactly a causal solution of the wave equation $c^{-2}\ddot u^h=\Delta u^h$. What has changed is the level of satisfaction of the boundary condition.

\paragraph{The semidiscrete system in the time domain.} It might help the reader understand the difficulties of dealing with equations like \eqref{eq:4.21} to see what happens in the simplest case. Assume that $\Gamma$ has been partitioned into $J$ non-overlapping elements $\{ \Gamma_1,\ldots,\Gamma_J\}$ and that we take
\[
X_h=\{ \mu^h :\Gamma \to \mathbb R\,:\, \mu^h|_{\Gamma_j}\in \mathcal P_0\quad\forall j\}=
\mathrm{span}\{ \chi_{\Gamma_1}, \ldots, \chi_{\Gamma_J}\},
\]
where $\mathcal P_0$ is the set of constant functions and we have used the symbol $\chi_{\Gamma_j}$ for the characteristic function of the element $\Gamma_j$. Let us assume that $\eta^h$ is actually a function $\mathbb R\to X_h$. We can then write
\[
\eta^h=\sum_{j=1}^J \eta_j(t) \chi_{\Gamma_j}(\mathbf y)
\]
for unknown causal scalar functions $\eta_j$. The equations \eqref{eq:4.21a} are equivalent to
\begin{equation}\label{eq:4.22}
\sum_{j=1}^J \int_{\Gamma_i} \int_{\Gamma_j} 
	\frac{\eta_j(t-c^{-1}|\mathbf x-\mathbf y|)}{4\pi |\mathbf x-\mathbf y|} 
	\mathrm d\Gamma(\mathbf y)\mathrm d\Gamma(\mathbf x)
=\int_{\Gamma_i} \beta(t,\mathbf x)\mathrm d\Gamma(\mathbf x) \qquad i=1,\ldots,J.
\end{equation}
This system is quite difficult to deal with in the time domain. The unknowns are subject to a continuum of delays which are then integrated against a weakly singular kernel. However, the Laplace transform of \eqref{eq:4.22} reveals much more about the problem. If $\mathrm E_j=\mathcal L\{ \eta_j\}$ and $\mathrm B=\mathcal L\{ \beta\}$, equations \eqref{eq:4.22} are transformed into
\begin{equation}\label{eq:4.23}
\sum_{j=1}^J \left(\int_{\Gamma_i}\int_{\Gamma_j} 
	\frac{e^{-\frac{s}{c} |\mathbf x-\mathbf y|}}{4\pi|\mathbf x-\mathbf y|}
	\mathrm d\Gamma(\mathbf y)\mathrm d\Gamma(\mathbf x)\right)
	\mathrm E_j(s) =
\int_{\Gamma_i} \mathrm B(s,\mathbf x)\mathrm d\Gamma(\mathbf x).
\end{equation}
It is interesting to note how the unknowns have left the integrals after the Laplace transform has been taken. The system \eqref{eq:4.23} corresponds to the $X_h$ Galerkin discretization of
\[
\int_{\Gamma} \frac{e^{-\frac{s}{c} |\cdot-\mathbf y|}}{4\pi|\cdot-\mathbf y|} 
	\mathrm E(s,\mathbf y)\mathrm d\Gamma (\mathbf y)
=\mathrm B(s,\cdot) \qquad \mbox{in $\Gamma$},
\]
which is an integral form of the equation
\[
\mathrm V(s/c)\mathrm E(s)=\mathrm B(s).
\]

\subsection{Full discretization}

\paragraph{First space, then time.} In Section \ref{sec:4.3} we have written a semidiscrete form of an integral formulation for the exterior scattering problem (by a sound-soft obstacle). These are the equations \eqref{eq:4.21} in the Laplace domain: we look for $\mathrm E^h:\mathbb C_+\to X_h$ satisfying
\begin{subequations}\label{eq:4.24}
\begin{equation}
\langle \mu^h , \mathrm V(s/c) \mathrm E^h(s) -\mathrm B(s)\rangle_\Gamma=0
	\qquad\forall \mu^h \in X_h, \quad \forall s\in \mathbb C_+.
\end{equation}
The solution of this equation is then postprocessed to compute a potential
\begin{equation}
\mathrm U^h(s)=\mathrm S(s/c)\mathrm E^h(s)
\end{equation}
and the exterior normal derivative
\begin{equation}
\Lambda^h(s)=-\smallfrac12\mathrm E^h(s)+\mathrm J(s/c)\mathrm E^h(s).
\end{equation}
\end{subequations}
It is very easy to describe a fully discrete method for \eqref{eq:4.24} (that is, \eqref{eq:4.21}) using this form. Note that there are three convolutions involved: one of them is a convolution equation and two other are needed for postprocessing of the solution. We just need to substitute the continuous differentiation parameter $s$ by $s_\kappa=\frac1\kappa \delta(e^{-s\kappa})$ within all the operators above: we now look for $\mathrm E^h_\kappa:\mathbb C_+\to X_h$ satisfying
\begin{subequations}\label{eq:4.25}
\begin{equation}\label{eq:4.25a}
\langle \mu^h , \mathrm V(s_\kappa/c) \mathrm E^h_\kappa (s) -\mathrm B(s)\rangle_\Gamma=0
	\qquad\forall \mu^h \in X_h, \quad \forall s\in \mathbb C_+,
\end{equation}
and postprocess twice
\begin{equation}
\mathrm U^h_\kappa(s)=\mathrm S(s_\kappa /c)\mathrm E^h_\kappa (s),
\end{equation}
\begin{equation}
\Lambda^h_\kappa (s)=-\smallfrac12\mathrm E^h_\kappa (s)+\mathrm J(s_\kappa /c)\mathrm E^h_\kappa(s).
\end{equation}
\end{subequations}

\paragraph{A look at the associated integral operators.} We are now going to pay attention to the integral equation in \eqref{eq:4.25}. Let us forget for a while about the space discretization (we will see later on how this is a legitimate move). We have an integral operator
\[
\eta \mapsto \mathrm V(s/c) \eta=
\int_{\Gamma}\frac{e^{-\frac{s}{c} |\cdot-\mathbf y|}}{4\pi |\cdot-\mathbf y|}\eta(\mathbf y)\mathrm d\Gamma(\mathbf y).
\]
When we substitute
\[
s\mapsto \smallfrac1\kappa (1-e^{-s\kappa}),
\]
that is, when we deal with the backward Euler CQ method, we need to expand
\[
\mathrm V(\smallfrac1{c\kappa} (1-\zeta))=\sum_{n=0}^\infty \omega_n^{\mathrm V}(c\kappa)\zeta^n,
\]
with operators given by the expressions
\begin{equation}\label{eq:4.26}
\omega_n^{\mathrm V}(c\kappa)\eta=\int_{\Gamma}
	\frac{e^{-\frac1{c\kappa}|\cdot-\mathbf y|}}{4\pi|\cdot-\mathbf y|} \frac1{n!} 
	\left(\frac{|\cdot-\mathbf y|}{c\kappa}\right)^n \eta(\mathbf y)\mathrm d\Gamma(\mathbf y).
\end{equation}
Note that the only operator of the sequence \eqref{eq:4.26} which has a weakly singular kernel is the one we need to invert in each time-step, that is, the one for $n=0$. This is the single layer operator for the very diffusive elliptic operator
\[
-\Delta +(c\kappa)^{-2}, 
\]
which seems to take the role of a transport operator. 
If instead of the backward Euler discretization we apply BDF2, the operators have a somewhat more complicated expression
\begin{equation}\label{eq:4.27}
\omega_n^{\mathrm V}(c\kappa)\eta=\int_{\Gamma}
	\frac{e^{-\frac3{2c\kappa}|\cdot-\mathbf y|}}{4\pi|\cdot-\mathbf y|} \frac1{n!} 
	\left(\frac{|\cdot-\mathbf y|}{2c\kappa}\right)^{n/2}
	H_n\left( \sqrt{\frac{2|\cdot-\mathbf y|}{c\kappa}}\right) 
	\eta(\mathbf y)\mathrm d\Gamma(\mathbf y),
\end{equation}
where $H_n$ is the $n$-th Hermite polynomial. Let us now go back to the fully discrete equation \eqref{eq:4.25a}. We will only pay attention to samples of this equation at the time steps. We then produce vectors $\boldsymbol\eta_n \in \mathbb R^J$ satisfying equations
\begin{equation}\label{eq:4.28}
\mathbf V_0 \boldsymbol\eta_n=\boldsymbol\beta_n-\sum_{m=1}^n \mathbf V_m \boldsymbol\eta_{n-m} \qquad n=0,1,\ldots,
\end{equation}
where $\mathbf V_n$ are the matrices with elements 
\[
\int_{\Gamma_i}\int_{\Gamma_j} 
\frac{e^{-\frac{\delta(0)}{c\kappa}|\mathbf x-\mathbf y|}}{4\pi|\mathbf x-\mathbf y|} 
P_n\left( \frac{|\mathbf x-\mathbf y|}{c\kappa}\right)
\mathrm d\Gamma(\mathbf y)\mathrm d\Gamma(\mathbf x)
\]
for adequate $P_n$ (see \eqref{eq:4.26} and \eqref{eq:4.27}), and where $\boldsymbol\beta_n$ is the vector with entries
\[
\int_{\Gamma_i} \beta(t_n,\mathbf x)\mathrm d\Gamma(\mathbf x).
\] 
It is a good moment to reiterate how CQ uses the Laplace transform of the operator but how data are sampled in time and the discrete solution is obtained by time-stepping.

\paragraph{First time, then space.} Let us go back again to the continuous convolutional system \eqref{eq:4.19} before Galerkin semidiscretization. Since the only term that has undergone spatial discretization is the equation $\mathcal V_c*\eta=\beta$, let us have a look at what happens if we first discretize in time using CQ and then in space using a Galerkin scheme. When we substitute $\mathrm V(s/c)\mathrm E(s)=\mathrm B(s)$ by $\mathrm V(s_\kappa/c)\mathrm E_\kappa(s)=\mathrm B(s)$ we are just moving from $\mathcal V_c*\eta=\beta$ to the sequence of problems
\begin{equation}\label{eq:4.29}
\omega_0^{\mathrm V}(c\kappa) \eta_n=
	\beta(t_n)
	-\sum_{m=1}^n \omega_m^{\mathrm V}(c\kappa) \eta_{n-m}.
\end{equation}
If we discretize all the equations \eqref{eq:4.29} with the same $X_h$-based Galerkin scheme, we obtain the sequence of linear system described in \eqref{eq:4.28}. Exactly the same. This means that for this family of problems, Convolution Quadrature and Galerkin discretization in space commute. 

\paragraph{A note on the all-steps-at-once method.} It is interesting to note that if we use Algorithm 3.IV to solve the convolution equations \eqref{eq:4.28} we are solving a collection of time-harmonic damped wave equations. In fact, if we want to compute $N$ steps of the process, we end up solving equations of the form
\[
\mathrm V\big( \smallfrac1{c\kappa} \delta (R \zeta_{N+1}^{-\ell}) \big)\widehat w_\ell=\widehat v_\ell, \qquad \ell=0,\ldots,N.
\]
These are the single layer operator equations associated to the operators
\[
\Delta-\omega_\ell^2, \quad \mbox{with}\quad \omega_\ell=\smallfrac1{c\kappa}\delta(R \zeta_{N+1}^{-\ell})).
\]

\subsection{Well posedness after discretization}\label{sec:4.5}

\paragraph{Why bother?} Following the previous pages, the reader can be led to believe that the application of CQ before or after space discretization is just a given. It is not however clear at all why we can even apply CQ after a Galerkin semidiscretization process and whether the semidiscrete operator inherits the properties of the continuous operator. In order to clarify concepts, let us compare the operators with what we get after semidiscretization in space. The first piece of good news is a coercivity estimate
\begin{equation}\label{eq:4.30}
\mathrm{Re} \left( e^{\imath\,\mathrm{Arg} s} \langle \overline\eta,\mathrm V(s)\eta\rangle_\Gamma\right)
\ge C \frac{\sigma\underline\sigma}{|s|^2} \|\eta\|_{-1/2,\Gamma}^2 
\qquad\forall \eta\in H^{1/2}(\Gamma)\qquad \forall s\in \mathbb C_+,
\end{equation}
where as in Section \ref{sec:4.2}, $\sigma:=\mathrm{Re}\,s$, $\underline\sigma:=\min\{1,\sigma\}$. The estimate \eqref{eq:4.30} is also often written as
\[
\mathrm{Re} \langle \overline\eta,s\mathrm V(s)\eta\rangle_\Gamma \ge C 
\frac{\sigma\underline\sigma}{|s|} \|\eta\|_{-1/2,\Gamma}^2 
\qquad\forall \eta\in H^{1/2}(\Gamma)\qquad \forall s\in \mathbb C_+,
\]
which shows how the operator that is actually coercive is $s\mathrm V(s)$, that is, application of $\mathrm V(s)$ followed (or preceded) by `differentiation.' Coercivity estimates are inherited by Galerkin discretizations. This means that if we compare $\mathrm A(s):=\mathrm V(s/c)^{-1}$ with the operator $\mathrm A_h(s):H^{1/2}(\Gamma) \to X_h$ that corresponds to solving
\[
\eta^h=\mathrm A_h(s)\beta \in X^h \qquad \langle \mu^h,\mathrm V(s) \eta^h-\beta\rangle_\Gamma =0\quad \forall \mu^h\in X_h,
\]
we have
\begin{equation}\label{eq:4.31}
\| \mathrm A(s)\|_{H^{1/2}(\Gamma)\to H^{-1/2}(\Gamma)}+
\| \mathrm A_h(s)\|_{H^{1/2}(\Gamma)\to H^{-1/2}(\Gamma)}
\le C \frac{|s|^2}{\sigma\underline\sigma} \qquad \forall s\in \mathbb C_+.
\end{equation}
(The constant $C$ is different in \eqref{eq:4.30} and \eqref{eq:4.31}, but we will adopt the bad habits of numerical analysts of calling all constants $C$.) The quality of the behavior with respect to $s$ and its real part in the right-hand-side of \eqref{eq:4.31} is highly relevant for the analysis of CQ as we will discuss in Section \ref{sec:7}. 

\paragraph{The postprocessed solutions.}
It would look like we were done if we had not postprocessed the solution in two different ways using the operators
\[
\mathrm S(s) \mathrm A_h(s) \qquad (-\smallfrac12\mathrm I+\mathrm J(s)) \mathrm A_h(s).
\]
In principle we can just go ahead and combine the estimate \eqref{eq:4.31} (due to coercivity of $s\mathrm V(s)$) with the estimates in \eqref{eq:4.16}. This yields the bounds
\begin{alignat*}{6}
\| \mathrm S(s)\mathrm A_h(s)\|_{H^{1/2}(\Gamma) \to H^1(\Omega_+)}  &
	 \le C \frac{|s|^3}{\sigma^2 \underline\sigma^2} &\qquad &\forall s\in \mathbb C_+,\\
\| (-\smallfrac12\mathrm I+\mathrm J(s) )\mathrm A_h(s) \|_{H^{1/2}(\Gamma) \to H^{-1/2}(\Gamma)} &
	\le C \frac{|s|^{7/2}}{\sigma^2 \underline\sigma^{5/2}} && \forall s\in \mathbb C_+.
\end{alignat*}
These estimates are quite pessimistic though.  With a different approach it can actually be proved that
\begin{subequations}
\begin{alignat}{6}
\label{eq:4.32a}
\| \mathrm S(s)\mathrm A_h(s)\|_{H^{1/2}(\Gamma) \to H^1(\Omega_+)}  &
	 \le C \frac{|s|^{3/2}}{\sigma \underline\sigma^{3/2}} &\qquad &\forall s\in \mathbb C_+,\\
\| (-\smallfrac12\mathrm I+\mathrm J(s) )\mathrm A_h(s) \|_{H^{1/2}(\Gamma) \to H^{-1/2}(\Gamma)} &
	\le C \frac{|s|^2}{\sigma \underline\sigma} && \forall s\in \mathbb C_+.
\end{alignat}
\end{subequations}
In view of \eqref{eq:4.31} and \eqref{eq:4.16}, the bound \eqref{eq:4.32a} is slightly shocking and should make the reader consider carefully the role of $s$ in the operators. In principle $s$ is differentiation. Therefore the bound for $\mathrm S(s)$ seems to say that application of the single layer potential comes with loss of time regularity in one unit, while the solution of the integral equation $\mathrm V(s)\eta=\beta$ (or of its Galerkin discretization) comes with a loss of two indices of smoothness in time. Their composition, however, loses only $3/2$ smoothness indices. (Take this idea of loss of regularity with multiples of $s$ with a grain of salt. Moving from bounds in the Laplace domain to time domain mapping properties is not optimal and yields some additional losses of smoothness in time.)  

\subsection*{Credits}

The presentation and analysis of the acoustic layer potentials in the time domain using their Laplace transforms can be traced back to the seminal work of Alain Bamberger and Toung Ha-Duong \cite{BaHa:1986a, BaHa:1986b}. They were the first to prove the coercivity estimate \eqref{eq:4.30}, which is the origin for the modern theory of acoustic layer potentials in the time domain. Here we have adopted the more systematic approach of \cite{LaSa:2009a, Sayas:2013}. The explicit distributional form of the three dimensional Huygens potentials, briefly mentioned at the end of Section \ref{sec:4.2}, first appears in \cite{LaSa:2009b}. The first use of CQ for a boundary integral equation, related to the heat equation, was given by Christian Lubich and Reinhold Schneider in \cite{LuSc:1992}. Two years later, Lubich himself made the first incursion of CQ applied to the boundary integral equations for the wave equation. The use of CQ for seveal kinds of elastic wave propagation phenomena using integral equations was quite extended in the engineering literature (see the monograph \cite{Schanz:2001} by Martin Schanz, one of the pioneers in the field) by the time the mathematical community went back to the topic.
The effect of Galerkin semidiscretization in space was studied in \cite{Lubich:1994}, but the postprocessing part (once the boundary integral equation is solved, input the result in a potential) was only studied in \cite{LaSa:2009a}. Other approaches for the study of layer potentials in the time domain are given in \cite{DoSa:2013, BaLaSa:2014}, using very different techniques that avoid the Laplace transform. The expansions for the operators related to CQ for the single layer operator equation were given in \cite{HaKrSa:2009}. For more precise information on the Sobolev space theory of boundary integral equations for steady state problems (which is required to understand the Laplace domain estimates of this chapter), the reader is recommended to explore Willian McLean's monograph \cite{McLean:2000}.


\section{Multistage convolution quadrature}

In this section we are going to introduce a new family of discretization methods for causal convolutions and convolution equations. The main difference with the multistep method will be in the fact that we will work simultaneously with several points in time (stages). The way we will develop the method, the Runge-Kutta (RK) steps will barely make an appearance. 
\begin{framed}\noindent
The general look of an RK-based discrete convolution for $y=f*g$ is
\[
\boldsymbol y_n = \sum_{m=0}^n W_m^{\mathrm F}(\kappa) \boldsymbol g_{n-m}
\]
where $\boldsymbol g_n:=(g(t_n+\kappa c_1),\ldots, g(t_n+\kappa c_p))^\top$ and $W_m^{\mathrm F}(\kappa)$ is a $p\times p$ matrix of operators, with values in the same space as $\mathrm F$. 
\end{framed}

\subsection{Some Runge-Kutta methods}\label{sec:5.1}

\paragraph{Vectorized notation.} We will accept the following (shorthand) vectorized form for evaluation of a function: 
\[
\boldsymbol c=\left[\begin{array}{c} c_1 \\ \vdots \\ c_p\end{array}\right]\quad\longmapsto\quad
g(t+\kappa\boldsymbol c)=\left[\begin{array}{c} g(t+\kappa c_1)\\ \vdots \\ g(t+\kappa c_p)\end{array}\right]\in \mathbb R^p.
\]
Similarly, if $f=f(t,y)$, then
\[
\boldsymbol c,\boldsymbol y\in \mathbb R^p \quad\longmapsto \quad f(t+\kappa\boldsymbol c,\boldsymbol y)=
\left[\begin{array}{c} f(t+\kappa c_1,y_1) \\ \vdots \\ f(t+\kappa c_p,y_p)\end{array}\right]
\in \mathbb R^p
\]

\paragraph{Implicit RK methods.} An implicit RK scheme is often presented through its Butcher tableau
\[
\begin{array}{c|c} \boldsymbol c & \boldsymbol A \\[1ex] \hline & \\[-1ex] & \boldsymbol b^\top\end{array} \qquad 
\boldsymbol b,\boldsymbol c\in \mathbb R^p, \qquad \boldsymbol A\in \mathbb R^{p\times p},
\]
with the conditions
\begin{equation}\label{eq:5.1}
\boldsymbol A\boldsymbol 1=\boldsymbol c, \qquad \boldsymbol b^\top \boldsymbol 1=1, \qquad \boldsymbol 1=(1,\ldots,1)^\top.
\end{equation}
The second condition in \eqref{eq:5.1} is necessary for convergence. The first condition in \eqref{eq:5.1} is related to the possibility of understanding time as an independent variable. All RK methods satisfy these conditions. Unfortunately the letter $s$ will be reserved for the variable in the Laplace transform, so we will use $p$ for the number of stages. The application of a step of an RK method for
\begin{equation}
\dot y=f(t,y)
\end{equation}
is based on the solution of a system of non-linear equations to compute the internal stages
\begin{equation}\label{eq:5.3}
\boldsymbol y_{n}=y_n\boldsymbol 1 + \kappa \boldsymbol A f(t_n+\kappa\boldsymbol c,\boldsymbol y_{n})
\end{equation}
followed by the computation of the next step
\begin{equation}\label{eq:5.4}
y_{n+1}=y_n+\kappa \boldsymbol b^\top f(t_n+\kappa\boldsymbol c,\boldsymbol y_{n}).
\end{equation}
The internal stages and the steps approximate
\[
\boldsymbol y_{n}\approx y(t_n+\kappa\boldsymbol c)\qquad y_n \approx y(t_n).
\]
Starting RK methods requires the knowledge of $y_0$. Unlike multistep methods, bringing vanishing information from the past will not change the initial step, which for us will always be $y_0=0$.

\paragraph{Some multistage differentiation formulas.} We next apply the above RK method to the search of causal solutions to the antiderivative problem
\[
\dot y=g \qquad\longleftrightarrow\qquad s\mathrm Y(s)=\mathrm G(s).
\]
The stages and steps for the method \eqref{eq:5.3}-\eqref{eq:5.4} translate into
\begin{equation}\label{eq:5.5}
\boldsymbol y_{n}=y_n\boldsymbol 1+\kappa\boldsymbol A \boldsymbol g_{n}, \qquad 
y_{n+1}=y_n+\kappa \boldsymbol b^\top \boldsymbol g_{n}, \qquad \boldsymbol g_{n}:=g(t_n+\kappa\boldsymbol c).
\end{equation}
It is clear from \eqref{eq:5.5} that causality of $g$ and the imposition of a causal discrete solution ($y_n=0$ and $\boldsymbol y_n=\boldsymbol 0$ for $n<0$) imposes $y_0=0$. It is also clear that for this particular equation (a quadrature), the internal stages do not play any role. We will still pay attention to them, since they are the quantities of our interest. We next write \eqref{eq:5.5} using the $\zeta$ transform. To do it, we introduce
\[
\mathrm Y(\zeta):=\sum_{n=0}^\infty y_n\zeta^n, \qquad 
\boldsymbol Y(\zeta):=\sum_{n=0}^\infty \boldsymbol y_n \zeta^n,\qquad
\boldsymbol G(\zeta):=\sum_{n=0}^\infty \boldsymbol g_n\zeta^n,
\]
and rewrite \eqref{eq:5.5} as
\[
\boldsymbol Y(\zeta)=\mathrm Y(\zeta)\boldsymbol 1+\kappa \boldsymbol A\boldsymbol G(\zeta), \qquad
\zeta^{-1}\mathrm Y(\zeta)=\mathrm Y(\zeta)+\kappa \boldsymbol b^\top \boldsymbol G(\zeta).
\]
Therefore
\[
\mathrm Y(\zeta)=\kappa\frac{\zeta}{1-\zeta}\boldsymbol b^\top  \boldsymbol G(\zeta),
\]
and from this
\begin{equation}\label{eq:5.6}
\boldsymbol Y(\zeta)=\kappa \left(\frac{\zeta}{1-\zeta}\boldsymbol 1\boldsymbol b^\top+\boldsymbol A\right) \boldsymbol G(\zeta).
\end{equation}
Equation \eqref{eq:5.6} is the discrete version of antidifferentiation $\mathrm Y(s)=s^{-1}\mathrm G(s)$, performed at the level of the stages. Therefore the $p\times p$ matrix
\begin{framed}
\[
\Delta(\zeta):= \left(\frac{\zeta}{1-\zeta}\boldsymbol 1\boldsymbol b^\top+\boldsymbol A\right)^{-1}
\]
\end{framed}
\noindent
will play the part of $\delta(\zeta)$ for the multistage case:
\[
s \mathrm Y(s) \mbox{ is discretized as } \smallfrac1\kappa \Delta(\zeta)\boldsymbol Y(\zeta).
\] 
In other words, $\frac1\kappa\Delta(\zeta)$ is the discrete transfer function for multistage differentiation. Note that the possibility of defining a discrete multistage differentiation operator by inverting the RK recurrence requires \eqref{eq:5.6} to be an invertible recurrence, which only happens when $\boldsymbol A$ is invertible.

\paragraph{A subclass of RK methods.} We consider the subclass of {\bf stiffly accurate} RK methods where the last row of $\boldsymbol A$ is $\boldsymbol b^\top$, that is,
\[
\boldsymbol e_p^\top \boldsymbol A=\boldsymbol b^\top, \qquad \boldsymbol e_p^\top=(0,\ldots,0,1),
\]
and therefore $c_m=1$ (multiply both sides by $\boldsymbol c$). Therefore, multiplying \eqref{eq:5.3} by $\boldsymbol e_p^\top$ and using \eqref{eq:5.4} it follows that
\begin{equation}\label{eq:5.7}
\boldsymbol e_p^\top \boldsymbol y_n=
y_n \boldsymbol e_p^\top\boldsymbol 1+\kappa \boldsymbol e_p^\top\boldsymbol A f(t_n+\kappa\boldsymbol c,\boldsymbol y_n)=
y_n+\kappa \boldsymbol b^\top f(t_n+\kappa\boldsymbol c,\boldsymbol y_n)=y_{n+1},
\end{equation}
which means that the last component of $\boldsymbol y_n$ is $y_{n+1}$ and we do not need to worry about the steps any more. 
For this subclass of methods we can recompute the discrete differentiation operator.  Using \eqref{eq:5.7}, we can write
\[
\boldsymbol y_n=y_n\boldsymbol 1+\kappa\boldsymbol A\boldsymbol g_n=
\boldsymbol 1\boldsymbol e_p^\top \boldsymbol y_{n-1}+\kappa\boldsymbol A\boldsymbol g_n,
\]
or, in the $\zeta$ domain,
\[
(\boldsymbol I-\zeta\boldsymbol 1\boldsymbol e_p^\top)\boldsymbol Y(\zeta)=\kappa\boldsymbol A\boldsymbol G(\zeta).
\]
In this case, differentiation is given by
\[
\smallfrac1\kappa\Delta(\zeta), \qquad \mbox{where}\qquad \Delta(\zeta):=\boldsymbol A^{-1}(\boldsymbol I-\zeta\boldsymbol 1\boldsymbol e_p^\top).
\]
Note that this is just an alternative formula for the matrix $\Delta(\zeta)$ defined in the more general case. 

\paragraph{Two examples.} The order three Radau IIa method is given by the table
\[
\begin{array}{c|cc}
1/3 & 5/12 & -1/12 \\
1 & 3/4 &1/4 \\
\hline
 & 3/4 & 1/4
\end{array}
\]
The order four Lobatto IIIc is given by the table
\[
\begin{array}{c|ccc}
 0 & 1/6 & -1/3 & 1/6 \\
 1/2 & 1/6 & 5/12 & -1/12 \\
 1 & 1/6 & 2/3& 1/6 \\
 \hline
  & 1/6 & 2/3 & 1/6
\end{array}
\]

\subsection{Elementary Dunford calculus}\label{sec:5.2}

\paragraph{Our next goal.} Now that we have a new approximation of the derivative (of $s$), we might want to define a new approximation of $\mathrm F(s)$. If $\mathrm F$ is scalar valued and entire (analytic in $\mathbb C$), it is not entirely difficult to define $\mathrm F(\boldsymbol B)$
for any $p\times p$ matrix $\boldsymbol B$ using the power series expansion for $\mathrm F$. When $\mathrm F$ is analytic only in part of $\mathbb C$ this is still doable, but not for every matrix: essentially we need all the eigenvalues of $\boldsymbol B$ to be in the domain of analyticity of $\mathrm F$. The process is however somewhat more complicated when we deal with operator-valued $\mathrm F$. 

\paragraph{Scalar functions of matrices.} (While this theory can be made much more general, we will keep it close to our assumptions on transfer functions.) Let $\mathrm F:\mathbb C_+\to \mathbb C$ be analytic and let $\lambda\in \mathbb C_+$. Then
\[
\mathrm F(\lambda)=\frac1{2\pi\imath} \oint_C (z-\lambda)^{-1} \mathrm F(z)\mathrm dz,
\]
where $C$ is a simple positively oriented closed path around $\lambda$. It does not take much imagination to figure out a definition for $\mathrm F(\boldsymbol\Lambda)$ where $\boldsymbol\Lambda=\mathrm{diag}(\lambda_1,\ldots,\lambda_p)$ is a diagonal matrix:
\begin{eqnarray*}
\mathrm F(\boldsymbol\Lambda)
	&:=&	\frac1{2\pi\imath} \oint_C (z\boldsymbol I-\boldsymbol\Lambda)^{-1} \mathrm F(z)\mathrm dz\\
	&=&	\mathrm{diag}\left(
				\frac1{2\pi\imath}\oint_C (z-\lambda_1)^{-1}\mathrm F(z)\mathrm dz,\ldots,
				\frac1{2\pi\imath}\oint_C (z-\lambda_p)^{-1}\mathrm F(z)\mathrm dz\right).
\end{eqnarray*}
The integral is done component by component, and the path $C$ has to enclose the values $\{\lambda_1,\ldots,\lambda_p\}$, that is, the spectrum of $\boldsymbol\Lambda$. If $\boldsymbol B=\boldsymbol P\boldsymbol\Lambda\boldsymbol P^{-1}$, where $\boldsymbol\Lambda$ is diagonal, then a simple computations yields
\[
\boldsymbol P\mathrm F(\boldsymbol\Lambda)\boldsymbol P^{-1}=
\frac1{2\pi\imath} \oint_C (z\boldsymbol I-\boldsymbol B)^{-1} \mathrm F(z)\mathrm dz,
\]
so it is just logical that we adopt the latter expression
\begin{equation}\label{eq:5.8}
\mathrm F(\boldsymbol B):=\frac1{2\pi\imath} \oint_C (z\boldsymbol I-\boldsymbol B)^{-1} \mathrm F(z)\mathrm dz,
\end{equation}
as a definition, even for non-diagonalizable matrices. As in previous cases, $C$ is a simple closed path in $\mathbb C_+$ surrounding the spectrum of $\boldsymbol B$.

\paragraph{Kronecker products.} Let $\boldsymbol B\in \mathbb C^{p\times q}$ and $F\in \mathcal B(X,Y)$. We then define
\[
\boldsymbol B \otimes F:=
\left[\begin{array}{ccc} 
	b_{11}F &\ldots & b_{1q}F \\
	\vdots & \ddots &\vdots\\
	b_{p1}F & \ldots & b_{pq}F
\end{array}\right]\in \mathcal B(X,Y)^{p\times q}\equiv \mathcal B(X^q,Y^p).
\]
\begin{framed}\noindent
This formula gives us the proper definition of an operator-valued function $\mathrm F:\mathbb C_+\to \mathcal B(X,Y)$ acting on a matrix:
\begin{equation}\label{eq:5.9}
\mathrm F(\boldsymbol B):=
	\frac1{2\pi\imath} \oint_C (z\boldsymbol I-\boldsymbol B)^{-1}\otimes \mathrm F(z)\mathrm dz.
\end{equation}
\end{framed}
\noindent
Once again $C$ is a closed path surrounding the spectrum of $\boldsymbol B$. The resulting $p\times p$ integrals take place in the Banach space $\mathcal B(X,Y)$, where $\mathrm F$ takes values. Equivalently, we can think of the integral as being computed in the Banach space $\mathcal B(X,Y)^{p\times p}\equiv \mathcal B(X^p,Y^p)$.

\paragraph{Key properties.} The fact that the Dunford calculus is given that name (calculus, not Dunford) is due to the fact that it interacts nicely with the algebra of operators. For instance, for every $\boldsymbol B$
\begin{equation}\label{eq:5.100}
\frac1{2\pi\imath} \oint_C (z\boldsymbol I-\boldsymbol B)^{-1}\otimes I_X\mathrm dz=
\left(\frac1{2\pi\imath} \oint_C (z\boldsymbol I-\boldsymbol B)^{-1}\mathrm dz\right)\otimes I_X=
\boldsymbol I\otimes I_X=I_{X^p},
\end{equation}
as long as $C$ surrounds the spectrum of $\boldsymbol B$. Also, if $\mathrm F:\mathbb C_+\to \mathcal B(X,Y)$ and $\mathrm G:\mathbb C_+\to\mathcal B(Z,X)$, then
\begin{equation}\label{eq:5.101}
(\mathrm F\mathrm G)(\boldsymbol B)=\mathrm F(\boldsymbol B)\mathrm G(\boldsymbol B)
\end{equation}
for every matrix with spectrum contained in $\mathbb C_+$. In particular, if $\mathrm F(s)$ is invertible for all $s\in \mathbb C_+$ and $\sigma(\boldsymbol B)\subset \mathbb C_+$, then by \eqref{eq:5.100} and \eqref{eq:5.101}, it is clear that $\mathrm F(\boldsymbol B)$ is invertible and 
\begin{equation}\label{eq:5.102}
\mathrm F(\boldsymbol B)^{-1}=\mathrm F^{-1}(\boldsymbol B).
\end{equation}

\paragraph{Some properties of Kronecker products.} We are next going to explore how to compute \eqref{eq:5.9} in the diagonalizable case by using only evaluations of $\mathrm F$ on the spectrum of $\boldsymbol B$. We need some preliminary work.
It is a simple exercise to prove that if $\boldsymbol B\in \mathbb C^{p\times q}$ and $\boldsymbol C\in \mathbb C^{q\times r}$, then
\begin{equation}\label{eq:5.10}
(\boldsymbol B\boldsymbol C)\otimes F=(\boldsymbol B\otimes I_Y)(\boldsymbol C\otimes F).
\end{equation}
The product in the right-hand-side of \eqref{eq:5.10} is the product of a matrix of operators in $\mathcal B(Y,Y)$ with a matrix of operators in $\mathcal B(X,Y)$. It can also be understood as a composition of an operator in $\mathcal B(Y^q,Y^p)$ with an operator in $\mathcal B(X^r,Y^q)$. Let now $\boldsymbol\Lambda$ be a $p\times p$ diagonal matrix. Then
\begin{equation}\label{eq:5.11}
(\boldsymbol\Lambda \boldsymbol C)\otimes F=(\boldsymbol\Lambda\otimes I_Y)(\boldsymbol C\otimes F)=
\left[\begin{array}{c} 
	\lambda_1 \mathrm{row}(\boldsymbol C,1)\otimes F \\
	\vdots\\
	\lambda_p \mathrm{row}(\boldsymbol C,p)\otimes F
\end{array}\right],
\end{equation}
where if $\boldsymbol C\in \mathbb C^{p,q}$, $\mathrm{row}(\boldsymbol C,i)$ is the $1\times q$ matrix containing the $i$-th row of $\boldsymbol C$ and therefore $\mathrm{row}(\boldsymbol C,i)\otimes F\in \mathcal B(X^q,Y)
\equiv\mathcal B(X,Y)^{1\times q}$. Using \eqref{eq:5.10} and \eqref{eq:5.11} we can compute
\begin{eqnarray*}
(\boldsymbol B\boldsymbol\Lambda\boldsymbol C)\otimes F
	&=&(\boldsymbol B\otimes I_Y)((\boldsymbol \Lambda\boldsymbol C)\otimes F)\\
	&=&\mathrm{col}(\boldsymbol B,1)\otimes \big( \lambda_1\mathrm{row}(\boldsymbol C,1)\otimes F\big)
		+\ldots+
		\mathrm{col}(\boldsymbol B,p)\otimes \big( \lambda_p\mathrm{row}(\boldsymbol C,p)\otimes F\big),
\end{eqnarray*}
or in short
\begin{equation}\label{eq:5.12}
(\boldsymbol B\boldsymbol\Lambda\boldsymbol C)\otimes F
=\sum_{i=1}^p 
\mathrm{col}(\boldsymbol B,i)\otimes \big( \lambda_i\mathrm{row}(\boldsymbol C,i)\otimes F\big).
\end{equation}
The outermost Kronecker product in \eqref{eq:5.12} corresponds to a column matrix $r\times 1$ with an operator in $\mathcal B(X^q,Y)$, outputting an operator in $\mathcal B(X^q,Y^r)$.

\paragraph{Operator-valued functions of a diagonalizable matrix.} Assume that $\boldsymbol B=\boldsymbol P\boldsymbol\Lambda\boldsymbol P^{-1}$ with $\boldsymbol\Lambda=\mathrm{diag}(\lambda_1,\ldots,\lambda_p)$, with $\lambda_i\in \mathbb C_+$ for all $i$, and that $\mathrm F:\mathbb C_+\to \mathcal B(X,Y)$ is analytic. Then, \eqref{eq:5.12} says that for $z\not\in\sigma(\boldsymbol B)=\{\lambda_i\}$
\begin{eqnarray*}
(z\boldsymbol I-\boldsymbol B)^{-1}\otimes \mathrm F(z)
&=& \sum_{i=1}^p 	\mathrm{col}(\boldsymbol P,i)\otimes\big((z-\lambda_i)^{-1}\mathrm{row}(\boldsymbol P^{-1},i)\otimes \mathrm F(z)\big)\\
&=& \sum_{i=1}^p 	\mathrm{col}(\boldsymbol P,i)\otimes\big(\mathrm{row}(\boldsymbol P^{-1},i)\otimes 
((z-\lambda_i)^{-1}\mathrm F(z))\big).
\end{eqnarray*}
Integrating on a path that surrounds the spectrum of $\boldsymbol B$ gives the following computable version of \eqref{eq:5.9}:
\[
\mathrm F(\boldsymbol B)=\sum_{i=1}^p\mathrm{col}(\boldsymbol P,i)
\otimes \big(\mathrm{row}(\boldsymbol P^{-1},i)\otimes \mathrm F(\lambda_i)\big).
\]

\subsection{RKCQ}\label{sec:5.3}

\paragraph{From discrete differentiation to discrete calculus.} In Section \ref{sec:5.1} we have defined a discrete differentiation symbol
\[
\smallfrac1\kappa\Delta(\zeta), \qquad \Delta(\zeta):=
	\left(\frac\zeta{1-\zeta}\boldsymbol 1\boldsymbol b^\top+\boldsymbol A\right)^{-1}.
\]
In Section \ref{sec:5.2} we have shown how to define an operator-valued function of a matrix variable using Dunford calculus. In particular, whenever this makes sense, we define
\[
\mathrm F(\smallfrac1\kappa\Delta(\zeta))=
	\frac1{2\pi\imath}\oint_C (z\boldsymbol I-\smallfrac1\kappa\Delta(\zeta))^{-1}\otimes\mathrm F(z) \mathrm d z,
\]
where $C$ is a path around the spectrum of $\smallfrac1\kappa\Delta(\zeta)$, which is assumed to be included in $\mathbb C_+$ for $|\zeta|<1$ (more about this at the end of this section). Then we use a Taylor expansion to obtain the coefficients:
\begin{equation}\label{eq:5.16}
\mathrm F(\smallfrac1\kappa\Delta(\zeta))=\sum_{n=0}^\infty W_n^{\mathrm F}(\kappa)\zeta^n.
\end{equation}

\paragraph{Multistage discrete convolutions.} To discretize $y=f*g$, we apply \eqref{eq:5.16} to obtain
\[
y(t_n+\kappa\boldsymbol c) \approx\boldsymbol y_n=\sum_{m=0}^n W_n^{\mathrm F}(\kappa)\boldsymbol g_{n-m},
\qquad \boldsymbol g_n:=g(t_n+\kappa\boldsymbol c).
\]
This same expression can be written with help of the $\zeta$ transform
\[
\boldsymbol Y(\zeta)=\mathrm F(\smallfrac1\kappa\Delta(\zeta))\boldsymbol G(\zeta).
\]
In the case of convolution equations $f*g=h$, we apply the same idea to obtain a lower triangular system of operator equations
\[
\sum_{m=0}^n W_m^{\mathrm F}(\kappa)\boldsymbol g_{n-m}=\boldsymbol h_n:=h(t_n+\kappa\boldsymbol c), \qquad n\ge 0,
\]
or, in more explicit form
\begin{equation}\label{eq:5.17}
W_0^{\mathrm F}(\kappa)\boldsymbol g_n=\boldsymbol h_n-\sum_{m=1}^nW_m^{\mathrm F}(\kappa)\boldsymbol g_{n-m}.
\end{equation}
Equation \eqref{eq:5.17} shows how we are inverting an operator equation associated to $\mathcal B(X^p;Y^p)$. Note that
\[
W_0^{\mathrm F}(\kappa)=\mathrm F(\smallfrac1\kappa\Delta(0))=\mathrm F(\smallfrac1\kappa\boldsymbol A^{-1})
\]
is invertible as shown in \eqref{eq:5.102}.

\paragraph{A note on the spectrum of $\boldsymbol A$.} A requirement for the correct definition of the RKCQ process is the possibility of producing the CQ coefficients $W_n^{\mathrm F}(\kappa)$. We recall that one of our prerequisites to define a multistage discrete derivative $\frac1\kappa\Delta(\zeta)$ was the existence of $\boldsymbol A^{-1}$. We will further assume the following property
\begin{equation}\label{eq:5.18}
\sigma(\boldsymbol A)\subset \mathbb C_+.
\end{equation}
Readers acquainted with the theory of $A$-stable RK methods will recognize that the invertibility of $\boldsymbol I-z\boldsymbol A$ for $\mathrm{Re}\,z\le 0$ is one of the hypotheses of $A$-stable methods. Actually this hypothesis and invertibility of $\boldsymbol A$ are equivalent to \eqref{eq:5.18}. Since $\Delta(\zeta)^{-1}=\frac{\zeta}{1-\zeta}\boldsymbol 1\boldsymbol b+\boldsymbol A$ is a small perturbation of $\boldsymbol A$ for small $\zeta$, hypothesis \eqref{eq:5.18} implies that for small $\zeta$ 
\[
\sigma(\smallfrac1\kappa\Delta(\zeta)) \subset \mathbb C_+
\]
and therefore the CQ coefficients are well defined.

\subsection{Stages, steps, and more}\label{sec:5.4}

\paragraph{Stages or steps.} As defined in Section \ref{sec:5.3} the RKCQ process works purely at the stage level. It samples data in the stages and then produces vectors of approximations in the internal stages. This is especially important when thinking of solving convolution equations, where we need to have as many data as unknowns. 

\paragraph{The associated convolution in continuous time.} Let us start by reviewing a simple fact of multistage CQ. Given the operator valued distribution $f$ and its Laplace transform $\mathrm F$, we had
\[
\mathrm F(\smallfrac1\kappa\delta(\zeta))=\sum_{n=0}^\infty \omega_n^{\mathrm F}(\kappa)\zeta^n,\qquad
\mathrm F(s_\kappa)=\sum_{n=0}^\infty \omega_n^{\mathrm F}e^{-st_n} \qquad 
s_\kappa:=\smallfrac1\kappa\delta(e^{-s\kappa}),
\]
and, in the time domain
\[
f_\kappa:=\sum_{n=0}^\infty \omega_n^{\mathrm F}\otimes \delta_{t_n}\qquad 
\mathcal L\{ f_\kappa\}(s)=\mathrm F(s_\kappa).
\]
This means that even if we only computed convolution at discrete times, there is a continuous convolutional operator in the background. 
The case of RKCQ is more complicated. We first try to replicate the previous formulas by defining
\[
\boldsymbol S_\kappa=\boldsymbol S_\kappa(s):=\smallfrac1\kappa\Delta(e^{-s\kappa}),\qquad
\mathrm F(\boldsymbol S_\kappa)=\sum_{n=0}^\infty W_n^{\mathrm F}(\kappa) e^{-s t_n},
\]
the latter $\mathcal B(X^p;Y^p)$-valued function of $s$ being the Laplace transform of the causal distribution
\[
\boldsymbol F_\kappa:=W_n^{\mathrm F}(\kappa) \otimes \delta_{t_n}
\]
This operator-valued distribution cannot be put in convolution with an $X$-valued distribution $g$. Instead, the distribution $g$ is first modified to the $X^p$-valued distribution
\[
\left[\begin{array}{c} g(\cdot+c_1\kappa) \\ \vdots \\ g(\cdot+c_p\kappa)\end{array}\right]=:g(\cdot+\kappa\boldsymbol c).
\]
Note that this can be understood as a convolution process, {\em but it is not causal}, because it pushes back the origin to be at the level of the different stages. Then the RKCQ process can be understood as the evaluation at the time steps $t_n$ of the continuous convolution
\[
\boldsymbol F_\kappa* g(\cdot+\kappa\boldsymbol c)=\sum_{n=0}^\infty W_n^{\mathrm F}(\kappa) g(\cdot-t_n+\boldsymbol c\kappa).
\]

\paragraph{Computation of steps in the simplest case.} While the emphasis in the multistage CQ process is on the stages, we might want to compute only values at the time steps $t_n$. In the case when $\boldsymbol e_p^\top\boldsymbol A=\boldsymbol b^\top$, we can compute
\[
(\boldsymbol e_p^\top\otimes I) \boldsymbol y_n=:y_{n+1}\approx y(t_{n+1}).
\]
This means that we have never computed an approximation $y_0$ (it is zero), unlike in the multistep case, where $y_0$ was computed from $g(0)$. 

\paragraph{The general case.} In order to give a definition of a multistage CQ method where the input are the vectors $\boldsymbol g_n$ and the output are quantities $y_{n+1}\approx (f*g)(t_{n+1})$, we need to go back to some computations in Section \ref{sec:5.1}. In particular we have shown that for the differential equation $y'=g$ (that is, for the operator $\mathrm F(s)=s^{-1}$), we could compute the steps as a postprocessing of the stages
\[
\mathrm Y(\zeta)=\frac{\zeta}{1-\zeta}\boldsymbol b^\top \Delta(\zeta)\boldsymbol Y(\zeta), \qquad 
\boldsymbol Y(\zeta)=\kappa \Delta(\zeta)^{-1}\boldsymbol G(\zeta).
\]
We can use this formula to extend the computation of steps for a general convolution $y=f*g$, by writing
\begin{equation}\label{eq:5.19}
\mathrm Y(\zeta)=\frac{\zeta}{1-\zeta} \big(\boldsymbol b^\top\Delta(\zeta)\otimes I)
\mathrm F(\smallfrac1\kappa\Delta(\zeta))\boldsymbol G(\zeta).
\end{equation}
This would suggest that we need to figure out a way of computing the discrete process described by \eqref{eq:5.19}. There is however a simpler formula to compute the steps from the previous step and the most recently computed stages. We start by computing
\begin{alignat*}{6}
(1-\zeta)\boldsymbol b^\top \boldsymbol A^{-1}\Delta(\zeta)^{-1}
	&=(1-\zeta)\boldsymbol b^\top\big( \boldsymbol I+\frac\zeta{1-\zeta}\boldsymbol A^{-1}\boldsymbol 1\boldsymbol b^T\big)\\
	&=(1-\zeta)\boldsymbol b^\top +\zeta(\boldsymbol b^\top \boldsymbol A^{-1}\boldsymbol 1)\boldsymbol b^\top
		=(1-\zeta\mu)\boldsymbol b^\top,
\end{alignat*}
where
\begin{equation}\label{eq:5.200}
\mu:=1-\boldsymbol b^\top\boldsymbol A^{-1}\boldsymbol 1.
\end{equation}
(More about this quantity later.) Therefore
\[
\frac\zeta{1-\zeta}\boldsymbol b^\top \Delta(\zeta)=\frac{\zeta}{1-\mu\zeta}\boldsymbol b^\top\boldsymbol A^{-1},
\]
which means that the discrete convolution
\[
\mathrm Y(\zeta)=\frac\zeta{1-\zeta}(\boldsymbol b^\top \Delta(\zeta)\otimes I)\boldsymbol Y(\zeta)
\]
is equivalent to
\begin{equation}\label{eq:5.201}
(\zeta^{-1}-\mu)\mathrm Y(\zeta)=(\boldsymbol d^\top\otimes I)\boldsymbol Y(\zeta),\qquad 
\boldsymbol d^\top:=\boldsymbol b^\top\boldsymbol A^{-1},
\end{equation}
which allows us to write
\begin{equation}\label{eq:5.202}
y_{n+1}=\mu y_n+(\boldsymbol d^\top \otimes I)\boldsymbol y_n
		=\mu y_n+d_1 y_{n,1}+\ldots+d_p y_{n,p}.
\end{equation}
The simple case is easily recovered from this formula by noticing that $\boldsymbol b^\top \boldsymbol A^{-1}=\boldsymbol e_p^\top$ and $\mu=0$.

\subsection{Implementation of RKCQ}\label{sec:5.5}

\paragraph{General idea.} Much of what we are going to sketch in this section follows closely what was explained in Section \ref{sec:3} for the multistage CQ scheme. We will not repeat many of the arguments there, and will just show some important steps. The key formula to keep in mind is the practical computation
\begin{equation}\label{eq:5.20}
\mathrm F(\boldsymbol B)=\sum_{i=1}^p\mathrm{col}(\boldsymbol P,i)
\otimes \big(\mathrm{row}(\boldsymbol P^{-1},i)\otimes \mathrm F(\lambda_i)\big), \quad 
\mbox{when $\boldsymbol B=\boldsymbol P \mathrm{diag}(\lambda_1,\ldots,\lambda_p)\boldsymbol P^{-1}$},
\end{equation}
that was derived in Section \ref{sec:5.2}.

\paragraph{Computation of RKCQ coefficients.} Using an integration contour $C_R:=\{ \zeta\in \mathbb C\,:\,|\zeta|=R\}$, with $R=\epsilon^{1/(2N+2)}$ and a trapezoidal rule of $N+1$ points, we can compute
\begin{alignat*}{6}
W_n^{\mathrm F}(\kappa)
	&=\frac1{n!}\frac{\mathrm d^n}{\mathrm d\zeta^n} \mathrm F\left(\frac1\kappa\Delta(\zeta)\right)\Big|_{\zeta=0}\\
	&=\frac1{2\pi\imath}\oint_{C_R} \zeta^{-n-1}\mathrm F(\smallfrac1\kappa\Delta(\zeta))\mathrm d\zeta \\
      & \approx \frac{R^{-n}}{N+1}\sum_{\ell=0}^N 
				\zeta_{N+1}^{n\ell}\mathrm F(\smallfrac1\kappa\Delta(\zeta_{N+1}^{-\ell})), \qquad n=0,\ldots,N.
\end{alignat*}
As usual $\zeta_{N+1}=e^{\frac{2\pi\imath}{N+1}}$.
The corresponding algorithm is to be compared with Algorithm 3.II.

\begin{framed}
\paragraph{Algorithm 5.I (computation of RKCQ coefficients).} Note that in the particular case of stiffly accurate methods, we can write
\[
\Delta(\zeta)=\boldsymbol A^{-1} -\zeta \boldsymbol C, \qquad \boldsymbol C:=\boldsymbol A^{-1} \boldsymbol 1 \boldsymbol e_p^\top.
\]
\begin{itemize}
\item[(a)] For $\ell=0,\ldots,N$, find the spectral decomposition 
\[
\boldsymbol P_\ell \boldsymbol\Lambda_\ell \boldsymbol P_{\ell}^{-1} = \smallfrac1\kappa\Delta(R \zeta_{N+1}^{-\ell})
\]
and use \eqref{eq:5.20} (looping over stages) to compute
\[
\widehat{\mathrm F}_{\ell} :=\mathrm F(\smallfrac1\kappa\Delta(R \zeta_{N+1}^{-\ell})).
\]
Note that $\widehat{\mathrm F}_{\ell}\in \mathcal B(X^p,Y^p)$.
\item[(b)] Apply the IDFT and scale
\[
W_n^{\mathrm F}(\kappa):=R^{-n} \left(\frac1{N+1}\sum_{\ell=0}^N \widehat{\mathrm F}_\ell\zeta_{N+1}^{n\ell}\right).
\]
\end{itemize}
This method works under the assumption that $\Delta(\zeta)$ is diagonalizable on the path $C_R$. 
\end{framed}

\paragraph{All-steps-at-once forward convolution.} Let us first sample the data
\[
\boldsymbol g_n:=g(t_n+\kappa\boldsymbol c), \qquad n=0,\ldots,N.
\]
We want to compute
\[
\boldsymbol u_n=\sum_{m=0}^n W_{n-m}^{\mathrm F}(\kappa) \boldsymbol g_m=
\sum_{m=0}^N W_{n-m}^{\mathrm F}(\kappa) \boldsymbol g_m, \qquad n=0,\ldots,N,
\]
(note that $W_n^{\mathrm F}(\kappa)=0$ for negative $n$), using approximations
\[
W_n^{\mathrm F}(\kappa) \approx \frac{R^{-n}}{N+1}\sum_{\ell=0}^N \widehat{\mathrm F}_\ell \zeta_{N+1}^{n\ell},
\qquad
\widehat{\mathrm F}_{\ell} :=\mathrm F(\smallfrac1\kappa\Delta(R \zeta_{N+1}^{-\ell})).
\]
Proceeding as in Section \ref{sec:3.3}, we approximate
\[
\boldsymbol u_n \approx R^{-n} \left(\frac1{N+1}\sum_{\ell=0}^N \widehat{\mathrm F}_\ell 
		\left( \sum_{m=0}^N R^m \boldsymbol g_m\zeta_{N+1}^{-m\ell}\right) \zeta_{N+1}^{\ell}\right).
\]

\paragraph{The finite dimensional case.}
Let us briefly detail what is to be done when $g:\mathbb R\to \mathbb R^{d_2}$ and $f:\mathbb R\to \mathbb R^{d_1\times d_2}$. In this case, it is advantageous to deal with samples at stage points as vectors $\boldsymbol g_n\in \mathbb R^{p\,d_2}$ organized in $p$ blocks of $d_2$ components. The key step is the multiplication
\[
\widehat{\mathrm F}_\ell\widehat{\boldsymbol h}_\ell, \qquad 
\widehat{\mathrm F}_\ell=\mathrm F(\smallfrac1\kappa\Delta(R\zeta_{N+1}^{-\ell})),
\]
for a given vector $\boldsymbol h\in \mathbb C^{p\,d_2}$. It is not difficult to see that when $\smallfrac1\kappa\Delta(R\zeta_{N+1}^{-\ell})=\boldsymbol P_\ell \boldsymbol\Lambda_\ell\boldsymbol P_\ell^{-1}$, then
\begin{equation}\label{eq:5.21}
\widehat{\mathrm F}_\ell\widehat{\boldsymbol h}_\ell=
(\boldsymbol P_{\ell}\otimes \boldsymbol I_{d_1})\mathrm{diag}(\mathrm F(\lambda_1),\ldots,\mathrm F(\lambda_p))
 (\boldsymbol P_{\ell}^{-1}\otimes \boldsymbol I_{d_2})\widehat{\boldsymbol h}_\ell.
\end{equation}
If $g$ takes values in $X$ and $F$ in $\mathcal B(X;Y)$, then we have to deal with samples $\boldsymbol g_n\in X^p$, with some modified vectors $\widehat{\boldsymbol h}_\ell\in X^p$ (see Algorithm 5.II) and that \eqref{eq:5.21} still applies if we substitute $\boldsymbol I_{d_1}$ and $\boldsymbol I_{d_2}$ by $I_Y$ and $I_X$ respectively. 

\begin{framed}
\paragraph{Algorithm 5.II (all-steps-at-once convolution).} We begin by sampling data at the discrete times
\[
\boldsymbol g_n:=g(t_n+\kappa \boldsymbol c)\in X^p, \qquad n=0,\ldots,N.
\]
\begin{itemize}
\item[(a)] Scale data:
\[
\boldsymbol h_m:=R^m \boldsymbol g_m, \qquad m=0,\ldots,N.
\]
\item[(b)] Compute the DFT:
\[
\widehat{\boldsymbol h}_\ell:=\sum_{m=0}^N \boldsymbol h_m \zeta_{N+1}^{-m\ell}, \qquad \ell=0,\ldots,N.
\]
These are formal DFTs, componentwise in $X^p$. When $X=\mathbb R^d$, these can be broken to $p\,d$ separate scalar DFTs.
\item[(c)] For every $\ell=0,\ldots,N$, find the spectral decomposition $\smallfrac1\kappa\Delta(R\zeta_{N+1}^{-\ell})=\boldsymbol P_\ell \boldsymbol\Lambda_\ell\boldsymbol P_\ell^{-1}$ and compute
\[
\widehat{\boldsymbol v}_\ell:=\widehat{\mathrm F}_\ell\widehat{\boldsymbol h}_\ell=
(\boldsymbol P_{\ell}\otimes I_Y)\mathrm{diag}(\mathrm F(\lambda_1),\ldots,\mathrm F(\lambda_p))
 (\boldsymbol P_{\ell}^{-1}\otimes I_X)\widehat{\boldsymbol h}_\ell.
\]
Note that the product by the block-diagonal matrix in the center can be done componentwise in $X$.
\item[(d)] Compute the IDFT:
\[
\boldsymbol v_n:=\frac1{N+1} \sum_{\ell=0}^N \widehat{\boldsymbol v}_\ell\zeta^{\ell n}_{N+1}, \qquad n=0,\ldots,N.
\]
(See the comments on step (b).)
\item[(e)] Scale back
\[
\boldsymbol u_n:=R^{-1} \boldsymbol v_n \in Y^p, \qquad n=0,\ldots,N.
\]
\end{itemize}
\end{framed}

If this is the last step of a sequence of convolutions, that is, if we are not going to apply any other convolution operator to this result, we can keep the last component of $\boldsymbol u_n$ as approximation in the point $t_{n+1}$. Note that because the RKCQ method counts intervals (groups of stages) and not steps, we are not computing an approximation at $t_0$, and the final time-step takes us to $t_{N+1}$ and not to $t_N$.

\paragraph{Convolution equations.} We will not repeat the argument for convolution equations. Algorithm 5.II can be easily modified to handle this new situation. The only step to be changed is (c), where we need a multiplication
\[
\widehat{\mathrm F}_\ell^{-1} \widehat{\boldsymbol v}_m=
\widehat{\mathrm F}_\ell\widehat{\boldsymbol h}_\ell=
(\boldsymbol P_{\ell}\otimes I_Y)\mathrm{diag}(\mathrm F(\lambda_1)^{-1},\ldots,\mathrm F(\lambda_p)^{-1})
 (\boldsymbol P_{\ell}^{-1}\otimes I_X)\widehat{\boldsymbol v}_\ell,
\]
that is, we need to solve $p$ equations associated to the operators $\mathrm F(\lambda_j)$. For a comparison, see how Algorithm 3.III is modified to Algorithm 3.IV.

\begin{framed}
\paragraph{Algorithm 5.III (computation of a piece of a convolution)} The algorithm to compute 
\[
\boldsymbol g_n:=\sum_{m=0}^Q W_n^{\mathrm F}(\kappa)\boldsymbol u_m,
\qquad n=Q+1,\ldots,M,
\]
renumbered in the form
\[
\widetilde{\boldsymbol g}_k:=\sum_{m=0}^Q W_{k+Q+1-m}^{\mathrm F}(\kappa) \boldsymbol u_m, \qquad k=0,\ldots,M-Q-1, 
\]
starting from vectors $\boldsymbol u_m \in X^p$ and outputtingvalues in $Y^p$, and using approximations
\begin{equation}\label{eq:5.22}
W_n^{\mathrm F}(\kappa) \approx \frac{R^{-n}}{N+1}\sum_{\ell=0}^N \widehat{\mathrm F}_\ell \zeta_{N+1}^{n\ell},
\qquad
\widehat{\mathrm F}_{\ell} :=\mathrm F(\smallfrac1\kappa\Delta(R \zeta_{N+1}^{-\ell}))
\end{equation}
(for positive and negative $n$) is derived in an entirely similar way to what we did in Section \ref{sec:3.4}. The parameter $N\ge M$ is a design parameter that influences the size of the computation, but also the precision to which the approximations \eqref{eq:5.22} are carried out. 
\begin{itemize}
\item[(a)] Scale and augment data
\[
\boldsymbol w_m:=
\left\{\begin{array}{ll} 
	R^m \boldsymbol u_m, & 0\le m\le Q,\\
	\boldsymbol 0, & Q+1\le m\le N.
\end{array}\right.
\]
\item[(b)] Compute the DFT $\widehat{\boldsymbol w}_\ell$ ($\ell=0,\ldots,N$) of the vectors in (a). See Algorithm 5.II(b) for a comment on this step.
\item[(c)] For every $\ell=0,\ldots,N$, find the spectral decomposition $\smallfrac1\kappa\Delta(R\zeta_{N+1}^{-\ell})=\boldsymbol P_\ell \boldsymbol\Lambda_\ell\boldsymbol P_\ell^{-1}$ and compute
\[
\widehat{\boldsymbol h}_\ell:=\zeta_{N+1}^{\ell(Q+1)} \widehat{\mathrm F}_\ell\widehat{\boldsymbol w}_\ell=
\zeta_{N+1}^{\ell(Q+1)}(\boldsymbol P_{\ell}\otimes I_Y)\mathrm{diag}(\mathrm F(\lambda_1),\ldots,\mathrm F(\lambda_p))
 (\boldsymbol P_{\ell}^{-1}\otimes I_X)\widehat{\boldsymbol w}_\ell.
\]
\item[(d)] Compute the IDFT of the sequence in (c), $\boldsymbol h_\ell$ ($\ell=0,\ldots,N$).
\item[(e)] Scale and chop the resulting sequence
\[
\widetilde{\boldsymbol g}_k:=R^{-k-Q-1} \boldsymbol h_\ell, \qquad k=0,\ldots,M-Q-1.
\]
\end{itemize}
\end{framed}

\subsection*{Credits}

For a deeper introduction to the Dunford calculus, the reader is referred to \cite{DaLi:1990b}.
Runge-Kutta convolution quadrature first appeared in a paper by Christian Lubich and Alexander Ostermann \cite{LuOs:1993}. Some further theoretical developments can be found in \cite{CaCuPa:2007}. The interest in RKCQ in the area of time-domain boundary integral equations is more recent \cite{Banjai:2010, BaMeSc:2012}. The algorithms shown in Section \ref{sec:5.5} can be found in  \cite{BaSc:2012}. The analysis of RKCQ applied to problems whose transfer function has the structure \eqref{eq:1.100}-\eqref{eq:1.101} was developed in \cite{BaLu:2011} and \cite{BaLuMe:2011}.


\section{A toy application}\label{sec:6}

In this section we will show a simple fully discrete (and very easy to code) example for the scattering of an acoustic wave by a smooth obstacle in the plane. All the operators will be given directly through their Laplace domain representations. The numerical method that we will present here can be understood as a fully discrete version of a Galerkin method. While the theory for the Galerkin method follows from existing arguments in the literature, the full discretization of the equations is not entirely justified. 

\subsection{Scattering by a smooth closed obstacle}

\paragraph{The geometry.} Consider a simple smooth closed curve in the plane parametrized by a 1-periodic function $\mathbf x:\mathbb R\to \Gamma\subset \mathbb R^2$ satisfying:
\[
\mathbf x(r)=\mathbf x(r+1) \quad \forall r, \qquad \mathbf x(r)\neq \mathbf x(\rho)\quad \mbox{if $r-\rho\not\in \mathbb Z$},
\qquad |\mathbf x'(r)|\neq 0 \quad \forall r.
\]
We assume that the parametrization gives a positive orientation to the curve, so that
\[
\mathbf n(r):=(x_2'(r),-x_1'(r))
\]
is a normal outward pointing vector at $\mathbf x(r)$.

\paragraph{One potential and two operators.} Given a 1-periodic density $\eta:\mathbb R\to\mathbb C$, directly defined in parametric space, the associated single layer potential at speed one (given in the Laplace domain) is
\begin{equation}\label{eq:6.1}
(\mathrm S(s)\eta)(\mathbf z):=
\frac\imath4\int_0^1 H^{(1)}_0(\imath s |\mathbf z-\mathbf x(\rho)|)\eta(\rho)\mathrm d\rho.
\end{equation}
Two operators can be used to represent the boundary values of the single layer potential, the single layer operator
\begin{equation}\label{eq:6.2}
(\mathrm V(s)\eta)(r):=
\frac\imath4\int_0^1 H^{(1)}_0(\imath s |\mathbf x(r)-\mathbf x(\rho)|)\eta(\rho)\mathrm d\rho,
\end{equation}
and the transposed double-layer operator
\begin{equation}\label{eq:6.3}
(\mathrm J(s)\eta)(r):=
	\frac{s}4\int_0^1 H^{(1)}_1(\imath s|\mathbf x(r)-\mathbf x(\rho)|) 
	\frac{(\mathbf x(r)-\mathbf x(\rho))\cdot\mathbf n(r)}{|\mathbf x(r)-\mathbf x(\rho)|}\eta(\rho)\mathrm d\rho.
\end{equation}
The functions $H^{(1)}_0$ and $H^{(1)}_1$ in \eqref{eq:6.1}-\eqref{eq:6.3} are the Hankel functions of the first kind and respective orders zero and one. The exterior boundary values for $U=\mathrm S(s)\eta$ are given by the expressions:
\begin{equation}\label{eq:6.4}
(U^+\circ\mathbf x)(r)=(\mathrm V(s)\eta)(r), \qquad
(\nabla U^+\circ\mathbf x)(r)\cdot\mathbf n(r)=-\smallfrac12\eta(r)+(\mathrm J(s)\eta)(r).
\end{equation}
Note how instead of using a unit normal vector, we are employing the non-normalized vector field $\mathbf n$ in \eqref{eq:6.3} and \eqref{eq:6.4}. This simplifies some expressions like the exterior normal derivative in \eqref{eq:6.4}.

\paragraph{The transient problem.} We are then going to bring a plane incident wave to the game (see Section \ref{sec:4.1})
\[
u^{\mathrm{inc}}(\mathbf z,t)=\psi(c(t-t_{\mathrm{lag}})-\mathbf z\cdot\mathbf d), \qquad |\mathbf d|=1.
\] 
This incident wave is read on points of the boundary to create
\begin{equation}\label{eq:6.4a}
\beta(t)(r)=\beta(r,t):=\psi(c(t-t_{\mathrm{lag}})-\mathbf x(r)\cdot\mathbf d),
\end{equation}
which is periodic in $r$ and assumed causal in $t$. The scattering problem then looks for $\eta:\mathbb R^2\to \mathbb R$, 1-periodic in its first variable and causal in the second such that
\[
(\mathcal V_c*\eta)(t)+\beta(t)=0 \qquad \forall t, \qquad \mbox{where}\quad \mathcal L\{ \mathcal V_c\}=\mathrm V(\cdot/c).
\]
(Note how it has been convenient to think of functions as being only functions of the time variable with output in a certain non-specified space of 1-periodic functions.) The density is then input in a potential expression
\[
U(t)=(\mathcal S_c*\eta)(t), \qquad \mbox{where}\quad \mathcal L\{\mathcal S_c\}=\mathrm S(\cdot/c),
\]
and is used to generate the exterior normal derivative
\[
\lambda(t)=-\smallfrac12\eta(t)+(\mathcal J_c*\eta)(t), \qquad \mbox{where}\quad \mathcal L\{\mathcal J_c\}=\mathrm J(\cdot/c).
\]

\subsection{Fully discrete equations}

\paragraph{The source geometry.} Let us choose a positive integer $N$ and consider $h:=1/N$ as the discrete mesh-size. The geometry is sampled in the following simple way:
\begin{equation}\label{eq:6.5}
\mathbf m_j:=\mathbf x(j \,h), \qquad \mathbf n_j:= h\mathbf n(j\,h), \qquad j\in \mathbb Z_N,
\end{equation}
where $\mathbb Z_N$ is the set of integers counted modulo $N$. (It is clear that \eqref{eq:6.5} defines only $N$ different points, since both $\mathbf x$ and $\mathbf n$ are 1-periodic functions.) A discrete density is now a vector $\boldsymbol\eta\in \mathbb C^N$, or, more properly (and pedantically) speaking, a function $\boldsymbol\eta:\mathbb Z_N \to \mathbb C$. The associated single layer potential is then given by a sum of sources
\begin{equation}
(\mathbf S(s)\boldsymbol\eta)(\mathbf z):=
	\frac\imath4 \sum_{j=1}^N H^{(1)}_0(\imath s|\mathbf z-\mathbf m_j|)\eta_j.
\end{equation}

\paragraph{The observation geometries.} Because of the logarithmic singularity of the Hankel function $H^{(1)}_0$ at the origin, we are not allowed to use a simple evaluation of \eqref{eq:6.2} at the same points were we have concentrated the density. To overcome this difficulty, and for reasons that will be discussed at the end of this section, we are going to choose two observation grids:
\begin{equation}\label{eq:6.7}
\mathbf m_i^\pm:=\mathbf x((i\pm\smallfrac16)h), \qquad \mathbf n_i^\pm:=h\mathbf n((i\pm\smallfrac16)h), \qquad
i\in \mathbb Z_N.
\end{equation}
For averaging any pair of discrete functions defined on the observation grids we will use the following notation:
\[
\sum_\pm a_i^\pm:=\smallfrac12 (a_i^++a_i^-).
\]
The discrete version of \eqref{eq:6.2} is given by
\begin{equation}\label{eq:6.8}
(\mathbf V(s)\boldsymbol\eta)_i:=
\sum_j \sum_\pm \frac\imath4  H^{(1)}_0(\imath s |\mathbf m_i^\pm-\mathbf m_j|) \eta_j,
\end{equation}
while for \eqref{eq:6.3} we use a first identical format
\begin{equation}\label{eq:6.7}
(\mathbf J^\circ(s)\boldsymbol\eta)_i:=
\sum_j \sum_\pm \frac{s}4 H^{(1)}_1(\imath s|\mathbf m_i^\pm-\mathbf m_j|)
\frac{(\mathbf m_i^\pm-\mathbf m_j)\cdot\mathbf n_i^\pm}{|\mathbf m_i^\pm-\mathbf m_j|} \eta_j,
\end{equation}
that we next correct in the following form
\begin{equation}
\mathbf J(s)=\mathbf Q\,\mathbf J^\circ(s), \qquad \mbox{where} \quad 
\mathbf Q_{ij}:=
	\left\{\begin{array}{ll} 
		\smallfrac{11}{12} & i=j, \\[1.5ex] \smallfrac1{24} & i=j\pm 1 \,(\mathrm{mod} N), \\[1.5ex] 0 & \mbox{otherwise}.
	\end{array}\right.
\end{equation}
There is a final tricky point that comes from observation of the second equation in \eqref{eq:6.4}. Since the density has been concentrated on the points $\mathbf m_j$, how can we observe it in the points $\mathbf m_i^\pm$? There is no simple reason for the following answer (see the final of this section for precise references on this), but this is it. A matrix 
\begin{equation}
\mathbf M_{ij}:=
	\left\{\begin{array}{ll} 
		\smallfrac79 & i=j, \\[1.5ex] \smallfrac19 & i=j\pm 1 \,(\mathrm{mod} N), \\[1.5ex] 0 & \mbox{otherwise},
	\end{array}\right.
\end{equation}
will play the role of the identity operator. Corresponding to these transfer functions, there are three time domain distributions 
\[
\mathcal L\{\boldsymbol{\mathcal S}_c\}=\mathbf S(\cdot/c), 
\qquad
\mathcal L\{\boldsymbol{\mathcal V}_c\}=\mathbf V(\cdot/c), 
\qquad
\mathcal L\{\boldsymbol{\mathcal J}_c\}=\mathbf J(\cdot/c). 
\]

\paragraph{The semidiscrete time domain problem.} The incident wave \eqref{eq:6.4a} is observed in the observation points to create a discrete causal function $\boldsymbol\beta:\mathbb R \to \mathbb R^N$
\begin{subequations}\label{eq:6.12}
\begin{equation}
\beta_i(t)=\sum_\pm \psi(c(t-t_{\mathrm{lag}})-\mathbf m_i^\pm\cdot\mathbf d),
\end{equation}
a causal discrete density $\boldsymbol\eta:\mathbb R \to \mathbb R^N$ is then computed by solving the convolution equations
\begin{equation}
(\boldsymbol{\mathcal V}_c* \boldsymbol\eta)(t)+\boldsymbol\beta(t)=0 \qquad \forall t,
\end{equation}
and then a potential is generated
\begin{equation}
U(t)=(\boldsymbol{\mathcal S}_c*\boldsymbol\eta)(t) \qquad \forall t,
\end{equation}
and as well as an approximation of the normal derivative, $\boldsymbol\lambda:\mathbb R\to \mathbb R^N$,
\begin{equation}
\mathbf M \boldsymbol\lambda(t)=-\smallfrac12\mathbf M \boldsymbol\eta(t)+
(\boldsymbol{\mathcal J}_c*\boldsymbol\eta)(t).
\end{equation}
\end{subequations}
It is interesting to notice than even after discretization the function $U$ is still a causal solution of the wave equation in $\mathbb R^2\setminus\Gamma$. The final fully discrete method comes from applying CQ to each of the three convolutions in \eqref{eq:6.12}. In the case of multistep CQ, the function $\boldsymbol\beta$ is sampled at equidistant times (the steps), and then three CQ processes (one convolution equation and two forward convolutions) are launched. In the case of RKCQ, the sampling is done at the stage points, the convolution equation is solved at the stage level, and finally two forward convolutions yield approximations at the step points. 

\newpage

\begin{figure}[H]
\begin{center}
\includegraphics[width=7.3cm]{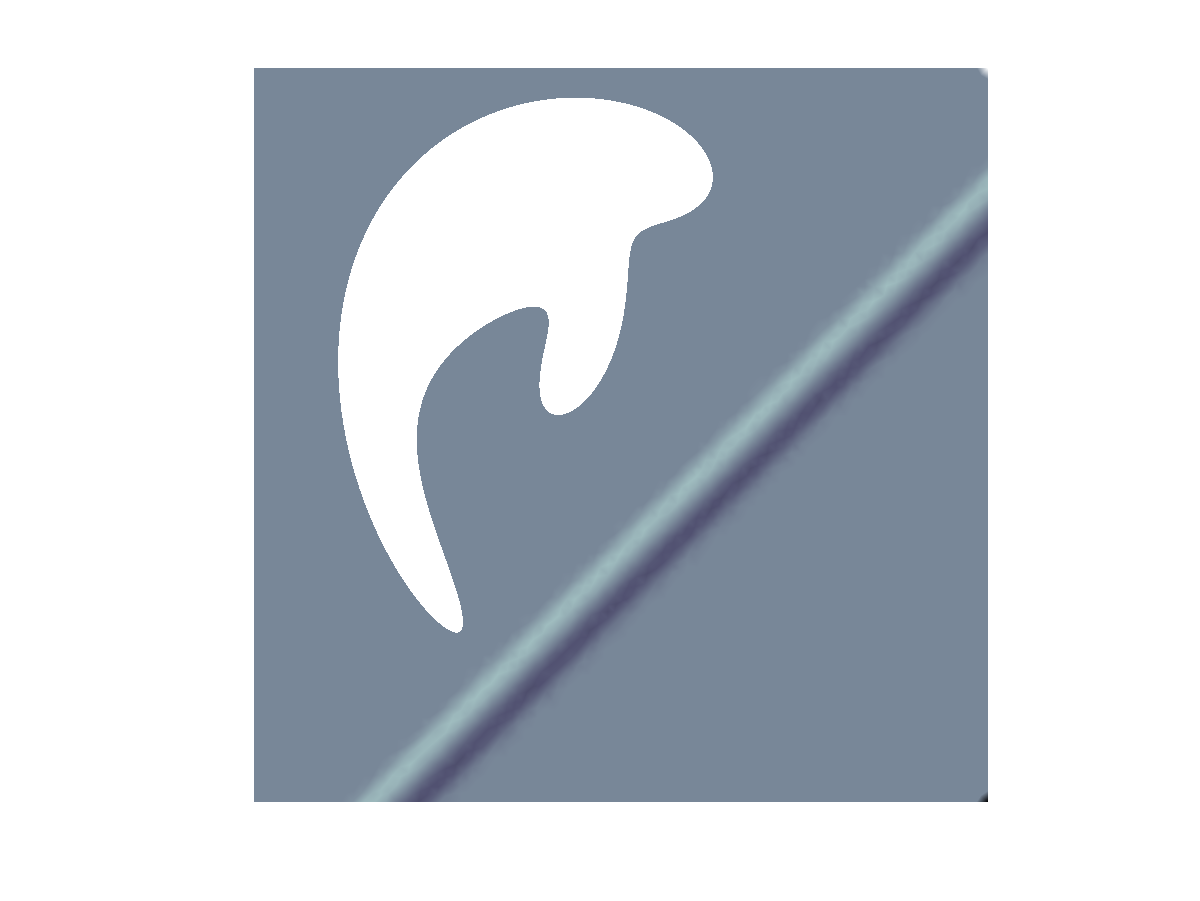} 
\includegraphics[width=7.3cm]{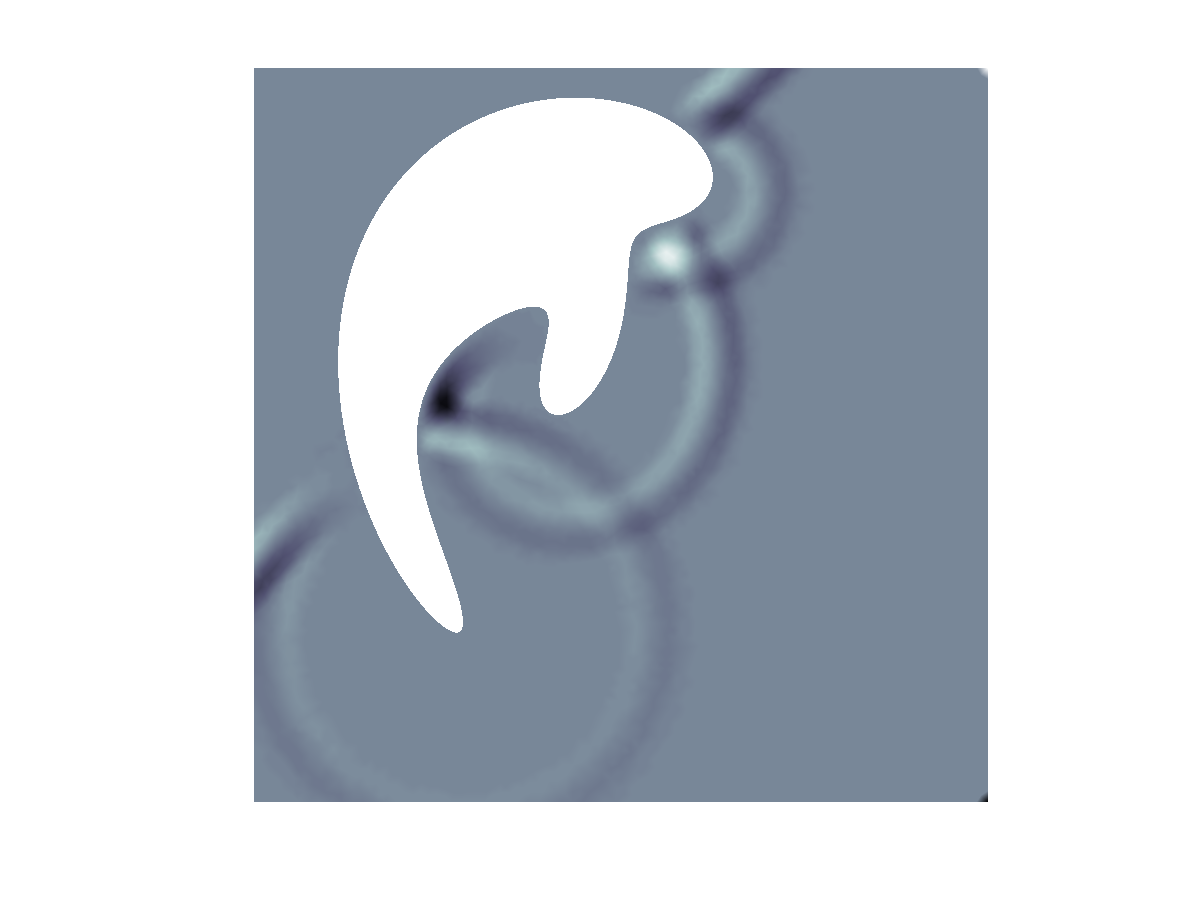} 
\includegraphics[width=7.3cm]{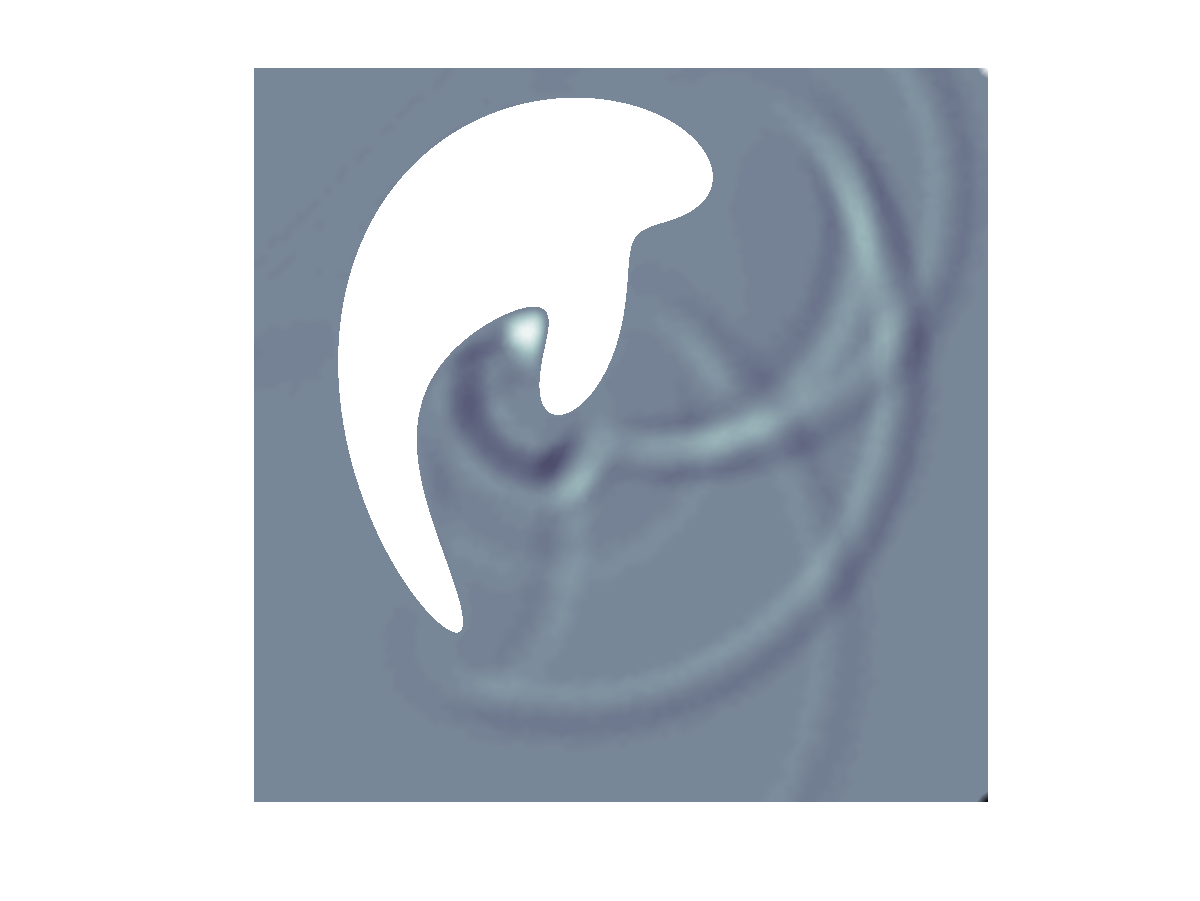} 
\includegraphics[width=7.3cm]{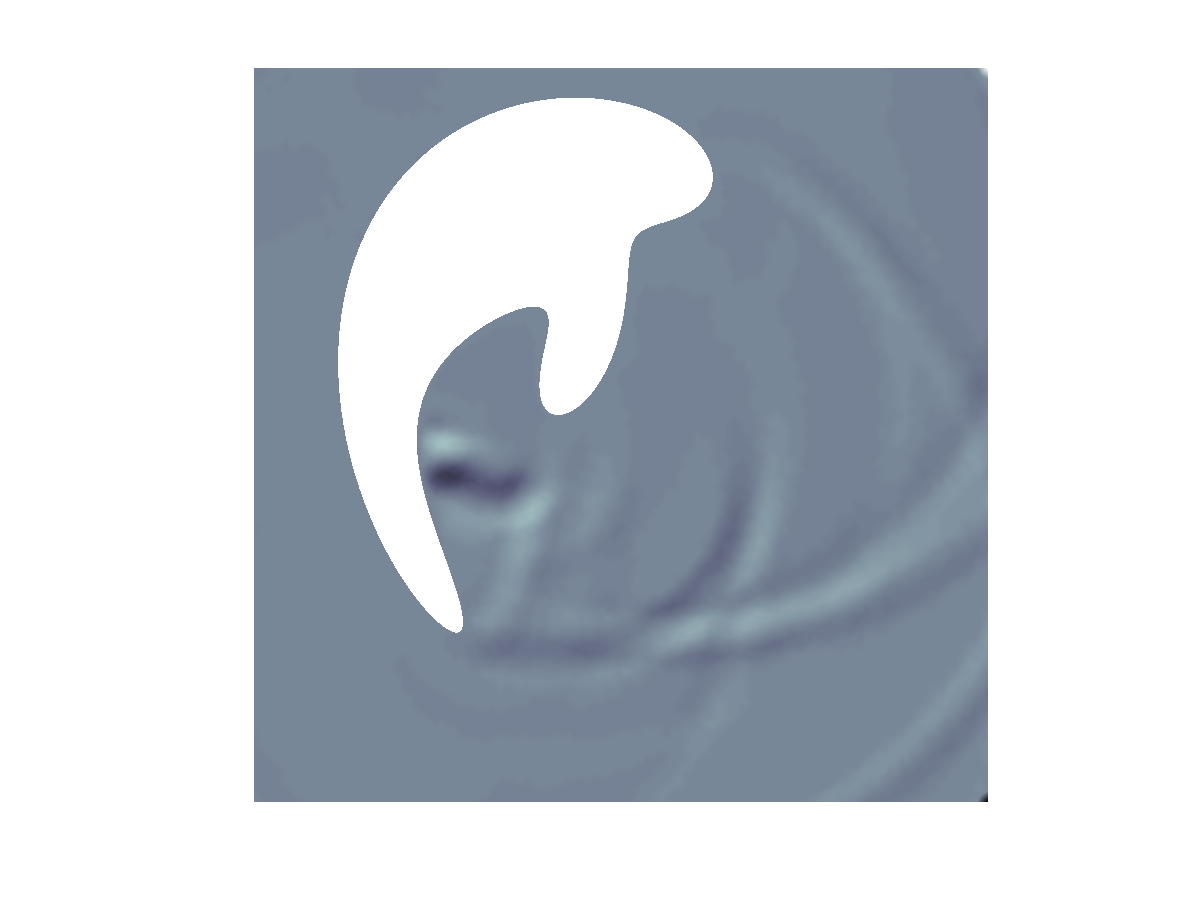} 
\includegraphics[width=7.3cm]{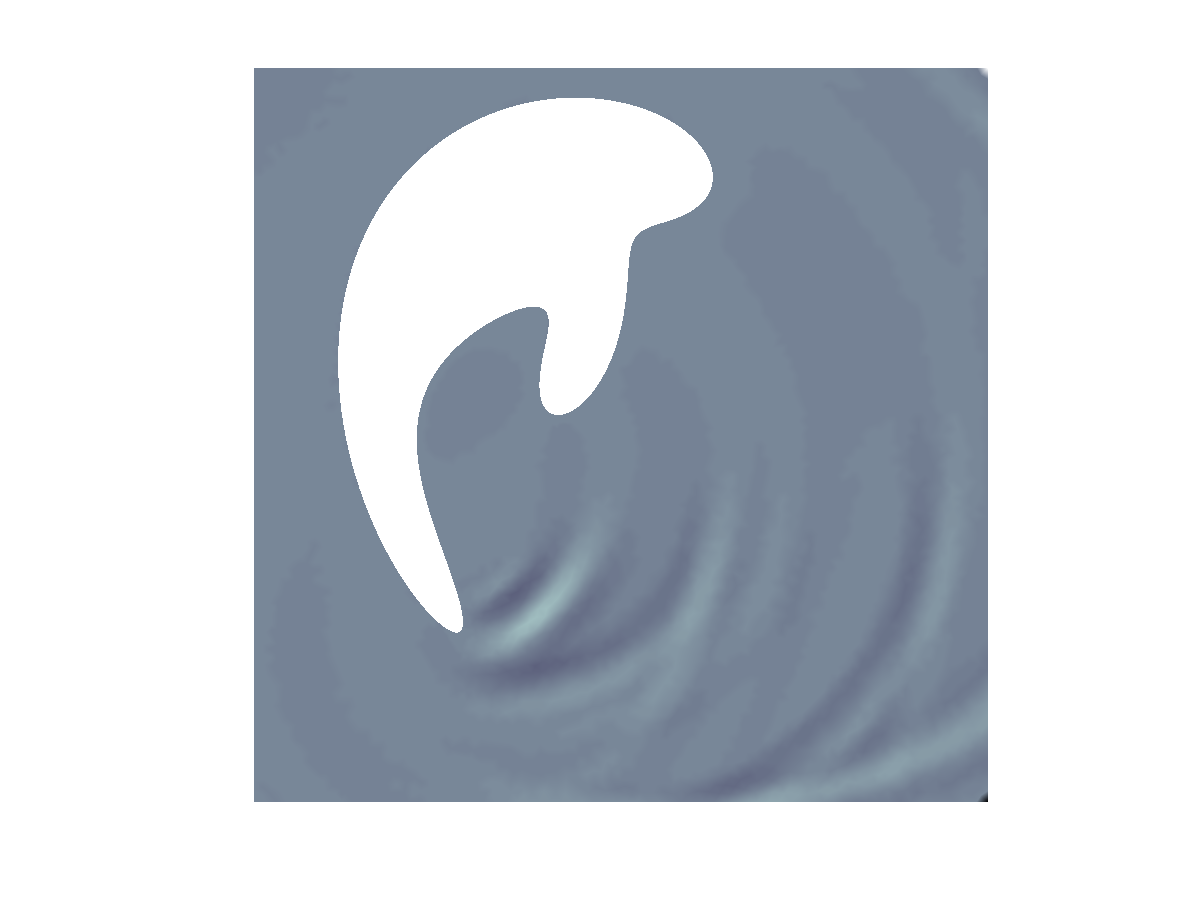} 
\includegraphics[width=7.3cm]{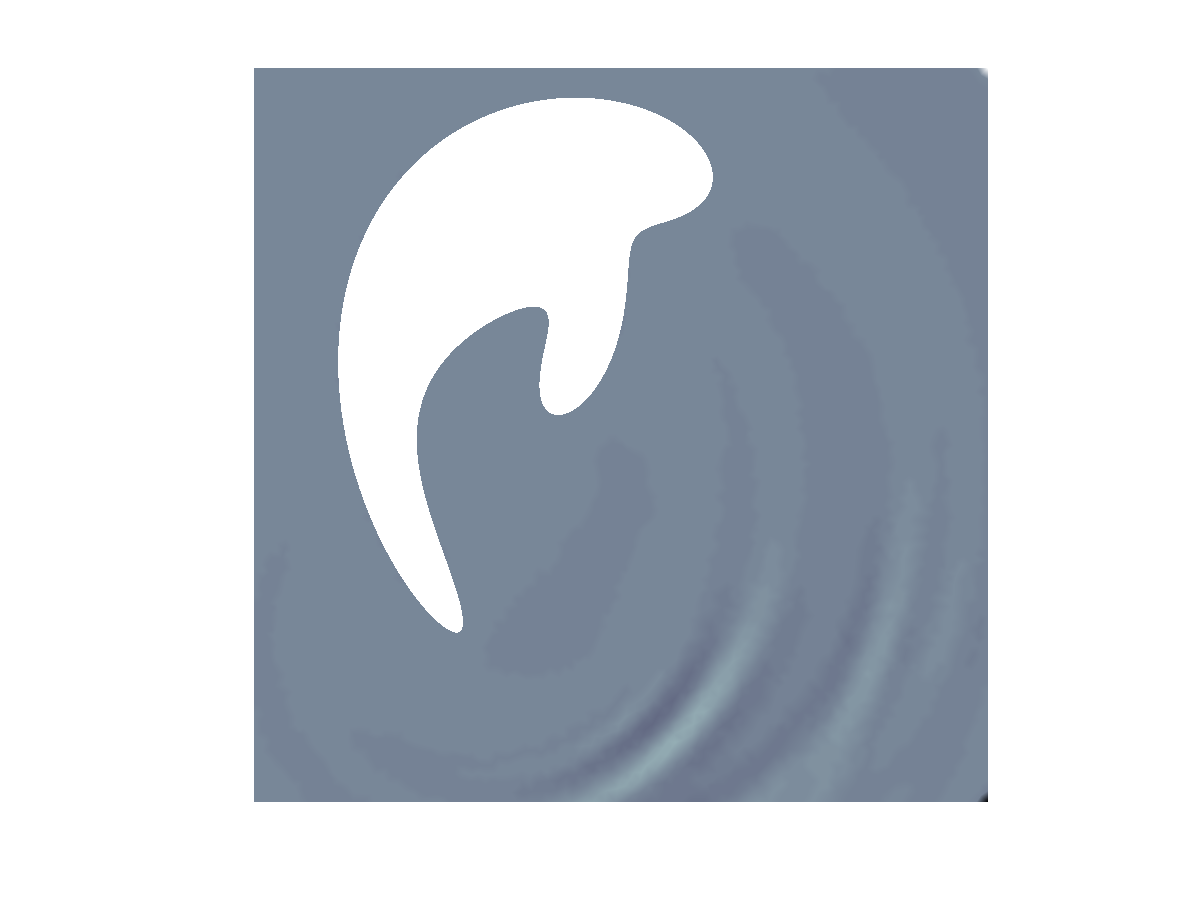} 
\end{center}
\caption{Six snapshots of a manatee-shaped sound-soft scatterer being hit by a plane wave, clearly visible in the first picture. The integral equations are solved using the method of Section \ref{sec:6} on the boundary of the scatterer. After that the potential is computed on points in the given frame. An order three Radau IIA CQ method was used in the time domain.} 
\end{figure}

\subsection*{Credits}

This section is based on the fully discrete Calder\'on Calculus developed in \cite{DoLuSa:2014}. The somewhat puzzling choices of parameters (the relative distances $\pm1/6$ to the observation grids, the matrices $\mathbf M$ and $\mathbf Q$) can be justified using careful Fourier analysis. Some intuitive explanation can be gathered from \cite{DoLuSa:2014}.


\section{The theory of convolution quadrature}\label{sec:7}

In this section we collect some convergence results for multistep and multistage CQ applied to convolutions whose symbol is defined in $\mathbb C_+$. The multistage results are taken from \cite{Lubich:1994} with a slight refinement (on the behavior of constants with respect to time) to be found in \cite{Sayas:2013}. The trapezoidal rule is not covered by that analysis but can be found in \cite{Banjai:2010}. The convergence of multistage CQ for problems relevant to the wave equation was developed in \cite{BaLu:2011} and \cite{BaLuMe:2011}. 

\subsection{Multistep CQ}

In this section we give a fast review of some results on convergence for multistep CQ for general operator valued convolutions. For ease of reference, let us recall some notation:
\[
\mathbb C_+:=\{ s\in \mathbb C\,:\, \mathrm{Re}\,s>0\}, \qquad \sigma:=\mathrm{Re}\,s, \qquad \underline\sigma:=\min\{ 1,\sigma\}.
\]

\paragraph{Hypotheses on $\delta(\zeta)$ for multistep methods.} The function $\delta:\mathcal U\to \mathbb C$ is analytic in $\mathcal U$ where
\[
\{ \zeta\in \mathbb C\,:\, |\zeta|\le 1\} \subset \mathcal U.
\]
We also require that
\[
\delta: B(0,1) \to \mathbb C_+,
\]
that is, $\mathrm{Re}\,\delta (\zeta)>0$ for all $\zeta$ such that $|\zeta|<1$. Finally, we require that there exists $q\ge 1$, $C_0>0$, and $\kappa_0>0$ such that
\[
\left| \smallfrac1\kappa \delta(e^{-\kappa})-1\right|\le C_0 \kappa^q \qquad \forall \kappa\le \kappa_0.
\]
These hypotheses are satisfied by the BE method $\delta(\zeta)=1-\zeta$ with $q=1$ and by the BDF2 method $\delta(\zeta)=\frac32-2\zeta+\frac12\zeta^2$ with $q=2$. Note that the trapezoidal rule is not covered by this analysis, and TR has to be analyzed using some different tricks.
Some work in the complex plane implies that the map $\mathbb C_+\ni s\mapsto s_\kappa:=\frac1\kappa \delta(e^{-s\kappa})$ satisfies
\begin{equation}\label{eq:7.1}
|s_\kappa| \le C_1 |s|, \qquad |s_\kappa-s|\le C_2 \kappa^q |s|^{q+1}, \qquad
\mathrm{Re}\,s_\kappa \ge C_3 \underline\sigma \qquad \forall s\in \mathbb C_+.
\end{equation}

\paragraph{Hypotheses on the transfer function.} Let now
\begin{subequations}\label{eq:7.2}
\begin{equation}
\mathrm F:\mathbb C_+\to \mathcal B(X;Y)
\end{equation} 
be analytic and satisfy
\begin{equation}
\| \mathrm F(s)\|_{X\to Y} \le C_{\mathrm F}(\sigma) |s|^\mu \qquad \forall s\in \mathbb C_+,
\qquad \mbox{with $\mu \ge 0$},
\end{equation}
where $C_{\mathrm F}:(0,\infty)\to (0,\infty)$ is non-increasing and
\begin{equation}
C_{\mathrm F}(\sigma) \le \frac{C_0}{\sigma^m}\qquad \forall \sigma\in(0,1],\quad \mbox{with $m\ge 0$}.
\end{equation}
\end{subequations}
Recall that Section \ref{sec:1.4} identified $\mathrm F$ with the Laplace transform of an operator-valued causal distribution which could be written as a distributional derivative of a continuous causal operator-valued function with polynomial growth. Thanks to \eqref{eq:7.1}, we can also identify the map $\mathbb C_+\ni s\mapsto \mathrm F(s_\kappa)$ with the Laplace transform of a causal distribution with values in $\mathcal B(X;Y)$.

\paragraph{A convergence result.} Let $\delta$ and $\mathrm F$ satisfy the above hypotheses. Let 
$g:\mathbb R\to X$ be causal and $\mathcal C^k$ with 
\[
k> \mu+q+2.
\]
Then
\begin{equation}\label{eq:7.3}
\| (f*g)(t)-(f_\kappa*g)(t)\|_Y \le C \kappa^q h(t) \int_0^t \| g^{(k)}(\tau)\|_X\mathrm d\tau,
\end{equation}
where
\[
h(t) =\left\{ \begin{array}{ll} t^{k-(\mu+q+1)}, &t\le 1, \\ t^{k-\mu+m}, & t\ge 1.\end{array}\right.  
\]
In \eqref{eq:7.3}, $f=\mathcal L^{-1}\{ \mathrm F\}$ and $f_\kappa=\mathcal L^{-1} \{ \mathrm F_\kappa\}$, with $\mathrm F_\kappa(s):=\mathrm F(s_\kappa)$. 

\subsection{Multistage CQ}

\paragraph{Order of convergence for an RK method.} Consider an RK method applied to the IVP
\[
\dot y=f(t,y), \quad 0\le t\le T, \qquad y(0)=y_0.
\]
The internal stages
\[
\boldsymbol y_{n}=y_n\boldsymbol 1 + \kappa \boldsymbol A f(t_n+\kappa\boldsymbol c,\boldsymbol y_{n})
\]
and the steps
\[
y_{n+1}=y_n+\kappa \boldsymbol b^\top f(t_n+\kappa\boldsymbol c,\boldsymbol y_{n}).
\]
create approximations
\[
\boldsymbol y_{n}\approx y(t_n+\kappa\boldsymbol c)\qquad y_n \approx y(t_n).
\]
We say that the {\bf stage order} of the RK method is $q$ when for a smooth enough solution $y$
\[
|\boldsymbol y_1-y(\kappa \boldsymbol c)|\le C h^{q+1}.
\]
We say that the {\bf classical order} of the RK method is $r$ when for a smooth enough solution $y$
\[
| y_1-y(t_1)|\le C h^{r+1}.
\]
Note that the two methods given at the end of Section \ref{sec:5.1} have respective classical orders 3 and 4, while both share stage order equal to 2.

\paragraph{Hypotheses on the RK method.} Consider an RK method and its stability function
\[
R(z)= 1+z \boldsymbol b^\top (\boldsymbol I-z\boldsymbol A)^{-1} \boldsymbol 1.
\]
We will assume that:
\begin{itemize}
\item[(a)] The matrix $\boldsymbol A$ is invertible. (This is needed right at the beginning of Section \ref{sec:5.1}, in order to give a definition to the discrete multistage differentiation operator.)
\item[(b)] (A-stability) For all $z$ such that $\mathrm{Re}\, z \le 0$, the matrix $\boldsymbol I-z\boldsymbol A$ is invertible and
\[
|R(z)| \le 1.
\]
As mentioned in Section \ref{sec:5.3}, hypotheses (a) and (b) imply that $\sigma(\boldsymbol A)\subset \mathbb C_+$.
\item[(c)] $R(\infty)=0$, that is
\[
\boldsymbol b^\top \boldsymbol A^{-1} \boldsymbol 1=1.
\]
(Note that the quantity $\mu=1-\boldsymbol b^\top \boldsymbol A^{-1} \boldsymbol 1=R(\infty)$ had appeared in \eqref{eq:5.200} at the time when we wanted to compute steps for RKCQ.)
\item[(d)] 
\[
|R(\imath \omega)| < 1\qquad \forall \omega\in \mathbb R\setminus\{0\}.
\]
\end{itemize}

\paragraph{A word on the RKCQ output.} While in principle the RKCQ produces a sequence of vectors $Y^p\ni \boldsymbol y_n\approx y(t_n+\kappa\boldsymbol c)$, with $y=f*g$, and a sequence of `scalars' $Y\ni y_n \approx y(t_n)$, like in the multistage case, there is a formula that extends these values to continuous times. We saw in Section \ref{sec:5.4} that RKCQ could be understood as outputting
\begin{equation}\label{eq:7.4}
\boldsymbol y_\kappa (t)=\sum_{m=0}^\infty W_m^{\mathrm F}(\kappa) g(t-t_m+\kappa\boldsymbol c),
\end{equation}
so that
\[
\boldsymbol y_\kappa(t_n)=\boldsymbol y_n.
\]
The time steps can be produced using a simple postprocessing of the stages
\[
y_{n+1}=(\boldsymbol b^\top \boldsymbol A^{-1}\otimes I) \boldsymbol y_n
\]
(recall \eqref{eq:5.201} and \eqref{eq:5.202} and notice that we are assuming that $\mu=0$). At the continuous level, this corresponds to the values at time $t_{n+1}$ of
\begin{equation}\label{eq:7.5}
y_\kappa = (\boldsymbol b^\top \boldsymbol A^{-1}\otimes I) \boldsymbol y_\kappa (\cdot-\kappa).
\end{equation}

\paragraph{Hypotheses on the transfer function.} The hypotheses on $\mathrm F$ are slightly different than those given for multistep CQ. We assume that
\begin{subequations}\label{eq:7.6}
\begin{equation}
\mathrm F:\mathbb C_+\to \mathcal B(X;Y)
\end{equation} 
is analytic and satisfies
\begin{equation}
\| \mathrm F(s)\|_{X\to Y} \le C_{\mathrm F}(\sigma) \frac{|s|^\mu}{\sigma^\nu} \qquad \forall s\in \mathbb C_+,
\qquad \mbox{with $\mu \ge 0$, $\nu\ge 0$},
\end{equation}
where $C_{\mathrm F}:(0,\infty)\to (0,\infty)$ is non-increasing and
\begin{equation}
C_{\mathrm F}(\sigma) \le \frac{C_0}{\sigma^m}\qquad \forall \sigma\in(0,1],\quad \mbox{with $m\ge 0$}.
\end{equation} 
\end{subequations}
Note that we have factored out $\sigma^\nu$ from $C_{\mathrm F}$ so that a bound with a power of $\sigma$ in the denominator is also valid as $\sigma\to \infty$. 

\paragraph{A convergence result.} Let $\mathrm F$ satisfy hypotheses \eqref{eq:7.2} and the RK method satisfy the previous hypotheses with stage order $q$ and classical order $r\ge q$. Assume that $g:\mathbb R\to X$ is causal and $\mathcal C^k$ with
\[
k> \mu+r+2. 
\]
If $y=f*g$ and $y_\kappa$ is the RKCQ approximation of $y$ using \eqref{eq:7.4}-\eqref{eq:7.5}, then
\[
\| y(t)-y_\kappa(t)\|_Y \le C \kappa^{\min\{r,q+1-\mu+\nu\}} h(t) \int_0^{t} \| g^{(k)}(\tau)\|_X\mathrm d\tau.
\]
Here $h$ is an increasing function of time, whose behavior is not entirely well understood, although it is unlikely that they behave worse than polynomially in time. For the particular case of operators satisfying
\[
\|\mathrm F(s)\| \le  C e^{-c\,\sigma}|s|^\mu \qquad \forall s \in \mathbb C_+, \qquad \mu\ge 0,\quad c>0,
\]
the full classical convergence order of the method is attained for smooth enough functions, since we can choose any arbitrarily large $\nu$ in the hypotheses for $\mathrm F$.

\bibliographystyle{abbrv}
\bibliography{referencesEHF}

\end{document}